\documentclass[twoside]{article}
\usepackage[margin=1in]{geometry}
\usepackage{amsmath,amsthm,amssymb,bm}
\usepackage{amsfonts}
\usepackage{mathrsfs}
\usepackage{mathtools}
\usepackage{tikz}
\newtheorem{remark}{Remark}
\usepackage{graphicx}
\usepackage[font=small,labelfont=bf]{caption}
\usepackage{subfig}
\usepackage{float}
\usepackage{bbold}
\usepackage{adjustbox}

\usepackage{cite}
\usepackage{nicefrac}
\usepackage{euscript,mathrsfs}
\usepackage{color}
\usepackage{multirow}
\usepackage{appendix}

\usepackage{fancyhdr}
\allowdisplaybreaks[1]

\numberwithin{equation}{section}
\numberwithin{figure}{section}
\numberwithin{table}{section}

\newcommand{\tA}{\tilde{A}}
\newcommand{\ta}{\tilde{a}}
\newcommand{\tw}{\tilde{w}}

\newcommand{\U}[1]{{\bm{#1}}}
\def\be{\begin{equation}}
\def\ee{\end{equation}}
\newcommand{\Q}{\mathbf{Q}}
\newcommand{\xx}{\bm{x}}
\newcommand{\uu}{\bm{u}}
\newcommand{\dtdx}{\frac{\Delta t}{\Delta x}}
\newcommand{\dtdy}{\frac{\Delta t}{\Delta y}}

\newcommand{\ihp}{i+\frac{1}{2}}
\newcommand{\jhp}{j+\frac{1}{2}}

\newcommand{\ihm}{i-\frac{1}{2}}
\newcommand{\jhm}{j-\frac{1}{2}}

\newcommand{\ip}{i+1}
\newcommand{\jp}{j+1}

\newcommand{\dx}{{\Delta x}}
\newcommand{\dy}{{\Delta y}}

\theoremstyle{plain}			

\newenvironment{acknowledgments}{{\flushleft \bf Acknowledgment:}}{}

\title{On a hybrid continuum-kinetic  model for complex fluids}

\author
{A. Chertock\thanks{Department of Mathematics, North Carolina State University, USA; {\tt chertock@math.ncsu.edu}}
\and P. Degond\thanks{Institut de Mathématiques de Toulouse; UMR5219; Université de Toulouse; CNRS;
UPS; F-31062 Toulouse Cedex 9, France; {\tt pierre.degond@math.univ-toulouse.fr}}
\and G. Dimarco\thanks{Department of Mathematics and Computer Science $\&$ Center for Modeling, Computing and Statistics of University of Ferrara, Italy; {\tt giacomo.dimarco@unife.it}}
\and M. Luk\'a\v{c}ov\'a-Medvid'ov\'a\thanks{Institute of Mathematics, University of Mainz, Germany;
{\tt lukacova@uni-mainz.de}}
\and A. Ruhi\thanks{Institute of Mathematics, University of Mainz, Germany; {\tt ankruh@gmail.com}}}

\date{}

\begin{document}

\maketitle

\begin{abstract}
In the present work, we first introduce a general framework for modelling complex multiscale fluids and then focus on the derivation and
analysis of a new hybrid continuum-kinetic  model. In particular, we
combine conservation of mass and momentum for an isentropic macroscopic model with a kinetic representation of the microscopic behavior.
After introducing a small scale of interest, we compute the complex stress tensor by means of the Irving–Kirkwood formula. The latter requires an expansion of the kinetic distribution
around an equilibrium state and a successive homogenization over the fast in time and small in space scale dynamics. For a new  hybrid continuum-kinetic  model the results of linear stability analysis indicate a conditional stability in the relevant low speed regimes
and linear instability for high speed regimes for higher modes. Extensive numerical experiments confirm that the proposed multiscale
model can reflect new phenomena of complex fluids not being present in standard Newtonian fluids. Consequently, the proposed general technique can be successfully
used to derive new interesting systems combining the macro and micro structure of a given physical problem.
\end{abstract}

\medskip

\medskip


\section{Introduction}\label{sec1}
Many important fluid flow problems are entirely multiscale: microscopic processes strongly influence macroscopic behavior of the fluid and need to be taken into account in order to accurately describe fluid dynamics.
Typical examples are granular \cite{granular} and high-speed rarefied flows \cite{cercignani69}, the plastic deformation in materials
\cite{plastic}, the viscoelastic \cite{visco} and biological type of fluids \cite{bio}. For this reason, in the last few decades, there
has been a huge interest in both modeling and numerical simulations of problems associated with multilevel physical models, which are able
to incorporate multiscale effects in different ways.

For the Newtonian fluids there have been many rigorous theoretical
studies of hydrodynamic limits and the relationship between microscopic molecular dynamics and/or mesoscopic kinetic models of the Boltzmann
type with macroscopic models such as the compressible Euler or Navier-Stokes equations, see, e.g.,
\cite{bardos93, dellar08, grad64, golse09, levermore10} and the references therein.
 On the other hand, for complex fluids
theoretical understanding is certainly less developed, and more research is needed.
More precisely, compared to standard fluids, the challenge associated with complex fluids lies in an accurate determination of
rheological relations that are typically obtained from physical or computational experiments. Consequently, in many
situations, e.g., soft matters or colloid-polymer mixtures, their full analytical description is not available.

In order to take small scale effects  into account we can apply either mesoscopic kinetic models or directly
microscopic models, such as molecular dynamics or dissipative particle dynamics, to reconstruct time evolution of
macroscopic quantities.
However, as is well known, an obvious drawback of meso- and microscopic descriptions, despite their higher accuracy, is sometimes
prohibitively high computational costs, which limit their direct application in many practical situations. To overcome this disadvantage
and make large-scale simulations possible the so-called {\em hybrid multiscale
methods} have been developed in the literature. The latter combine the advantages of both descriptions: the accuracy
of microscopic models to be able to consider complex rheological relationships with the efficiency of macroscopic models based on using
classical conservation laws.

A prototype of hybrid multiscale methods is
the {\em heterogeneous multiscale method} proposed by E, Enquist et al.~\cite{Weinan2007, weinan11, E2003, E2004, E_2005, Ren_2007}, see also \cite{yelash21, yelash17,  Ren_2005,  yelash22, Yasuda_2010} for its application for complex polymeric fluids. In this context let us also refer to
{triple-decker atomistic-mesoscopic-continuum method} \cite{Fedosov_2009}, {the seamless multiscale methods} \cite{E2007,e_2007},
{the equation-free multiscale methods} \cite{Kevrekidis_2003,Kevrekidis_2009} or {the internal-flow multiscale method} \cite{lockerby2, lockerby1}. In \cite{koumoutsakos} a overview of multiscale flow simulations using particles is presented.
For classical gas dynamics, similar ideas were employed, for instance, in
\cite{DP,DDP,Crestetto,Hadji}.

However, such general methods can have several disadvantages: the precise definition of the area in which the
microscopic model should be applied is problem-dependent and needs to be specified. If the microscopic description is applied
(almost) everywhere, the computational cost may become even larger than the cost required to numerically solve a micro- or
mesoscale model. Clearly,  bridging the large range of dynamically coupled scales is a fundamental challenge that was and still is
a driving force in the development of new mathematical algorithms.

The goal of the present paper is to derive and analyse a new hybrid continuum-kinetic model for complex fluids.
We start first by introducing a general methodological approach, which allows us to derive a class of new multiscale models. The basic idea is based on two simple
ingredients. First, we assume that the problems we are studying have
fast and slow scale dynamics. In addition, we assume a scale separation between the microscopic (fast) and macroscopic (slow) dynamics and thus, the phenomena may be imagined to act at different domain scales.
Second, we assume that the effects of the fast scale dynamics on microscopic scale can be captured
at the macroscopic level, at least approximately, by homogenization of the microscopic properties of the fluid over a finite size domain.
In this work, we consider a prototype situation of isentropic flows governed by
the conservation of mass and momentum at the continuum level. In order to model rheology of a complex fluid,
the non-Newtonian stress tensor is obtained by an upscaling homogenizing procedure using
the kinetic relaxation type equations.  We will study properties of the derived hybrid multiscale model
using a linear stability analysis and successively conduct several numerical
experiments to illustrate the accuracy and efficiency of the proposed hybrid model.

The rest of the paper is organized in the following way. In Section~\ref{sec2}, we derive the hybrid continuum-kinetic model
including the non-Newtonian stress tensor. The latter is obtained by the Irving-Kirkwood formula which represents upscale  microscopic effects.
In Section~\ref{sec3}, we perform a linear stability
analysis which shows that for the regimes of interest, the low speed flows, the model is linearly stable. In Section~\ref{sec4}, several
numerical examples are conducted that demonstrate the validity of the proposed model for different prototype situations arising
in complex flows. Finally, Section~\ref{sec5} is dedicated to the discussion of the obtained results and future developments.

\section{The hybrid model approach}
\label{sec2}
In this section, we first detail the general framework which will be used to derive a class of hybrid multiscale models. We then focus on
the derivation and analysis of a prototype case, a hybrid continuum-kinetic model for complex fluids. This is realized by
homogenization over a fixed size cell of the microscopic domains composed by the fluid  molecules. The latter process leads to the
Irving-Kirkwood formula for complex macroscopic stress tensor.

\begin{figure}[!ht]
\centering
\includegraphics[trim=300 350 200 200,clip=true,width=1\textwidth]{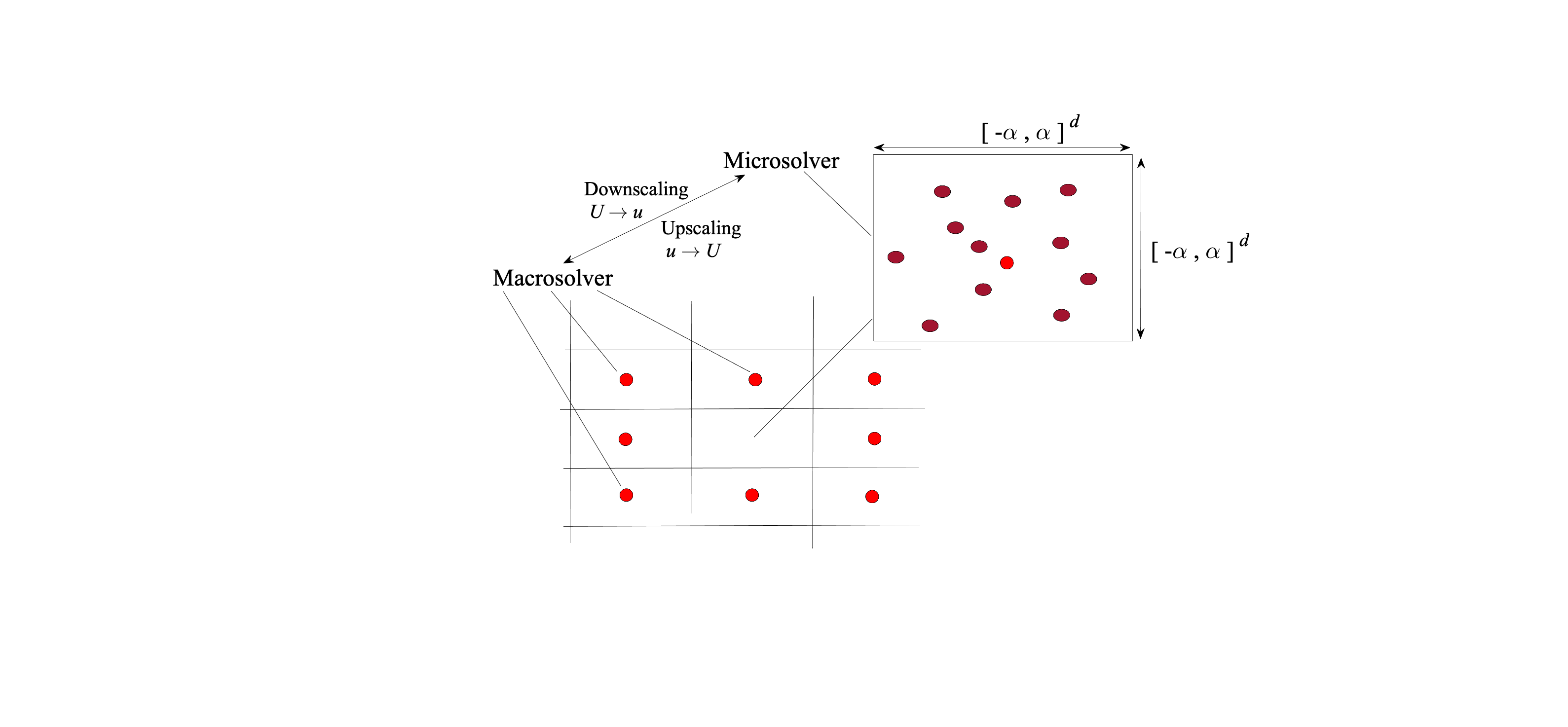}
\caption{\sf Sketch of the hybrid multiscale approach.}
\label{fig:model}
\end{figure}

We start by considering the following general setting. We assume that a microscopic process used to describe time evolution of the
state of the system of interest is known. For instance, molecular dynamics or kinetic mesoscopic equations are able to provide such
information accurately enough. We also assume that we have at our disposal a macroscopic model in which
a missing information will be provided by means of a microscopic model.
These two models, micro/meso and macro,
can be related through a reconstruction (upscaling) and a projection (downscaling) operator which permits to commute from one system to the
other. The upscaling operator averages the micro properties of the fluid up to a coarser description, while the downscaling operator uses
the coarser description to obtain the unknown variables at the microscopic level. A classical example of the above operators
are downscaling/upscaling operators between kinetic and macroscopic description.  In this case, the upscaling procedure is obtained
by the integration of the distribution function multiplied by the so-called collision invariants over the velocity space while the
downscaling operator is recovered from the knowledge of the macroscopic variable defining a so-called equilibrium distribution. In such a
setting, our aim is to be able to give a description of the state of the system by working on a given macroscopic grid, defined a-priori,
and through the use of a macroscopic model bringing some information from the microscopic/mesoscopic dynamics. Schematically, the proposed
method works as depicted in Figure~\ref{fig:model}. We upscale the microsolver information, which is determined by  time evolution
of the unknown $\bm{u}$, through a homogenization on a box of fixed size $[-\alpha,\alpha]^d$ with $d$ the spatial dimension and $\alpha$ the
characteristic length of the microscopic variation. The size of the box depends on the problem under consideration. This
micro information is successively used into the macroscopic solver, with unknown $\bm{U}$, over the macroscopic grid to update the solution of the complex flow.

In what follows, we describe the details of such procedure in the case, where the macroscopic model is represented by the isentropic fluid
equations for complex fluids while the microscopic model is the BGK kinetic equation \cite{bhatnagar54}.
From the depicted scenario, it is possible to imagine
alternative types of combination of micro and macro dynamics which will be discussed in future investigations.

\subsection{A prototype case and derivation of the model}\label{sec2.1}
In this section, we consider a fluid which, by hypothesis, can be described at the macroscopic level by the compressible isentropic fluid
equations. In this case, the system is governed by two equations describing the conservation of mass and momentum.
This system reads
\begin{equation}
\label{eq2}
\begin{split}
&\partial_t \rho +\nabla_{\bm{x}} \cdot (\rho \bm{u})= 0,\\
&\partial_t (\rho \bm{u}) + \nabla_{\bm{x}} \cdot (\rho \bm{u} \otimes \bm{u}) =  \nabla_{\bm{x}} \cdot \mathbb{T},
\end{split}
\end{equation}
where $\rho$ and $\bm{u}$ stand for the fluid density and velocity, respectively.  The so-called Cauchy stress tensor $\mathbb{T}$
describes specific rheological fluid properties. All unknowns are functions of space $\bm{x}\in\mathbb{R}^d$ and time $t>0$.
For inviscid fluids,
$\mathbb{T}=-p \mathbb{I}$, where $\mathbb{I}$ is the identity tensor. This expression leads to {\em the Euler equations}. Considering isentropic fluids, pressure $p = p(\rho)$ is a given function of density $\rho.$

If viscous effects are taken into account, the Cauchy stress tensor reads
\begin{equation*}
\label{eq3}
\mathbb{T} = -p \mathbb{I} + \mathbb{S},
\end{equation*}
where $\mathbb{S}$ stands for the viscous stress.  For Newtonian fluids, the latter
is given by the Newton rheological rule (NRR)
\begin{equation*}
\label{NRR}
\mathbb{S}=\mu\left(\nabla_{\bm{x}}\bm{u}+\nabla_{\bm{x}}\bm{u}^T-\frac{2}{d}\nabla_{\bm{x}}\cdot\bm{u}\mathbb{I}\right)+
\lambda\nabla_{\bm{x}}\cdot\bm{u}\mathbb{I},
\end{equation*}
with $\mu > 0$ and $\lambda \geq 0$ being the constant shear and bulk viscosity coefficients, respectively.
This relation leads to the compressible {\em Navier-Stokes equations}.
 In complex fluids, however,
the rheological relation for the Cauchy stress tensor is more general and typically obtained by computational or physical experiments.
Here, we propose instead that the Cauchy stress tensor
takes into account microscopic effects in a homogenized way. This new model is derived from considering the kinetic equations underpinning
fluid models of the type \eqref{eq2} as a microscopic model.

More specifically, we make the hypothesis that the so-called BGK (Bhatnagar-Gross-Krook) equation is a suitable microscopic model which can be upscaled
to yield a missing information on a complex Cauchy stress tensor at the macroscopic level. The considered BGK equation reads
\cite{bhatnagar54}
\begin{equation}
\label{cell_problem}
\partial_t f + \bm{v}\cdot\nabla_{\bm{x}}f = Q(f),
\end{equation}
where $f$ is the probability density function of fluid molecules at position $\bm{x}$ having velocity $\bm{v}\in\mathbb{R}^d$ at time $t$ and $Q(f)$ is a collision operator (modelling the molecular interactions) to be described later on. This hypothesis can be relaxed giving
rise to similar models with, however, different coefficients.

We assume now as depicted in Figure \ref{fig:model} that $f$ has slow variations at the domain scale $D$ and fast variations at scale
$\alpha\ll D$. Given that hypothesis, we consequently set  $f = \tilde{f}(\bm{x},\frac{\bm{x}}{\alpha},\bm{v},t)$ where
$\tilde{f}(\bm{x},\bm{y},\bm{v},t)$ is 2-periodic with respect to $\bm{y}$, with unit cell $[-1,1]^{d}$. We also suppose that the collision
operator has magnitude $1/\alpha$, i.e. that the microscopic spatial effect are balanced by the collision dynamics among molecules.
Substituting this representation into \eqref{cell_problem}, yields the following microscopic model
\begin{equation}
\label{cell_problem_2}
\partial_t \tilde{f} + \bm{v}\cdot\nabla_{\bm{x}} \tilde{f} + \frac{1}{\alpha} \bm{v}\cdot\nabla_{\bm{y}} \tilde{f}
= \frac{1}{\alpha} Q(\tilde f).
\end{equation}
Now, we introduce two new functions, one for the density $ \tilde{\rho} = \tilde{\rho}(\bm{x},\bm{y},t)$ and one for the vector velocity
$ \tilde{\bm{u}} = \tilde{\bm{u}}(\bm{x},\bm{y},t)$ to be precisely defined later on. $Q(\tilde f)$ is
the BGK operator, i.e. a relaxation operator towards a given Maxwellian distribution $M_{(\tilde \rho, \tilde{\bm{u}})}$ of parameter
$\tilde{\rho}$ and $\tilde{\bm{u}}$. We observe at this stage that the way in which $ \tilde{\rho}$ and $ \tilde{\bm{u}}$ will be defined
will permit to close the resulting system of macroscopic equations. Denoting by $T$ the (constant) temperature and $R$ the universal gas
constant, $M_{(\tilde \rho, \tilde{\bm{u}})}$ is given by
\begin{equation}
M_{(\tilde \rho, \tilde{\bm{u}})}(\bm{v})
= \frac{\tilde \rho}{(2\pi R T)^{d/2}}\exp{\left(-\frac{|\bm{v}- \tilde{\bm{u}}|^2}{2R T}\right)},
\label{max}
\end{equation}
and the BGK operator is given by
\begin{equation}\label{BGK}
 Q(\tilde f) = -\frac{1}{\tilde \tau} ( \tilde f - M_{(\tilde \rho, \tilde{\bm{u}})} ),
\end{equation}
where $\tilde \tau$ is a relaxation time. Finally, \eqref{cell_problem_2} takes the form
\begin{equation}
\label{cell_problem_3}
\alpha \big( \partial_t \tilde{f} + \bm{v}\cdot\nabla_{\bm{x}} \tilde{f} \big) +  \bm{v}\cdot\nabla_{\bm{y}} \tilde{f}
= -\frac{1}{\tilde \tau} ( \tilde f - M_{(\tilde \rho, \tilde{\bm{u}})} ).
\end{equation}
The microscopic model \eqref{cell_problem_3} can be used to provide a missing detailed information for the macroscopic model \eqref{eq2}. To
this end,  the Cauchy stress tensor $\mathbb{T}$ can be reconstructed from $\tilde f$ by an upscaling homogenization
procedure through the Irving-Kirkwood (or Kramer) formula \cite{Irving_1950}:
\begin{equation}
\mathbb{T} (\bm{x}, t) = - \frac{1}{2^d} \int\displaylimits_{\left[ -1,1 \right]^{d}} \int_{\mathbb{R}^d} (\bm{v}-\tilde{\bm{u}}(\bm{x},
\bm{y},t))\otimes (\bm{v}-\tilde{\bm{u}}(\bm{x},\bm{y},t)) \tilde f(\bm{x},\bm{y},\bm{v},t) d\bm{v} d\bm{y}
\label{eq:expressCST}
\end{equation}
with $\tilde{f}$ being the solution of \eqref{cell_problem_3}. In other words, we measure the microscopic effects by homogenization of the kinetic
model over the microscopic box of size $[-\alpha,\alpha]^d$. Problem \eqref{cell_problem_3} can be viewed as a cell problem for
$\bm{y}$ in the unit cell $\left[ -1,1 \right]^{d}$ and $\bm{v} \in {\mathbb R}^d$. The size of this box is left as a free parameter, which
may depend on the system under consideration. An interesting setting, which will be numerically explored in Section \ref{sec4}, consists of
considering the size of the box into relation with the macroscopic description of the flow at the numerical level. This can be done, for
instance, by fixing a ratio between the box where the microscopic effects are measured and the size of the mesh employed at the macroscopic
level.

In order to have a well posed problem, we finally supplement \eqref{cell_problem_3} with the periodic boundary conditions for the
probability density function:
$$
\tilde f(\bm{x},\bm{y} + \bm{n},\bm{v},t) = \tilde f(\bm{x},\bm{y},\bm{v},t), \quad  \forall \bm{n} \in (2\mathbb{Z})^{d}.
$$
Our aim now is to find an approximate solution of the cell problem by means of the Chapman-Enskog expansion \cite{cercignani69} for small
relaxation times $\tilde \tau$. This will allow us to upscale the microscopic effects at the macroscopic level without resorting to the
numerical resolution of the microscopic model which would lead to an expensive computation. This step is performed in the next section.

\subsection{Chapman-Enskog expansion for the prototype hybrid model}\label{sec:hybrid}
In this part, we present a perturbation analysis of the distribution
function $\tilde f$ over the box $[-1,1]^d.$   We assume in the rest of the paper that $\tilde{\tau} = O(\varepsilon),$ where $\varepsilon \ll 1$ is the
Knudsen number, i.e. the ratio of the relaxation parameter $\tilde{\tau}$ to a characteristic hydrodynamic temporal scale. The latter
represents the macroscopic scale of interest in our analysis. More precisely, we write $\tilde{\tau} = \varepsilon \hat{\tau}$, where
$\hat{\tau} = O(1)$ when $\varepsilon \rightarrow 0$ in the following. We proceed by introducing the Chapman-Enskog expansion \cite{cercignani69}
that will be truncated at the first-order terms
\begin{equation}
\label{CE_Expansion}
\tilde f^\varepsilon = \tilde f^{(0)} + \varepsilon \tilde f^{(1)} + {\mathcal O}(\varepsilon^2),
\end{equation}
where $\tilde f^\varepsilon$ is the solution of \eqref{cell_problem_3} when $\tilde \tau$ is replaced by $\varepsilon \hat{\tau}$.
The so-called cell problem \eqref{cell_problem_2} is then rewritten, using the Einstein summation convention under the above scaling as
\begin{equation}
\label{cell_problem_nondimensional}
\Big( \alpha \big( \partial_t \tilde f^\varepsilon + v_i \partial_{x_i} \tilde f^\varepsilon \big)
+  v_i \partial_{y_i} \tilde f^\varepsilon \Big) (\bm{x},\bm{y},\bm{v},t)
= -\frac{1}{\varepsilon\hat{\tau}} \Big[\tilde f^{\varepsilon}(\bm{x},\bm{y},\bm{v},t)
- M_{(\tilde \rho^{\varepsilon}(\bm{x},\bm{y},t),\tilde{\bm{u}}^{\varepsilon}(\bm{x},\bm{y},t))}(\bm{v}) \Big].
\end{equation}
Substituting \eqref{CE_Expansion} into \eqref{cell_problem_nondimensional} and equating the same powers of $\varepsilon$ yields
\begin{equation*}
\label{CE_HM}
\tilde{f}^{(0)} =  M_{(\tilde{\rho}, \tilde{\bm{u}})},\quad
\tilde{f}^{(1)} = -\hat{\tau} \Big( \alpha \big( \partial_t \tilde f^{(0)} + v_i \partial_{x_i} \tilde f^{(0)} \big) +  v_i \partial_{y_i}
\tilde f^{(0)} \Big).
\end{equation*}
We proceed now, using the first order expansion \eqref{cell_problem_nondimensional}, to the computation of the Cauchy stress tensor
appearing in equation \eqref{eq2} in the above described setting. To this end, we denote
\begin{eqnarray}
\mathcal{T}_{ij}&=&\int_{\mathbb{R}^d} (v_j-\tilde{u}_j)(v_i-\tilde{u}_i)\tilde{f}^{(0)} d\bm{v}  +
\varepsilon\int_{\mathbb{R}^d} (v_j-\tilde{u}_j)(v_i-\tilde{u}_i)\tilde{f}^{(1)}d\bm{v} + {\mathcal O}(\varepsilon^2)
\label{Equation_for_P}  \\
&=:& \mathcal{T}_{ij}^{(0)} + \mathcal{T}_{ij}^{(1)} + {\mathcal O}(\varepsilon^2) ,
\nonumber
\end{eqnarray}
the contribution to the stress tensor coming from the Chapmann-Enskog expansion and we consider in the following computations each term
separately.

First, using the equation of state for perfect gases $\tilde{p}=RT\tilde{\rho}$, we obtain that the leading term in \eqref{Equation_for_P}
can be written as
\begin{subequations}
\begin{equation}
\begin{aligned}
\mathcal{T}^{(0)}_{ij} &=\int_{\mathbb{R}^d} (v_j-\tilde{u}_j)(v_i-\tilde{u}_i)\tilde{f}^{(0)} d\bm{v}\\
&=\int_{\mathbb{R}^d} (v_j-\tilde{u}_j)(v_i-\tilde{u}_i)  \frac{\tilde{\rho}}{(2\pi R T)^{d/2}}\exp{\left(-\frac{|\bm{v}-\tilde{\bm{u}}|^2}
{2R T}\right)} d\bm{v}
=\tilde{\rho} R T \delta_{ij} = \tilde{p}\delta_{ij}
\end{aligned}
\label{T0}
\end{equation}
or in the matrix form as
\begin{equation}
\mathcal{T}^{(0)}=\tilde{p}\mathbb{I}.
\label{T0a}
\end{equation}
\end{subequations}
The second term in \eqref{Equation_for_P} reads:
\begin{eqnarray}
\mathcal{T}^{(1)}_{ij} &=& \varepsilon\int_{\mathbb{R}^d} (v_j-\tilde{u}_j)(v_i-\tilde{u}_i)\tilde{f}^{(1)} d\bm{v} \nonumber \\
&=&-\hat{\tau}\varepsilon \int_{\mathbb{R}^d} (v_j-\tilde{u}_j)(v_i-\tilde{u}_i)
\Big( \alpha \big( \partial_t \tilde f^{(0)} + v_\ell \partial_{x_\ell} \tilde f^{(0)} \big) +  v_\ell \partial_{y_\ell} \tilde f^{(0)}
\Big) d\bm{v}.
\label{T1}
\end{eqnarray}
We analyze now each term of equation \eqref{T1} separately. We start by applying the chain rule first to
\begin{equation}\label{chain}
\partial_{{y_\ell}}\tilde{f}^{(0)}=\partial_{\tilde{\rho}}\tilde{f}^{(0)}\partial_{{y_\ell}}\tilde{\rho}
+\partial_{\tilde{u}_{m}}\tilde{f}^{(0)}\partial_{y_\ell}\tilde{u}_{m}.
\end{equation}
We then observe that we can write the following relation
\begin{equation}
\begin{aligned}
v_\ell \frac{\partial \tilde{f}^{(0)}}{\partial y_\ell} &=v_\ell \frac{\partial \tilde{f}^{(0)}}{\partial \tilde{\rho}}\frac{\partial
\tilde{\rho}}{\partial y_\ell} + v_\ell \frac{\partial \tilde{f}^{(0)}}{\partial \tilde{u}_m}
\frac{\partial \tilde{u}_m}{\partial y_\ell} \\
&=v_\ell \frac{\tilde{f}^{(0)}}{\tilde{\rho}} \frac{\partial \tilde{\rho}}{\partial y_\ell} + \frac{v_\ell (v_m - \tilde{u}_m)}{RT}
\tilde{f}^{(0)} \frac{\partial \tilde{u}_m}{\partial y_\ell} =: A_1 + A_2,
\end{aligned}
\label{T12}
\end{equation}
to separate different contributions to the stress tensor, i.e. the one, denoted by $A_1$, is due to the mass variation and the other one,
denoted by $A_2$, is due to the variation of the mean velocity in the cell $[-1,1]^d$. Now, by introducing the notation
$\langle \cdot \rangle \equiv \int_{\mathbb{R}^d} \cdot d\bm{v}$ to indicate integration over the velocity space, one can compute the
moments for the terms $A_1$ and $A_2$ defined in \eqref{T12} separately. To that aim, with the notation
$\tilde A_k = \langle (v_j-\tilde{u}_j)(v_i-\tilde{u}_i) A_k \rangle$, with $k=1, \, 2$, for the second order moments of $A_1$ and $A_2$
we obtain:
\begin{equation*}
\begin{aligned}
\tilde A_1&:=\ \frac{1}{\tilde{\rho}}\frac{\partial \tilde{\rho}}{\partial y_\ell} \langle (v_j - \tilde{u}_j) (v_i - \tilde{u}_i)v_\ell
\tilde{f}^{(0)} \rangle\\
& \qquad = \frac{1}{\tilde{\rho}}\frac{\partial \tilde{\rho}}{\partial y_\ell} \langle (v_j - \tilde{u}_j) (v_i - \tilde{u}_i) (v_\ell -
\tilde{u}_\ell) \tilde{f}^{(0)} \rangle +  \frac{1}{\tilde{\rho}}\frac{\partial \tilde{\rho}}{\partial y_\ell} \tilde{u}_\ell \langle (v_j
- \tilde{u}_j) (v_i - \tilde{u}_i) \tilde{f}^{(0)} \rangle\\
& \qquad = \frac{\partial \tilde{\rho}}{\partial y_\ell} \tilde{u}_\ell R T \delta_{ij},\\
\tilde A_2&:=\ \frac{1}{RT}\frac{\partial \tilde{u}_m}{\partial y_\ell}\langle (v_j - \tilde{u}_j)
(v_i - \tilde{u}_i) (v_m - \tilde{u}_m)v_\ell \tilde{f}^{(0)} \rangle\\
& \qquad = \frac{1}{RT}\frac{\partial \tilde{u}_m}{\partial y_\ell}\langle (v_j - \tilde{u}_j) (v_i - \tilde{u}_i) (v_m - \tilde{u}_m)
(v_\ell - \tilde{u}_\ell) \tilde{f}^{(0)} \rangle+\frac{\tilde u_m}{RT}\frac{\partial \tilde{u}_m}{\partial y_\ell}\langle (v_j -
\tilde{u}_j) (v_i - \tilde{u}_i) (v_m - \tilde{u}_m)\tilde{f}^{(0)} \rangle \\
& \qquad = \frac{\partial \tilde{u}_m}{\partial y_\ell} \tilde{\rho} R T \big( \delta_{ij}\delta_{m \ell}+\delta_{im}\delta_{j\ell}+
\delta_{i\ell}\delta_{jm} \big).
\end{aligned}
\end{equation*}
We now analyze the term involving the time derivative $\frac{\partial \tilde{f}^{(0)}}{\partial t}$ in \eqref{T1}. Similarly to
\eqref{chain} we have
\begin{equation*}\label{time_f0}
\frac{\partial \tilde{f}^{(0)}}{\partial t} = \frac{\tilde{f}^{(0)}}{\tilde{\rho}} \frac{\partial \tilde{\rho}}{\partial t} + \frac{v_m -
\tilde{u}_m}{RT} \tilde{f}^{(0)} \frac{\partial \tilde{u}_m}{\partial t},
\end{equation*}
and analogously to \eqref{T12} with $\bm{y}$ replaced by $\bm{x}$ we have
\begin{equation*}
v_\ell \frac{\partial \tilde{f}^{(0)}}{\partial x_\ell} =v_\ell \frac{\tilde{f}^{(0)}}{\tilde{\rho}}
\frac{\partial \tilde{\rho}}{\partial x_\ell} + \frac{v_\ell (v_m - \tilde{u}_m)}{RT} \tilde{f}^{(0)}
\frac{\partial \tilde{u}_m}{\partial x_\ell}.
\label{T12y}
\end{equation*}
This leads to the following contribution for what concerns the time derivative of the distribution $\tilde f^0$
\begin{equation*}
\begin{split}
&\tilde A_3:=\Big\langle (v_j-\tilde{u}_j)(v_i-\tilde{u}_i) \frac{\partial \tilde{f}^{(0)}}{\partial t} \Big\rangle
=\Big\langle (v_j-\tilde{u}_j)(v_i-\tilde{u}_i) \left(\frac{\tilde{f}^{(0)}}{\tilde{\rho}} \frac{\partial \tilde{\rho}}{\partial t} +
\frac{v_m - \tilde{u}_m}{RT} \tilde{f}^{(0)} \frac{\partial \tilde{u}_m}{\partial t}\right) \Big\rangle=\\
&=\Big\langle (v_j-\tilde{u}_j)(v_i-\tilde{u}_i) \frac{\tilde{f}^{(0)}}{\tilde{\rho}} \frac{\partial \tilde{\rho}}{\partial t}
\Big\rangle+\Big\langle (v_j-\tilde{u}_j)(v_i-\tilde{u}_i) \frac{v_m - \tilde{u}_m}{RT} \tilde{f}^{(0)}
\frac{\partial \tilde{u}_m}{\partial t} \Big\rangle=RT\frac{\partial \tilde \rho}{\partial t} \, \delta_{ij}
= \frac{\partial \tilde p}{\partial t} \, \delta_{ij},\nonumber
\end{split}
\end{equation*}
while the contribution to the stress tensor coming from the space variation of the quantities $\tilde \rho$, $\tilde u_m$ at the
macroscopic scale $\bm{x}$ is
\begin{equation*}
\tilde A_4:= \Big\langle (v_j-\tilde{u}_j)(v_i-\tilde{u}_i) v_\ell \frac{\partial \tilde{f}^{(0)}}{\partial x_\ell} \Big\rangle
= \frac{\partial \tilde{\rho}}{\partial x_\ell} \tilde{u}_\ell R T \delta_{ij}+
\frac{\partial \tilde{u}_m}{\partial x_\ell} \tilde{\rho} R T \big( \delta_{ij}\delta_{m \ell}
+\delta_{im}\delta_{j\ell}+\delta_{i\ell}\delta_{jm} \big).
\end{equation*}
Collecting all these terms, we obtain the following form of $\mathcal{T}^{(1)}_{ij}=-\hat{\tau}\varepsilon\left(\tilde A_1+\tilde A_2+\alpha (\tilde A_3+\tilde A_4)\right)$:
\begin{eqnarray}
\mathcal{T}^{(1)}_{ij} &=& -\varepsilon \hat{\tau} RT \tilde{\rho} \Big[  \frac{\partial \tilde{u}_m}{\partial y_\ell} \left( \delta_{ij}
\delta_{m \ell} + \delta_{im} \delta_{j \ell} + \delta_{i \ell}\delta_{jm} \right) + \delta_{ij} \frac{\tilde{u}_\ell}{\tilde{\rho}}
\frac{\partial \tilde{\rho}}{\partial y_\ell}  \Big] \nonumber \\
&&- \varepsilon \alpha \hat \tau  \Big[ \frac{\partial \tilde p}{\partial t} \delta_{ij} + RT \tilde \rho \Big( \frac{\partial \tilde u_m}
{\partial x_\ell} \left( \delta_{ij}\delta_{m \ell} + \delta_{im} \delta_{j \ell} + \delta_{i \ell}\delta_{jm} \right) + \delta_{ij}
\frac{\tilde{u}_\ell}{\tilde{\rho}}\frac{\partial \tilde{\rho}}{\partial x_\ell} \Big) \Big].
\label{T1a}
\end{eqnarray}

In the matrix form, the contribution to the stress tensor given by the truncated Chapmann-Enskog expansion can be written as
\begin{eqnarray}
\mathcal{T}^{(1)} &=& - \varepsilon \hat{\tau}RT\tilde{\rho} \Big[ \nabla_{\bm{y}} \tilde{\bm{u}} + \nabla_{\bm{y}} \tilde{\bm{u}}^T  +
\Big( \nabla_{\bm{y}}\cdot \tilde{\bm{u}}
+ \big( \frac{\tilde{\bm{u}}}{\tilde{\rho}} \cdot \nabla_{\bm{y}} \big) \tilde{\rho} \Big) \mathbb{I} \Big] \nonumber \\
&&\nonumber- \varepsilon \alpha \hat \tau  \Big[ \partial_t \tilde p \, \mathbb{I} + RT \tilde \rho \Big[ \nabla_{\bm{x}} \tilde{\bm{u}} +
\nabla_{\bm{x}} \tilde{\bm{u}}^T  + \Big( \nabla_{\bm{x}} \cdot \tilde{\bm{u}}
+ \big( \frac{\tilde{\bm{u}}}{\tilde{\rho}} \cdot \nabla_{\bm{x}} \big) \tilde{\rho} \Big) \mathbb{I} \Big]= \\
&=:& \mathcal{T}^{(1,1)} + \mathcal{T}^{(1,2)}.
\label{T1b}
\end{eqnarray}
Combining now \eqref{eq:expressCST}, \eqref{T0a}, \eqref{T1b}, and discarding the terms ${\mathcal O}(\varepsilon^2)$, yields
\[
\mathbb{T} (\bm{x},t)=-\frac{1}{2^d}\int_{[-1,1]^d}\left(\mathcal{T}^{(0)}+\mathcal{T}^{(1)}\right) d\bm{y}
\]
and after substituting plugging this into \eqref{eq2} we finally obtain
\begin{equation}
\begin{aligned}
& \partial_t \rho +\nabla_{\bm{x}} \cdot (\rho \bm{u}) = 0,\\
& \partial_t (\rho \bm{u}) + \nabla_{\bm{x}} \cdot (\rho \bm{u} \otimes \bm{u})=\\
& \hspace{1cm} \nabla_{\bm{x}} \cdot \Big[  \frac{1}{2^d}
\int\limits_{[-1,1]^{d}}  \Big\{ - \tilde{p} \mathbb{I}  + \varepsilon \hat{\tau} RT \tilde{\rho}
\Big( \nabla_{\bm{y}} \tilde{\bm{u}} + \nabla_{\bm{y}} \tilde{\bm{u}}^T
+ \big[ \nabla_{\bm{y}} \cdot \tilde{\bm{u}} + \big( \frac{\tilde{\bm{u}}}{\tilde{\rho}}\cdot \nabla_{\bm{y}} \big) \tilde{\rho} \big]
\mathbb{I} \Big) \\
& \hspace{2cm} \varepsilon \hat{\tau} \alpha \big\{
\partial_t \tilde{p} \,  \mathbb{I}  + RT \tilde{\rho}
\Big( \nabla_{\bm{x}} \tilde{\bm{u}} + \nabla_{\bm{x}} \tilde{\bm{u}}^T
+ \big[ \nabla_{\bm{x}} \cdot \tilde{\bm{u}} + \big( \frac{\tilde{\bm{u}}}{\tilde{\rho}}\cdot \nabla_{\bm{x}} \big) \tilde{\rho} \big]
\mathbb{I} \Big) \big\} \Big\} d\bm{y} \Big] .
\end{aligned}
\label{CE_HM_Momentum}
\end{equation}

We proceed in the next section by proposing a closing strategy for the system \eqref{CE_HM_Momentum}. This will be established by relating
the perturbations $\tilde{\bm{u}}, \tilde{\rho}, \tilde{p}$ to the mean velocity $\bm{u}$ and density $\rho$ of the fluid at the
macroscopic scale. With this aim, we will assume that $\tilde{\bm{u}}$ and $\tilde{\rho}$ have approximate polynomial variations w.r.t.
$\bm{y}$ in the unit cell. More precisely, we assume that $\tilde{\rho}(\bm{x},\bm{y},t)\approx \rho(\bm{x}+\alpha\bm{y},t)$ and similarly
for $\tilde{\bm{u}}\approx \bm{u}(\bm{x}+\alpha\bm{y},t).$  We will use the Taylor expansion and successive truncation with respect to the
parameter $\alpha$.

\subsection{System closure}
In order to get a closed system of equations we proceed  with inserting our assumption on a linear variation of the perturbed quantities
$\tilde{\bm{u}}$ and $\tilde{\rho}$ into \eqref{CE_HM_Momentum}. We assume that both $\rho$ and $\bm{u}$ are sufficiently regular to  be
Taylor-expanded from the center of the box. In our expansion, we keep the terms up to order ${\mathcal O}(\alpha^2)$. Thus, $\rho
(\bm{y},t), \bm{u}(\bm{y},t)$ can be written as:
\begin{equation}
\begin{aligned}
&\tilde{\rho}(\bm{x},\bm{y},t) = \rho(\bm{x},t) + \alpha (\bm{y}\cdot \nabla_{\bm{x}})\rho(\bm{x},t) + \frac{\alpha^2}{2} (\bm{y}\cdot
\nabla_{\bm{x}})^2 \rho(\bm{x},t) + {\mathcal O}(\alpha^3), \\
&\tilde{\bm{u}}(\bm{x},\bm{y},t) = \bm{u}(\bm{x},t) + \alpha (\bm{y}\cdot \nabla_{\bm{x}})\bm{u}(\bm{x},t) + \frac{\alpha^2}{2} (\bm{y}
\cdot \nabla_{\bm{x}})^2 \bm{u}(\bm{x},t) + {\mathcal O}(\alpha^3).
\end{aligned}
\label{ru2}
\end{equation}
Substituting \eqref{ru2} into \eqref{T0} and \eqref{T1a}, integrating over the microscopic box with respect to $\bm{y}$ and noting that
odd terms with respect to $\bm{y}$ cancel by antisymmetry, we obtain:
\begin{equation}
\begin{aligned}
-\frac{1}{2^d} \int\limits_{[-1,1]^d}\mathcal{T}^{(0)}\,d\bm{y}
&=-\frac{1}{2^d}\int\limits_{[-1,1]^d}\tilde{\rho}(\bm{x},\bm{y},t) R T\mathbb{I}\,d\bm{y}\\
&=-\frac{1}{2^d}\int\limits_{[-1,1]^d} \left[ \rho + \alpha (\bm{y}\cdot \nabla_{\bm{x}})\rho + \frac{\alpha^2}{2}(\bm{y}\cdot
\nabla_{\bm{x}})^2\rho \right]RT\mathbb{I} \,d\bm{y}  + {\mathcal O}(\alpha^3)\\
&=-\left(\rho + \frac{\alpha^2}{6}\Delta_{\bm{x}}\rho \right) RT\mathbb{I}  + {\mathcal O}(\alpha^3).
\end{aligned}
\label{P_0_quadratic}
\end{equation}
Now taking the $\bm{y}$-derivative of the second equation in \eqref{ru2} relative to the expansion of the velocity $\tilde{\bm{u}}$, we
have
$$
\nabla_{\bm{y}}\tilde{\bm{u}} = \alpha \left( \nabla_{\bm{x}} \bm{u} + \alpha(\bm{y}\cdot\nabla_{\bm{x}})\nabla_{\bm{x}}\bm{u} \right)  +
{\mathcal O}(\alpha^3) ,
$$
with similar relations holding true for $(\nabla_{\bm{y}} \tilde{\bm{u}})^T$, $(\nabla_{\bm{y}} \cdot \tilde{\bm{u}})$, and $
(\nabla_{\bm{y}} \tilde \rho)$. Using again that odd terms in $\bm{y}$ cancel by antisymmetry, we get
\begin{equation*}
\begin{split}
&- \frac{1}{2^d} \int\limits_{[-1,1]^d} \mathcal{T}^{(1,1)}\,d\bm{y}
 =\frac{\varepsilon \hat{\tau}RT}{2^d} \int\limits_{[-1,1]^d}   \tilde{\rho} \big( \nabla_{\bm{y}} \tilde{\bm{u}} + (\nabla_{\bm{y}}
 \tilde{\bm{u}})^T + (\nabla_{\bm{y}} \cdot \tilde{\bm{u}}) \, \mathbb{I} \big)\,d\bm{y} +\frac{\varepsilon \hat{\tau}RT}{2^d}
 \int_{[-1,1]^d} (\tilde{\bm{u}}\cdot \nabla_{\bm{y}}) \tilde{\rho} \, \mathbb{I}   \,d\bm{y} \\
& = \frac{\varepsilon \hat{\tau}RT \alpha}{2^d} \int\limits_{[-1,1]^d} \Big[ \rho + \alpha(\bm{y}\cdot\nabla_{\bm{x}}) \rho
\Big] \times
\\
& \phantom{mmmmmmmmmmmmm}\Big[  \nabla_{\bm{x}} \bm{u} + (\nabla_{\bm{x}} \bm{u})^T + (\nabla_{\bm{x}} \cdot \bm{u})\, \mathbb{I}+
\alpha(\bm{y}\cdot \nabla_{\bm{x}}) \Big( \nabla_{\bm{x}} \bm{u} + (\nabla_{\bm{x}} \bm{u})^T + (\nabla_{\bm{x}} \cdot \bm{u})\mathbb{I}
\Big) \Big]\,d\bm{y} \\
&+ \frac{\varepsilon \hat{\tau}RT \alpha}{2^d} \int\limits_{[-1,1]^d}\Big[ \bm{u} + \alpha(\bm{y}\cdot\nabla_{\bm{x}}) \bm{u}
\Big] \cdot \Big[ \nabla_{\bm{x}}\rho + \alpha(\bm{y}\cdot \nabla_{\bm{x}})(\nabla_{\bm{x}}\rho)  \Big] \,
\mathbb{I} d\bm{y} + {\mathcal O}(\alpha^3)\\
&= \varepsilon \hat{\tau}RT \alpha \Big\{ \rho
\Big[ \nabla_{\bm{x}} \bm{u} + (\nabla_{\bm{x}} \bm{u})^T + (\nabla_{\bm{x}} \cdot \bm{u}) \, \mathbb{I} \Big]
+ \bm{u} \cdot \nabla_{\bm{x}}\rho \, \mathbb{I} \Big\} + {\mathcal O}(\alpha^3).
\end{split}
\end{equation*}
Using similar computations for the term $\mathcal{T}^{(1,2)}$ in \eqref{T1b}, we have
\begin{equation*}
\begin{aligned}
- \frac{1}{2^d} \int\limits_{[-1,1]^d} \mathcal{T}^{(1,2)}\,d\bm{y}
& =\frac{1}{2^d} \varepsilon \hat{\tau} RT \alpha \int\limits_{[-1,1]^d}  \Big\{ \partial_t \tilde \rho \, \mathbb{I} + \tilde{\rho}
\big( \nabla_{\bm{x}} \tilde{\bm{u}} + (\nabla_{\bm{x}} \tilde{\bm{u}})^T + (\nabla_{\bm{x}} \cdot \tilde{\bm{u}}) \, \mathbb{I} \big) +
(\tilde{\bm{u}}\cdot \nabla_{\bm{x}}) \tilde{\rho} \, \mathbb{I}  \Big\} \,d\bm{y}\\
& = \frac{1}{2^d} \varepsilon \hat{\tau}RT \alpha \int\limits_{[-1,1]^d} \Big\{ \Big( \partial_t \rho + \alpha (\bm{y} \cdot
\nabla_{\bm{x}}) \frac{\partial \rho}{\partial t} \Big) \, \mathbb{I} + \Big[ \rho + \alpha(\bm{y}\cdot\nabla_{\bm{x}}) \rho
\Big] \, \Big[  \nabla_{\bm{x}} \bm{u} + (\nabla_{\bm{x}} \bm{u})^T \\
& \hspace{2cm}   + (\nabla_{\bm{x}} \cdot \bm{u})\, \mathbb{I} + \alpha(\bm{y}\cdot \nabla_{\bm{x}}) \Big( \nabla_{\bm{x}} \bm{u} +
(\nabla_{\bm{x}} \bm{u})^T + (\nabla_{\bm{x}} \cdot \bm{u})\mathbb{I} \Big) \Big] \\
&\hspace{2cm} + \Big[ \bm{u} + \alpha(\bm{y}\cdot\nabla_{\bm{x}}) \bm{u}
\Big] \cdot \Big[ \nabla_{\bm{x}}\rho + \alpha(\bm{y}\cdot \nabla_{\bm{x}})(\nabla_{\bm{x}}\rho) \Big] \, \mathbb{I} \Big\} d\bm{y} +
{\mathcal O}(\alpha^3)\\
&= \varepsilon \hat{\tau}RT \alpha \Big\{ \partial_t \rho \, \mathbb{I} + \rho
\Big[ \nabla_{\bm{x}} \bm{u} + (\nabla_{\bm{x}} \bm{u})^T + (\nabla_{\bm{x}} \cdot \bm{u}) \, \mathbb{I} \Big]
+ \bm{u} \cdot \nabla_{\bm{x}}\rho \, \mathbb{I} \Big\} + {\mathcal O}(\alpha^3).
\end{aligned}
\end{equation*}
We finally assume that $\varepsilon = \sqrt{\alpha}$. This means that the spatial inhomogeneity due to the microscopic interactions of
order $\alpha$ are much smaller and specifically equal to the square of the relaxation length, i.e. the mean free-path which is of order
$\varepsilon$ as stated at the beginning of Section \ref{sec:hybrid}. Thus, we drop all terms of order $\alpha^2$ in front of terms of
order $\varepsilon \alpha$. This is the case of the term involving $\Delta_x \rho$ in \eqref{P_0_quadratic}. With this closure and using
that $\partial_t \rho = - \nabla_{\bm{x}} \cdot (\rho \bm{u})$  the system \eqref{eq2} is written:
\begin{equation}
\begin{aligned}
&\partial_t \rho+\nabla_{\bm{x}}\cdot(\rho \bm{u})=0,\\
&\partial_t (\rho \bm{u})+\nabla_{\bm{x}}\cdot(\rho \bm{u} \otimes \bm{u}) + \nabla_{\bm{x}} p  \\
&\hspace*{2.2cm} = \alpha^{3/2} \hat \tau RT \ \nabla_{\bm{x}}\cdot \Big\{ 2 \rho \Big( \nabla_{\bm{x}} \bm{u} + (\nabla_{\bm{x}}
\bm{u})^T \Big) + \nabla_{\bm{x}} \cdot (\rho \bm{u}) \, \mathbb{I} \Big\},
\end{aligned}
\label{HM2_vector_form}
\end{equation}
with $p = \rho RT.$ We note that when terms of order $\alpha^{3/2}$ are neglected, we recover the standard Euler equations.
\begin{remark}~
\begin{itemize}
\item
Using alternative closure relations with respect to the one defined in \eqref{ru2} leads to alternative systems of balance laws which may
be used to describe different phenomena related to the microscopic dynamics.
\item
Keeping additional terms in the Taylor expansion \eqref{ru2}/considering a different scaling relation with respect to
$$
\varepsilon = \sqrt{\alpha}
$$
leads also to a different balance laws of dispersive type: third order terms appears.
\item
Different kinetic models can be applied instead of the BGK relaxation equation \eqref{BGK}. For example, molecular dynamics interacting
through short range potential can be used giving rise, for instance, to alternative Piola-Kirchhoff stress tensors.
\end{itemize}
\end{remark}

\section{Linear stability analysis}\label{sec3}
In this section, we study the linear stability of system \eqref{HM2_vector_form}. Let us assume $\rho = \rho_0 + \beta \rho_1 + \cdots$,
and $\bm{u} = \bm{u}_0 + \beta \bm{u}_1 + \cdots$ with $\beta \ll 1$, where $(\rho_0,\bm{u}_0)$ represents a uniform steady state and
$(\rho-\rho_0,\bm{u}-\bm{u}_0)$ is a small perturbation from this equilibrium state. We then substitute the above expressions into
\eqref{HM2_vector_form} and, dropping terms of order $\beta^2$ or more, we obtain a system of linearized equations for $\rho_1(\bm{x},t)$
and $\bm{u}_1(\bm{x},t)$. The scope of this section is to investigate the stability properties of the resulting linear system. The
corresponding equations for $(\rho_1,\bm{u}_1)$ are
\begin{equation}
\begin{aligned}
&\partial_t \rho_1 + \rho_0 \nabla_{\bm{x}} \cdot \bm{u}_1 + \bm{u}_0 \cdot \nabla_{\bm{x}} \rho_1=0,\\
&\rho_0 \left(\partial_t \bm{u}_1+(\bm{u}_0\cdot\nabla_{\bm{x}})\bm{u}_1\right) + RT \nabla_{\bm{x}} \rho_1 \\
&\hspace*{2.2cm} = \alpha^{3/2} \hat \tau RT\  \nabla_{\bm{x}}\cdot \Big\{ 2 \rho_0 \Big( \nabla_{\bm{x}} \bm{u}_1 + (\nabla_{\bm{x}}
\bm{u}_1)^T \Big) + \Big(\rho_0 \nabla_{\bm{x}} \cdot \bm{u}_1 + \bm{u}_0 \cdot \nabla_{\bm{x}} \rho_1 \Big) \,  \mathbb{I} \Big\}.
\end{aligned}
\label{Linearized_HM2}
\end{equation}
Using the Fourier transform in both space and time, we can write
\begin{equation}
\begin{bmatrix}
\rho_1\\
\bm{u}_1
\end{bmatrix}(\bm{x},t) = \begin{bmatrix}
\bar{\rho}_1\\
\bar{\bm{u}_1}
\end{bmatrix} e^{i (\bm{k}\cdot\bm{x} - \omega t)}
\label{eq:Fourier}
\end{equation}
with $\bar\rho_1$ and $\bar{\bm{u}}_1$ being the Fourier coefficients and $\bm{k}\in \mathbb{R}^3$, $\omega\in\mathbb{C}$. Using such
transformation, one can infer that the system is stable about the
stationary solution $(\rho_0, \bm{u}_0)$ if and only if the imaginary part of $\omega$ is nonpositive for all non-trivial solution
\eqref{eq:Fourier}. Moreover, one can state that the system is stable if and only if it is stable for all
$(\rho_0, \bm{u}_0) \in [0,\infty) \times {\mathbb R}^d$.

For simplicity, before proceeding with the computations, we redefine the following quantities $RT \equiv T$,
$\alpha^{3/2} \hat \tau \equiv \tau$ and $(\bar\rho_1,\bar{\bm{u}}_1) \equiv (\rho, \bm{u})$. Substituting now \eqref{eq:Fourier} in
\eqref{Linearized_HM2} and using these redefined quantities, we get
\begin{equation*}
\begin{aligned}
&-i \omega \rho + \rho_0 i (\bm{k} \cdot \bm{u}) + i (\bm{k}\cdot \bm{u}_0) \rho =0,\\
&\rho_0 \Big(-i \omega \bm{u} + i(\bm{k} \cdot \bm{u}_0) \bm{u} \Big) + T i \bm{k} \rho  \\
&\hspace*{2.2cm} = \tau T  i \bm{k}^T  \Big\{ 2 \rho_0 \Big( i \bm{k}\otimes \bm{u} + i \bm{u} \otimes \bm{k} \Big) +
\Big( \rho_0 i (\bm{k} \cdot \bm{u}) + i (\bm{k}\cdot \bm{u}_0) \rho \Big) \mathbb{I} \Big\}.
\end{aligned}
\end{equation*}
Simplifying the above expressions, we obtain
\begin{equation}
\begin{aligned}
&\Big( - \omega + \bm{k}\cdot\bm{u}_0 \Big)  \rho + \rho_0 \bm{k}\cdot \bm{u} =0,\\
&\Big( - \omega + \bm{k}\cdot\bm{u}_0 \Big) \bm{u} + \frac{T}{\rho_0} \Big( 1 - i \tau (\bm{k}\cdot \bm{u}_0) \Big) \bm{k} \rho  -
2 i \tau T \Big( |\bm{k}|^2 \bm{u} + (\bm{k}\cdot \bm{u}) \bm{k}  \Big)= 0.
\end{aligned}
\label{Perturbed_HM2}
\end{equation}
Now, we want to analyze equation \eqref{Perturbed_HM2}. To this aim, let us first assume that $\bm{k}=0$. This choice directly implies
$\omega \not = 0$ otherwise, the perturbed solution is just like the unperturbed equation: constant in time and spatially uniform. From
the above hypothesis, we immediately obtain $-\omega \rho = 0$ and $-\omega \bm{u} = 0$. This means $\rho =0$ and $\bm{u} = 0$ and thus,
the perturbed solution is again constant in time and space, i.e. there does not exist any non-trivial solution to the system
\eqref{Perturbed_HM2} in this setting. Thus, in the sequel, we assume that $\bm{k}\not=0$.

Let thus, under the hypothesis $\bm{k}\not=0$, the orthogonal projection onto $\{\bm{k}\}^{\perp}$ be denoted as
$\mathtt{P}_{\bm{k}^{\perp}}$. Applying $\mathtt{P}_{\bm{k}^{\perp}}$ to the second equation in \eqref{Perturbed_HM2}, yields
\begin{equation}
\Big( -\omega + \bm{k}\cdot\bm{u}_0 - 2 i \tau T |\bm{k}|^2 \Big) \mathtt{P}_{\bm{k}^{\perp}} \bm{u} = 0.
\label{3.5}
\end{equation}
It should be observed that the expression inside the parentheses in \eqref{3.5} can not be zero, otherwise, by taking the imaginary part,
one concludes that $k=0$, which is a contradiction. Hence, the perturbed velocity is in the direction of $\bm{k}$ and thus
$\mathtt{P}_{\bm{k}^{\perp}} \bm{u} = 0$. As a consequence, we can write
\begin{equation}\label{uperp}
\bm{u} = v \frac{\bm{k}}{|\bm{k}|},
\end{equation}
where $v \in \mathbb{R}$. After substituting \eqref{uperp} into equation \eqref{Perturbed_HM2}, the system for $(\rho, v)$ reads
\begin{equation*}
\begin{aligned}
&\Big( - \omega + \bm{k}\cdot\bm{u}_0 \Big)  \rho + \rho_0 |\bm{k}| v =0,\\
&\Big( - \omega + \bm{k}\cdot\bm{u}_0 \Big) v + \frac{T}{\rho_0} \Big( 1 - i \tau (\bm{k}\cdot \bm{u}_0) \Big) |\bm{k}| \rho
- 4 i \tau T |\bm{k}|^2 v = 0.
\end{aligned}
\end{equation*}
We introduce now the angle $\theta$ between $\bm{u}_0$ and $\bm{k}$, hence $\bm{k}\cdot\bm{u}_0 = k u_0 \cos \theta$, where we denote
$k = |\bm{k}|$ and $u_0 = |\bm{u}_0|$. This leads to
\begin{equation*}
\begin{aligned}
&\Big( - \omega + k u_0  \cos \theta \Big)  \rho + \rho_0 k v =0,\\
&\frac{T}{\rho_0} \Big( 1 - i \tau k u_0 \cos \theta \Big)  k \rho+ \Big( - \omega + k u_0 \cos \theta  - 4 i \tau T k^2 \Big) v = 0,
\end{aligned}
\end{equation*}
or, in matrix form:
\begin{equation*}
\begin{bmatrix}
-\omega + k u_0 \cos \theta & \rho_0 k\\
\frac{T}{\rho_0} \Big( 1 - i \tau k u_0 \cos \theta \Big)  k & - \omega + k u_0  \cos \theta  - 4 i \tau T k^2
\end{bmatrix}\begin{bmatrix}
\rho\\
v
\end{bmatrix}=0
\end{equation*}
Thus, there exists a non trivial solution if and only if the determinant of the above system is equal to zero. This gives
\begin{equation*}
\left| \begin{array}{llcc}
-\omega +  k u_0 \mbox{cos}\theta & \rho_0 k\\
\frac{T}{\rho_0} \Big( 1 - i \tau k u_0 \mbox{cos}\theta \Big) k  & - \omega + k u_0 \mbox{cos}\theta  - 4 i \tau T k^2
\end{array}\right| =0 .
\end{equation*}
Letting $X = - \omega + k u_0 \cos \theta$, the previous equation leads to the following characteristic equation
\begin{equation}
X^2 - 4 i \tau T k^2 X - T k^2 (1-i \tau k u_0 \cos \theta) = 0.
\label{Characteristic_polynomial}
\end{equation}
Let us consider the set of solutions of \eqref{Characteristic_polynomial}. We want to guarantee that the imaginary part of $\omega$ remains non-positive such that the linear stability holds true. To analyse \eqref{Characteristic_polynomial}, we refer to the Routh-Hurwitz criterion \cite{ishida2013linear}. This criterion is summarized as follows: let a polynomial equation be of the form
\begin{equation}\label{RHC}
(a_0 + i b_0)\mu^n + (a_1 + i b_1)\mu^{n-1} + \cdots + (a_n + i b_n) = 0.
\end{equation}
Then, all solutions of the above equation satisfy $\mbox{Re}(i\mu) < 0$ if and only if
\begin{eqnarray*}\label{delta}
(-1)\Delta_2 &=& -\left| \begin{array}{cc}
a_0 & a_1\\
b_0 & b_1
\end{array}\right| >0, \\\label{delta1}
(-1)^2\Delta_4 &=& \left| \begin{array}{cccc}
a_0 & a_1 & a_2 & a_3\\
b_0 & b_1 & b_2 & b_3\\
0 & a_0 & a_1 & a_2\\
0 & b_0 & b_1 & b_2
\end{array}\right| >0, \\
\vdots \nonumber\\
(-1)^n\Delta_{2n} &=& (-1)^n\left| \begin{array}{ccccccccc}
a_0 & a_1 & \cdots & a_{n-1} & a_n & 0 & \cdots & \cdots & 0\\
b_0 & b_1 & \cdots & b_{n-1} & b_n & 0 & \cdots & \cdots & 0\\
0 & a_0 & \cdots & a_{n-2} & a_{n-1} & a_n & 0 & \cdots & 0\\
0 & b_0 & \cdots & b_{n-2} & b_{n-1} & b_n & 0 & \cdots & 0\\
\vdots & \vdots & \ddots & \vdots & \vdots & \vdots & \vdots & \ddots & \vdots \\
0 & 0 & \cdots & a_0 & a_1 & a_2 & a_3 & \cdots & a_n\\
0 & 0 & \cdots & b_0 & b_1 & b_2 & b_2 & \cdots & b_n
\end{array}\right| >0.
\end{eqnarray*}
Now, observing that $\textrm{Re}(iX) = \textrm{Im} (\omega)$, the stability conditions follow immediately. We then apply the Routh-Hurwitz
criterion directly to \eqref{Characteristic_polynomial}. To this aim, we compare the coefficients of \eqref{Characteristic_polynomial}
with those of \eqref{RHC}. This gives
\begin{equation}
\begin{array}{ll}
a_0 = 1,& b_0 = 0\\
a_1 = 0 ,& b_1 = - 4 \tau T k^2\\
a_2 = -T k^2 ,& b_2 = \tau T k^3 u_0 \cos \theta.
\end{array}
\end{equation}
Using the above quantities in the Routh-Hurwitz criterion we get
\begin{equation*}
-\Delta_2 = 4 \tau T k^2 > 0
\end{equation*}
and
\begin{equation*}
\Delta_4 = \tau^2 T^2 k^6 ( 16 T - u_0^2 \cos^2 \theta).
\end{equation*}
Thus, we can conclude that if $u_0^2 \geq 16 T$, then, $\Delta_4 \geq 0$, i.e. the model is stable for all values of $\theta$. On the
other hand, if $u_0^2 > 16 T$, then, there is a threshold $\cos \theta_c$ defined by
$$ \cos^2 \theta_c = \frac{16 T}{u_0^2}, $$
such that the model is stable if $|\cos \theta| \leq \cos \theta_c$ and unstable if $|\cos \theta| > \cos \theta_c$. In this case, the
model is unstable for waves propagating in directions close to the direction of the unperturbed velocity.

To summarize the above stability analysis we have showed that when the flow is subsonic the system is certainly linearly stable. Moreover,
even for supersonic flows the system remains stable, in fact for the model considered the sound speed corresponds to $c=T$. Only for
hypersonic flows (Mach number $M>4$) we have the appearance of linearly unstable modes.  However, for the applications we have in mind, we
always consider regimes where the fluid velocity is small compared with the thermal velocity, and so, the stability criterion is always
verified.

\section{Numerical methods and experiments}\label{sec4}
In this section, we study the behavior of the hybrid continuum-kinetic model \eqref{HM2_vector_form} derived in Section~\ref{sec2} by
comparing it with the isentropic Euler equations, the isentropic Navier-Stokes equations, and the BGK model. In what follows, we shall
refer to our hybrid multiscale model \eqref{HM2_vector_form} as HMM. We expect the HMM to improve the results of the standard macroscopic
model and to approach the behavior of the BGK equation, at least in specific regimes that have been used in its derivation.

In the rest of this section we first describe the numerical methods used for computations and then illustrate the performance of the
derived hybrid continuum-kinetic model in a number of experiments that clearly illustrate that the HMM is capable to better describe the
physics of complex fluids in several regimes.

\subsection{The BGK model and its numerical discretization}
The HMM has been derived using the  Chapmann-Enskog expansion starting from a steady state BGK model, cf. Section~\ref{sec:hybrid}. It is
therefore natural to expect the HMM to be close to the underlying kinetic equation in several regimes of the relaxation parameter $\tau$.
Consequently, we consider the following time dependent BGK equation
\be \label{bgk}
\partial_t f+\bm{v}\cdot\nabla_{\bm{x}} f=\frac{1}{\tau}(M_{(\rho, \bm{u})}-f),
\ee
where as assumed in the derivation of the HMM, the Maxwellian distribution \eqref{max} has constant temperature.  In particular, we
consider $T=1$ and also fix the gas constant $R=1$. The shape of this equilibrium distribution is given in equation \eqref{max}. Now, the
mean velocity $\bm{u}$ and the density $\rho$ are those obtained from integration of the distribution $f$ in velocity space:
\begin{equation*}
\rho=\int_{\mathbb{R}^d} f d\bm{v}, \qquad \bm{u}=\int_{\mathbb{R}^d} \bm{v}f d\bm{v}.
\end{equation*}
For the sake of comparison with the other macroscopic models, we choose and fix the dimension of the velocity space as $d=2$ both for the
BGK as well as for the HMM. Then, in order to numerically approximate \eqref{bgk} (see \cite{ACTA} for details), we first replace the
unbounded velocity space with a suitable sufficiently large bounded set. This implies the truncation of the tails of the distribution
function, which normally lives in a non-compactly supported set. Successively, we replace our continuous model by a so-called Discrete
Velocity Model (DVM) by discretizing this new bounded space by means of a finite number of discrete points representing the discrete
velocities that the particles can assume. The result of this procedure is that the continuous BGK model is replaced by $N$ linear
transport equations coupled through a suitable discretization of the relaxation operator $(M_{(\rho, \bm{u})}-f)$. We now introduce the
method and the notations, taking inspiration from \cite{Mieussens}. We work on a Cartesian grid $\mathcal{V}$ with
\begin{equation*}
\mathcal{V}=\left\{ \bm{v}_{\bm{k}}=\bm{k}\Delta \bm{v}+a, \ \bm{k}=(k^{(1)},k^{(2)}), \ \bm{a}=(a_1,a_2)\right\},
\label{disc_space}
\end{equation*}
where $\bm{a}$ is an arbitrary vector, $\Delta \bm{v}$ is a constant mesh size in velocity and where the components of the index $\bm{k}$
have some given bounds $K^{(1)}, K^{(2)}$. In this setting, the continuous distribution function $f$ is replaced by the vector
$f_{\mathcal{K}}(\bm{x},t)$ of size $N$, where $N$ is chosen as a compromise between accuracy and computational cost. Each component of
this vector is assumed to be an approximation of the distribution function $f$ at location $\bm{v}_{\bm{k}}$:
$$
f_{\mathcal{K}}(\bm{x},t)=(f_{\bm{k}}(\bm{x},t))_{\bm{x}},\qquad f_{\bm{k}}(\bm{x},t) \approx
f(\bm{x},\bm{v}_{\bm{k}},t).
$$
Thus, the discrete ordinate kinetic model consists of the following system of ODEs to be solved
\be
\partial_t f_{\bm{k}} + \bm{v}_{\bm{k}} \cdot\nabla_{\bm{x}}f_{\bm{k}} = \frac{1}{\tau} (M_{\bm{k}}-f_{\bm{k}}), 
\label{eq:DM1}
\ee
with $M_{\bm{k}}\approx M(\bm{x},\bm{v}_{\bm{k}},t)$ being a suitable approximation of $M_{\rho,\bm{u}}(\bm{x},\bm{v}_{\bm{k}},t)$.

The system \eqref{eq:DM1} is discretized in space using standard WENO approaches of order three and we do not detail them here. With
respect to the time discretization, it should be observed, that the Maxwellian distribution in \eqref{eq:DM1} depends on the distribution
function $f$ through its moments, cf. \eqref{max}, and hence the time integration of this ODE system is implemented using an
implicit-explicit (IMEX) Runge-Kutta method; see, e.g., \cite{Dimarco_stiff2} and references therein. A general formulation of the IMEX
Runge-Kutta method for \eqref{eq:DM1} can be written as
\begin{equation}
	\begin{aligned}
	F_{\bm{k}}^{(i)} &= \displaystyle f_{\bm{k}}^{n}-\Delta t \sum_{j=1}^{i-1} \ta_{ij} \bm{v}_{\bm{k}}\cdot\nabla_{\bm{x}} F_k^{(j)}+
	\Delta t\sum_{j=1}^{\nu} a_{ij}\frac{1}{\tau} \left(M^{(j)}_{\bm{k}}-F_{\bm{k}}^{(j)}\right), \\ 
	f_{\bm{k}}^{n+1} &= \displaystyle f_{\bm{k}}^{n}-\Delta t
	\sum_{i=1}^{\nu}\tw_{i}\bm{v}_{\bm{k}}\cdot\nabla_{\bm{x}} F_{\bm{k}}^{(i)}+\Delta
	t\sum_{i=1}^{\nu}w_{i}\frac{1}{\tau}\left(M_{\bm{k}}^{(i)}-F_{\bm{k}}^{(i)}\right), \ 
	\end{aligned}
	\label{eq:GIMEX}
\end{equation}
where the matrices $ \tA=(\ta_{ij} )$, $\ta_{ij} = 0$ for $j \geq i$ and $A = (a_{ij})$ are $\nu\times \nu$ matrices such that the
resulting scheme is explicit in $\bm{v}_{\bm{k}}\cdot\nabla_{\bm{x}} f$, and implicit in $(M_{\rho,\bm{u}}-f)$. Here, we use the so-called
second order in time ARS(2,2,2) scheme \cite{Ascher}, for which the coefficient vectors $\tw =( \widetilde{w}_{1},..,\tw_{\nu})^{T}$,
$w =(w_{1},..,w_{\nu})^{T}$ are determined by the following double Butcher tableau:
\[
	\begin{array}{c|ccc}
		0 & 0 & 0 & 0   \\
		\gamma_1   & \gamma_1 & 0 & 0 \\
		1   & \gamma_2 & 1-\gamma_2 & 0\\
		\hline
		& \gamma_2 &  1-\gamma_2 & 0
	\end{array}\qquad
	\begin{array}{c|ccc}
		0 & 0 & 0 & 0 \\
		\gamma_1   & 0 & \gamma_1 & 0 \\
		1   & 0 & 1-\gamma_1 & \gamma_1\\
		\hline
		& 0 & 1-\gamma_1 & \gamma_1
	\end{array}
\]
with $\gamma_1=1-1/\sqrt{2}$ and $\gamma_2=1-1/(2\gamma_1)$. The above scheme belongs to a particular class of IMEX methods for which the
implicit tableau is simply diagonally implicit, i.e. it is such that $a_{ij}=0$ if $j>i$. This permits a direct evaluation of the implicit
terms without resorting to the inversion of non linear systems despite the nonlinearity of the function which defines the equilibrium
state $M_{\rho,\bm{u}}$. Indeed, let us remark that the stage evaluation \eqref{eq:GIMEX} can be rewritten as
\begin{equation}
F_{\bm{k}}^{(i)} = \displaystyle f_{\bm{k}}^{n}-\Delta t \sum_{j=1}^{i-1} \ta_{ij} \bm{v}_{\bm{k}}\cdot\nabla_{\bm{x}} F_{\bm{k}}^{(j)}
+\Delta t\sum_{j=1}^{i-1} a_{ij}\frac{1}{\tau} \left(M^{(j)}_{\bm{k}}-F_{\bm{k}}^{(j)}\right)+\Delta t \frac{a_{ii}}{\tau}
\left(M^{(i)}_{\bm{k}}-F_{\bm{k}}^{(i)}\right),
\label{eq:GIMEX2}
\end{equation}
where the only implicit term is the diagonal factor $\frac{a_{ii}}{\tau} \left(M^{(i)}_{\bm{k}}-F_{\bm{k}}^{(i)}\right)$, in which
$M_{\bm{k}}^{(i)}$ depends on the density and momentum of the distribution function $(\rho,\bm{u})$. These macroscopic quantities can be
obtained from equation (\ref{eq:GIMEX2}) by discrete integration in velocity space against $\bm{\phi}_{\bm{k}}=(1,\bm{v}_{\bm{k}})$:
\begin{equation}
\begin{aligned}
(\rho^{(i)}(\bm{x}),\bm{u}^{(i)}(\bm{x}))&=\sum_{\bm{k}}(1,\bm{v}_{\bm{k}})F^{(i)}_{\bm{k}}(\bm{x})\, \Delta \bm{v}:=
\langle \bm{\phi}_{\bm{k}}  F_{\bm{k}}^{(i)}\rangle_{\mathcal{K}}\\
&=\langle \bm{\phi}_{\bm{k}} f_{\bm{k}}^{n}\rangle_{\mathcal{K}}-\Delta t \sum_{j=1}^{i-1} \ta_{ij}\langle \bm{\phi}_{\bm{k}}
\bm{v}_{\bm{k}}\cdot\nabla_x(F_{\bm{k}}^{(j)})\rangle_{\mathcal{K}}.
\end{aligned}
\label{mac}
\end{equation}
As a consequence from the calculation performed in \eqref{mac}, $(\rho^{(i)}(x),\bm{u}^{(i)}(x))$, and thus $M_{\bm{k}}^{(i)}$, can be
explicitly evaluated and then the scheme (\ref{eq:GIMEX})-(\ref{eq:GIMEX2}) is, in fact, explicitly solvable.

The choice of this specific IMEX Runge-Kutta ODE solver reflects the facts that we need to compare different model acting at different
scales and the scheme should allow to handle consistently the passage from the kinetic to the fluid equations. In particular, when the
scaling parameter $\tau\to 0$, the scheme is formally equivalent to a discretization of the inviscid equation  \eqref{eq2}. For
small but non zero values of $\tau$, one can also expect the scheme to be close to the new HMM. An analysis of the above scheme can be
performed showing that indeed it possesses the property of being consistent with the limit macroscopic inviscid model \eqref{eq2} when
$\tau\to 0$; see, e.g., \cite{Dimarco_stiff2}.

\subsection{A numerical method for the hybrid continuum-kinetic model (HMM)}
We continue by introducing the numerical method for the HMM. We give the details of the two-dimensional discretization since the following
numerical experiments will be restricted to two space dimensions. Let us fix, as for the BGK case, the gas constant $R=1$ as well as the
temperature $T=1$, set also $\alpha^{3/2}\hat\tau\equiv \tau_{HMM}.$  HMM equations \eqref{HM2_vector_form} rewritten in the vector form
read
\begin{equation}\label{Euler_eq}
\U{Q}_t + \U{F}_{x} + \U{G}_{y} = \mathbf{S},
\end{equation}
where the vectors of conservative variables $\Q$, the fluxes $(\U{F}(\Q),\U{G}(\Q)$ and the source terms $\U{S}(\Q)$ are
\begin{equation} \label{eq:euler_sys}
\begin{aligned}
	\U{Q} = \begin{bmatrix} \rho \\ \rho u_x \\ \rho u_y   \end{bmatrix}, \quad &
	\U{F}(\Q) = \begin{bmatrix} \rho u_x \\ \rho u_x^2 + p \\ \rho u_x u_y         \end{bmatrix}, \quad
	\U{G}(\Q) =\begin{bmatrix} \rho u_y \\ \rho u_y u_x     \\ \rho u_y^2 + p       \end{bmatrix}, 	\\[1ex]
	&\U{S}(\Q) =\tau_{HMM} \begin{bmatrix} 0 \\ (\xi_{xx})_x+(\xi_{yx})_y
		\\ (\xi_{xy})_x+(\xi_{yy})_y
	\\
	\end{bmatrix}.
	\end{aligned}
\end{equation}
Here, $\uu=(u_x,u_y)$ is the macroscopic velocity and
\begin{equation}
\begin{aligned}
&\xi_{xx}=4\rho\frac{\partial u_x}{\partial x}+\frac{\partial (\rho u_x)}{\partial x}+\frac{\partial (\rho u_y)}{\partial y},\
\xi_{yy}=4\rho\frac{\partial u_y}{\partial y}+\frac{\partial (\rho u_y)}{\partial y}+\frac{\partial (\rho u_x)}{\partial x},\\
&\xi_{xy}=\xi_{yx}=2\rho\left(\frac{\partial u_x}{\partial y}+\frac{\partial u_y}{\partial x}\right),
\end{aligned}
\label{eqn.tauVisc}
\end{equation}
We also use the notation $\xx=(x,y) \in \Omega$, where $\Omega$ is the computational domain. We assume that
$\Omega=[x_{\min},x_{\max}]\times[y_{\min},y_{\max}]$ is paved with $N_x \times N_y$ uniform cells of size $\Delta x \times \Delta y$.
A cell is labelled by two indices ${i,j}$, one for each direction, while when referring to an interface between the cells we use
respectively $({\ihp,j})$ or $({i,\jhp})$. Using this notation, the cell center is located at point $\U{x}_{i,j}=(x_i,y_j)$ while a face
center lies at point
$$
\U{x}_{\ihp,j}=\left( \frac12 (x_i+x_{i+1}),y_j \right),\ \U{x}_{i,\jhp}=\left( x_i,\frac12 (y_j+y_{j+1})\right).
$$
We denote any generic cell-centered quantity $m_{i,j}$, i.e. in the following density, momentum, pressure and mean velocity. A
conservative finite volume scheme is adopted along with the first order explicit Euler scheme for time integration, and we define the
generic explicit operator $F[m_{i,j}^{n}]$ applied to $m_{i,j}^{n}$ as
\begin{equation*}
F[m_{i,j}^{n}] = m_{i,j}^n-\dtdx\left(f_{\ihp,j}^{m}-f_{\ihm,j}^{m}\right)-\dtdy\left(g_{i,\jhp}^{m}-g_{i,\jhm}^{m}\right),
\end{equation*}
where the numerical fluxes $f_{\ihp,j}^{m}$ and $g_{i,\jhp}^{m}$ are of Rusanov type:
\begin{eqnarray}
f_{\ihp,j}^{m}&=&\frac{1}{2} \left( f(m_{i+1,j}^n) + f(m_{i,j}^n) \right)
-\frac{1}{2} | \lambda_{i+1/2,j}^{x,\text{max}} | \left( m_{i+1,j}^n - m_{i,j}^n  \right), \label{eq:q3d}\\	
g_{i,\jhp}^{m} &=& \frac{1}{2}\left( g(m_{i,j+1}^n) + g(m_{i,j}^n) \right)
	- \frac{1}{2} | \lambda_{i,j+1/2}^{y,\text{max}} | \left( m_{i,j+1}^n - m_{i,j}^n  \right) \label{eq:q3d1}
\end{eqnarray}
with the eigenvalues at the interfaces given by
\begin{eqnarray}
	| \lambda_{\ihp,j}^{x,\text{max}} | = \max \left( |\lambda_{i+1,j}^{x,n}|,|\lambda_{i,j}^{x,n}| \right), \
	| \lambda_{i,\jhp}^{y,\text{max}} | = \max \left( |\lambda_{i,j+1}^{y,n}|,|\lambda_{i,j}^{y,n}| \right).\nonumber
\end{eqnarray}
In this setting, the density $\rho^{n}_{i,j}$ can be directly computed as
\begin{equation*}\label{eqn.rhomult}
	\rho_{i,j}^{n+1}= F[\rho^{n}_{i,j}],
\end{equation*}
and the components of the momentum equations satisfy
\begin{eqnarray*}
	\label{eqn_momx}
	&(\rho {u})_{i,j}^{n+1}= F[(\rho {u})_{i,j}^{n}]
	+\frac{\Delta t}{\Delta x}\left( (\xi_{xx})_{\ihp,j}^{n} - (\xi_{xx})_{\ihm,j}^{n} \right)+\frac{\Delta t}{\Delta y}
	\left( (\xi_{yx})_{i,\jhp}^{n} - (\xi_{yx})_{i,\jhm}^{n} \right), \\
	\label{eqn_momy}
	&(\rho {v)}_{i,j}^{n+1}=  F[(\rho {v})_{i,j}^{n}]
	+\frac{\Delta t}{\Delta x}\left( (\xi_{xy})_{\ihp,j}^{n} - (\xi_{xy})_{\ihm,j}^{n} \right)+\frac{\Delta t}{\Delta y}
	\left( (\xi_{yy})_{i,\jhp}^{n} - (\xi_{yy})_{i,\jhm}^{n} \right).
	\label{eqn_momz}
\end{eqnarray*}
According to the definitions \eqref{eqn.tauVisc} the computation of the stress components $(\xi_{xx},\xi_{xy}$, $\xi_{yx},\xi_{yy})$
requires the knowledge of the discrete velocity gradients. Here, they are computed on each boundary face of the control volume using a
trapezoidal quadrature rule, that is,
\begin{equation}
	\nabla_{\xx} \uu_{\ihp,j} = \frac{1}{4} \left( \nabla_{\xx} \uu_{\ihp,\jhm} + \nabla_{\xx} \uu_{\ihp,\jhp}
	+ \nabla_{\xx} \uu_{\ihp,\jhm} + \nabla_{\xx} \uu_{\ihp,\jhp}  \right),
	\label{eqn.face_grad}
\end{equation}
where the two-dimensional gradient of the velocity field $\uu(t,\xx)$ is obtained at each corner defined by subscript index
$({\ihp,\jhp})$ as follows:
\begin{equation}
	\nabla_{\xx} \uu_{\ihp,\jhp} = \frac{1}{2}\left[\begin{array}{l}	
		\frac{\uu_{\ip,j}-\uu_{i,j}}{\dx} + \frac{\uu_{\ip,\jp}-\uu_{i,\jp}}{\dx}   \\[10pt]		
		\frac{\uu_{i,\jp}-\uu_{i,j}}{\dy} + \frac{\uu_{\ip,\jp}-\uu_{\ip,j}}{\dy}
	\end{array}\right].\nonumber
\end{equation}
Finally, to increase the accuracy in space, the cell-centred quantities $m_{i,j}$ used in the definition of the numerical fluxes
\eqref{eq:q3d}-\eqref{eq:q3d1} are replaced by a high order polynomial interpolation through the so-called WENO reconstruction of order
three.

\subsubsection{A numerical method for the Navier-Stokes model}\label{NS1}
We briefly recall here the Navier-Stokes (NS) model which can obtained from the standard Chapmann-Enskog expansion \cite{cercignani69} in
the case of isentropic flows and is used for comparison purposes in the rest of the section. The NS model reads
\begin{equation}
	\begin{aligned}
		&\partial_t \rho+\nabla_{\bm{x}}\cdot(\rho \bm{u})=0,\\
		&\partial_t (\rho \bm{u})+\nabla_{\bm{x}}\cdot(\rho \bm{u} \otimes \bm{u}) + \nabla_{\bm{x}} p
		= \varepsilon RT \ \nabla_{\bm{x}}\cdot \left(\rho \Big( \nabla_{\bm{x}} \bm{u} + (\nabla_{\bm{x}}
		\bm{u})^T \Big)\right).
	\end{aligned}
	\label{NS}
\end{equation}
The numerical method follows the path of the discretization of the HMM model. Thus, rewriting the system \eqref{NS} as \eqref{Euler_eq},
where $\U{Q}, \U{F}, \U{G}$ have the same meaning of \eqref{eq:euler_sys} while the source term becomes
\begin{equation*}
\U{S}(\Q) =\varepsilon RT \begin{bmatrix} 0 \\ (\xi^{NS}_{xx})_x+(\xi^{NS}_{yx})_y\\
(\xi^{NS}_{xy})_x+(\xi^{NS}_{yy})_y\\
\end{bmatrix}.\label{eq:S}
\end{equation*}
with
\begin{eqnarray*}
	\xi^{NS}_{xx}=2\rho\frac{\partial u_x}{\partial x}, \ \xi^{NS}_{yy}=2\rho\frac{\partial u_y}{\partial y},\
	\xi^{NS}_{xy}=\xi^{NS}_{yx}=\rho\left(\frac{\partial u_x}{\partial y}+\frac{\partial u_y}{\partial x}\right).
	\label{eqn.tauNS}
\end{eqnarray*}
We then proceed, as done previously, using a first order explicit Euler scheme in time and the same Rusanov flux with WENO reconstruction
for the hyperbolic fluxes and the same second order discrete velocity gradients \eqref{eqn.face_grad} are used for the viscous terms. This
ends the description of the model and of the numerical scheme.

\subsection{Highly oscillating fluid}\label{oscill}
\begin{figure}[!ht]
	\centering
	\includegraphics[width=0.45\textwidth]{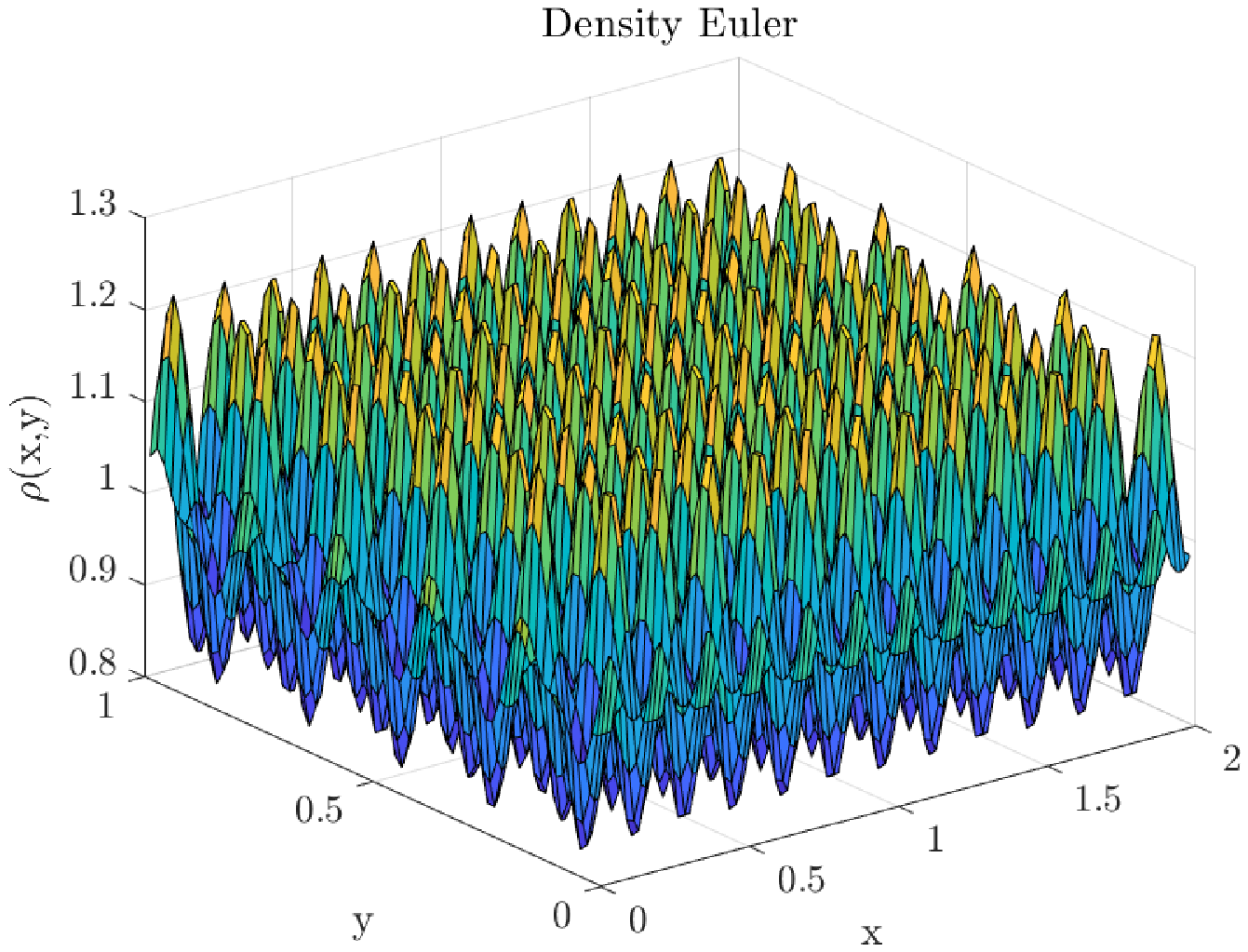}
	\includegraphics[width=0.45\textwidth]{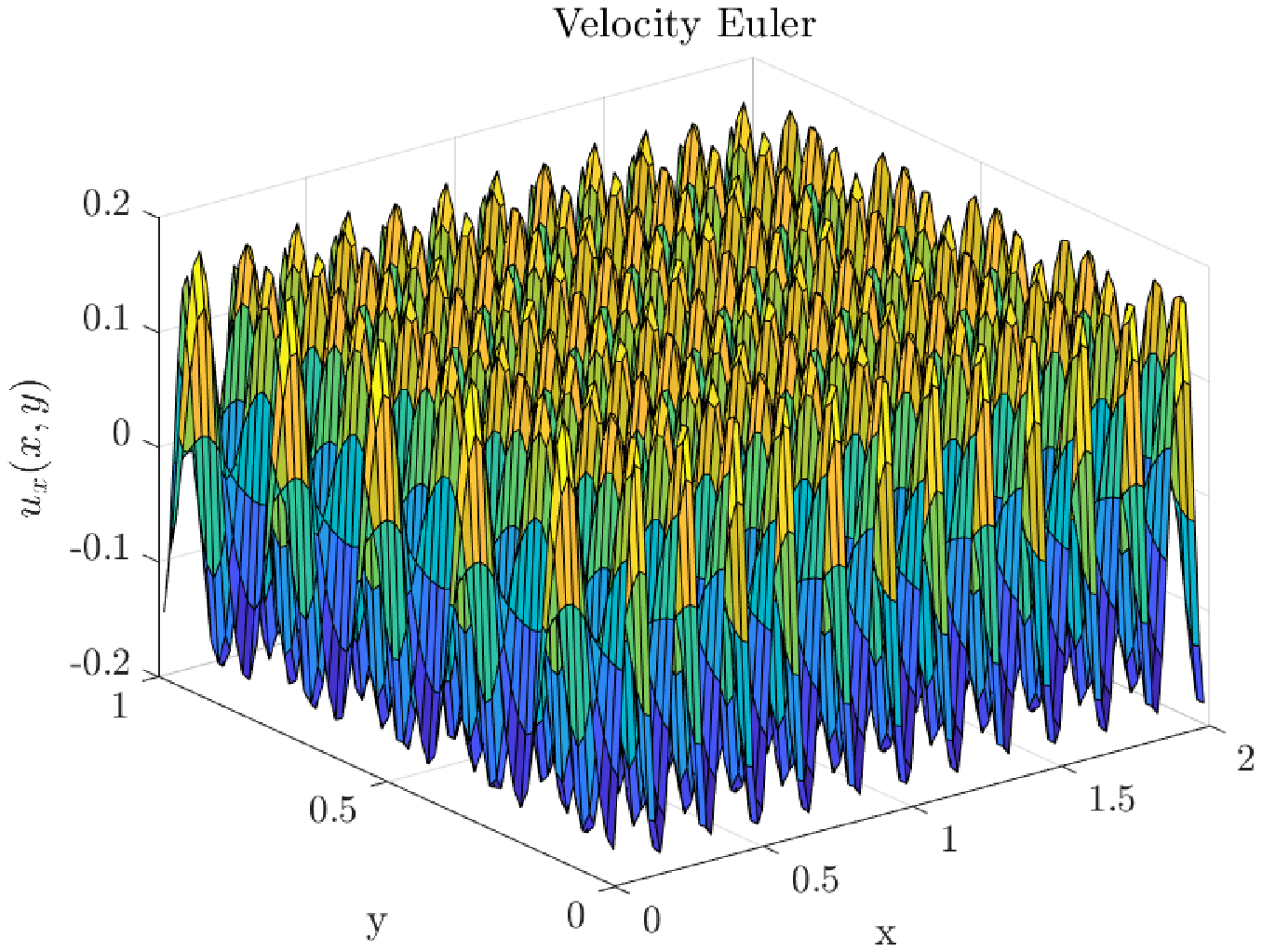}
	\caption{\textbf{Highly oscillating fluid}. Profiles of density (left) and  first component of the velocity $u_x(x,y,t)$ (right) at
	final time for the isentropic Euler equations. 	\label{fig:oscill1}}
\end{figure}

\begin{figure}[!ht]
	\centering
	\includegraphics[width=0.45\textwidth]{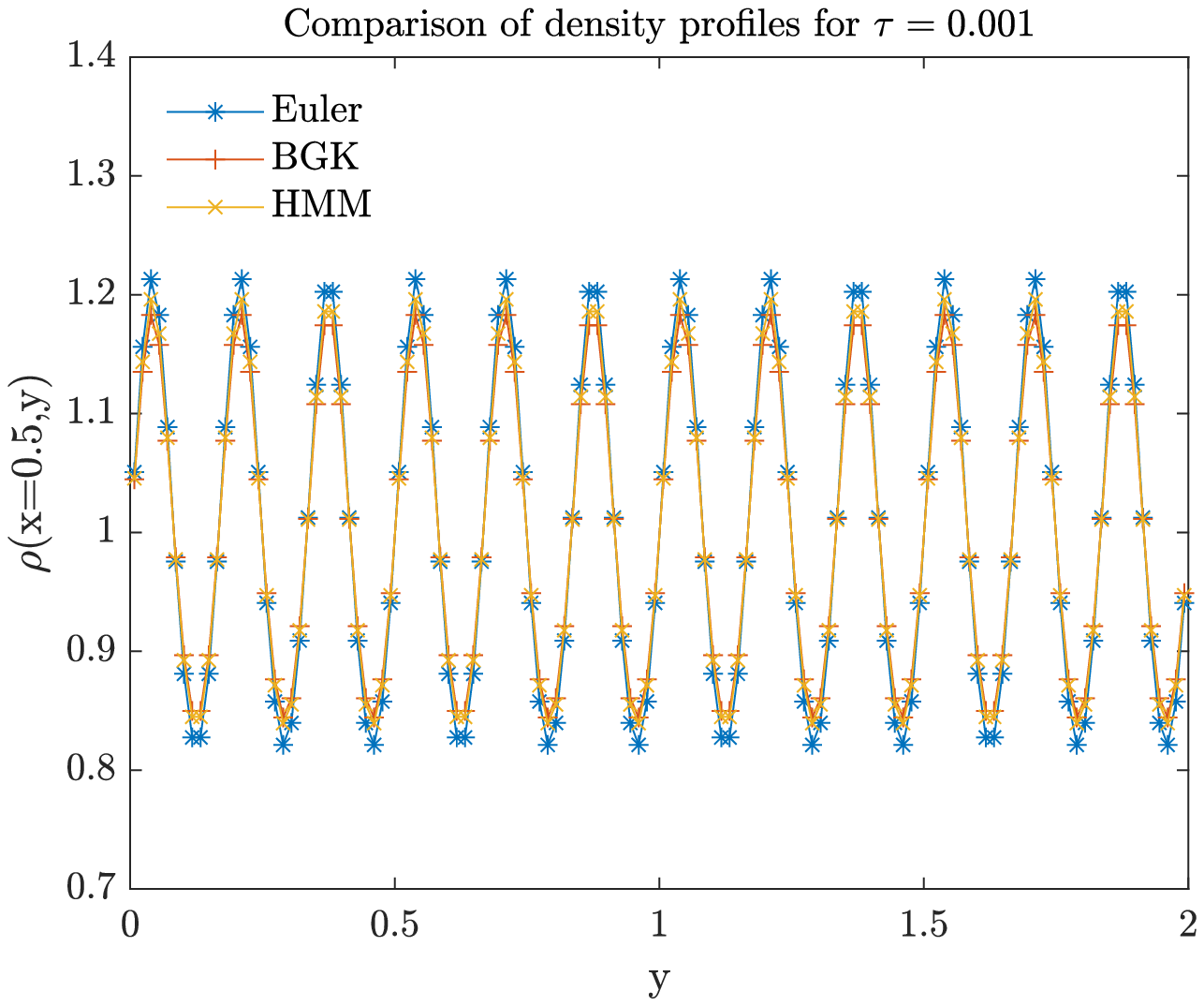}
	\includegraphics[width=0.45\textwidth]{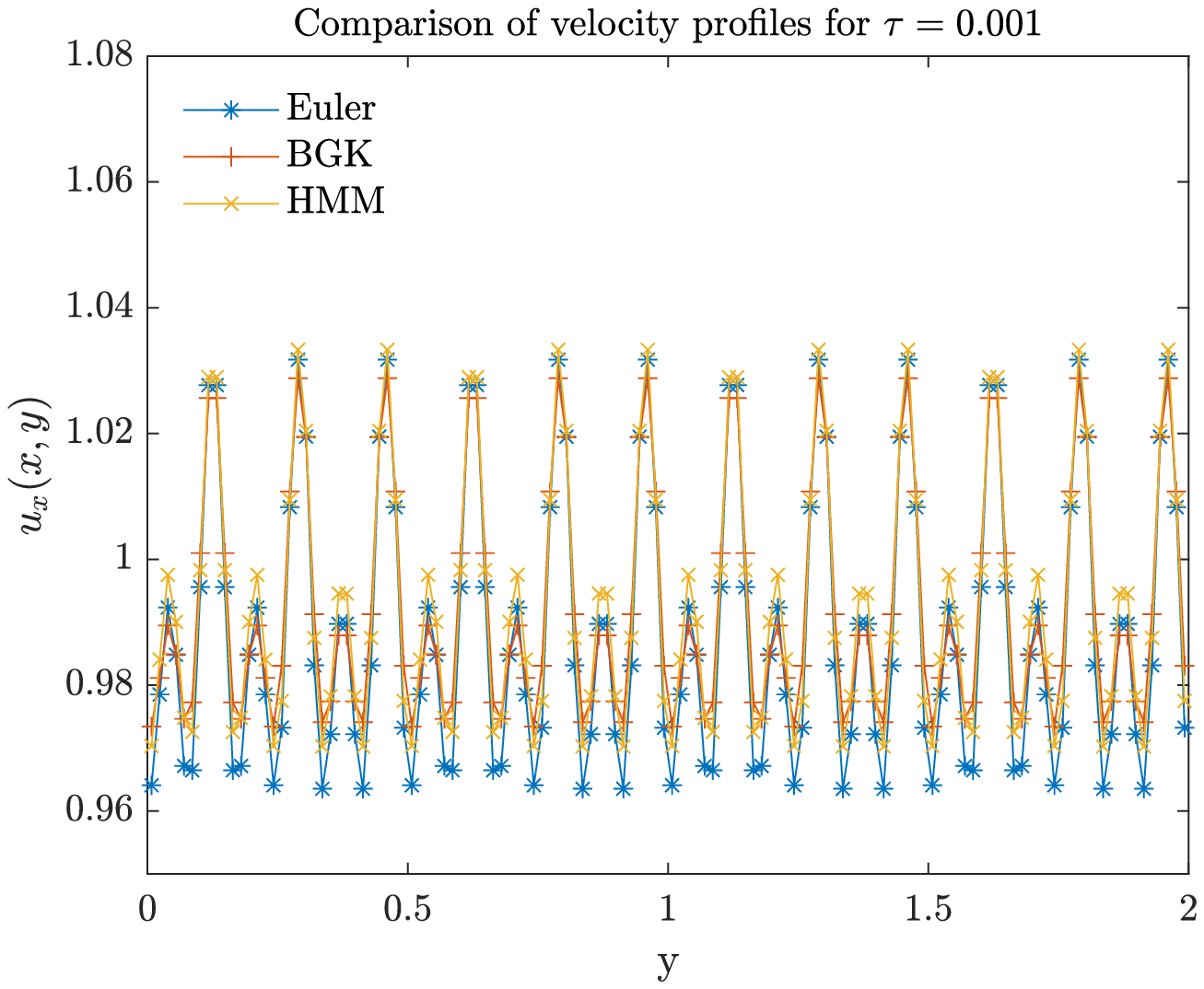}\\
	\includegraphics[width=0.45\textwidth]{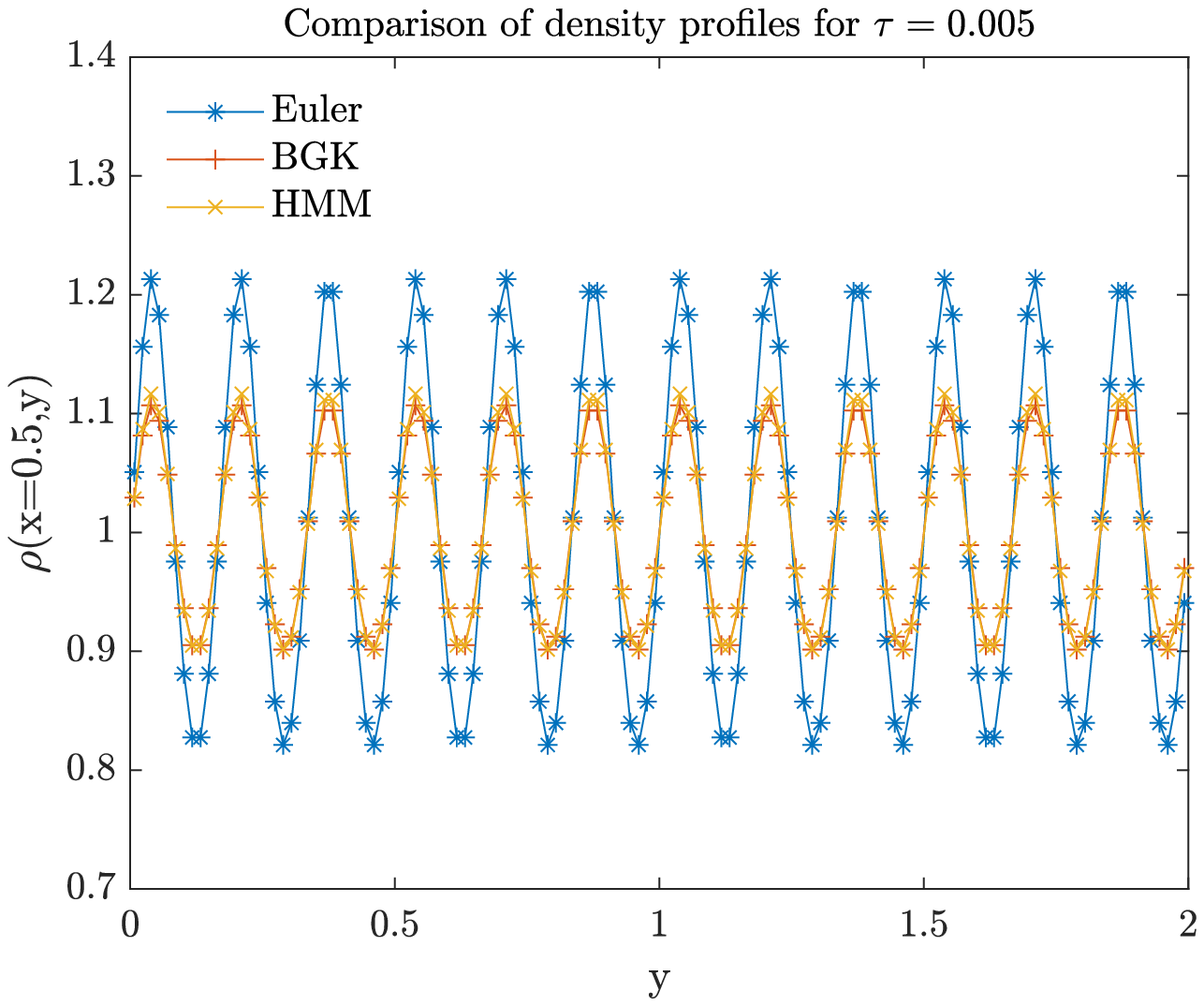}
	\includegraphics[width=0.45\textwidth]{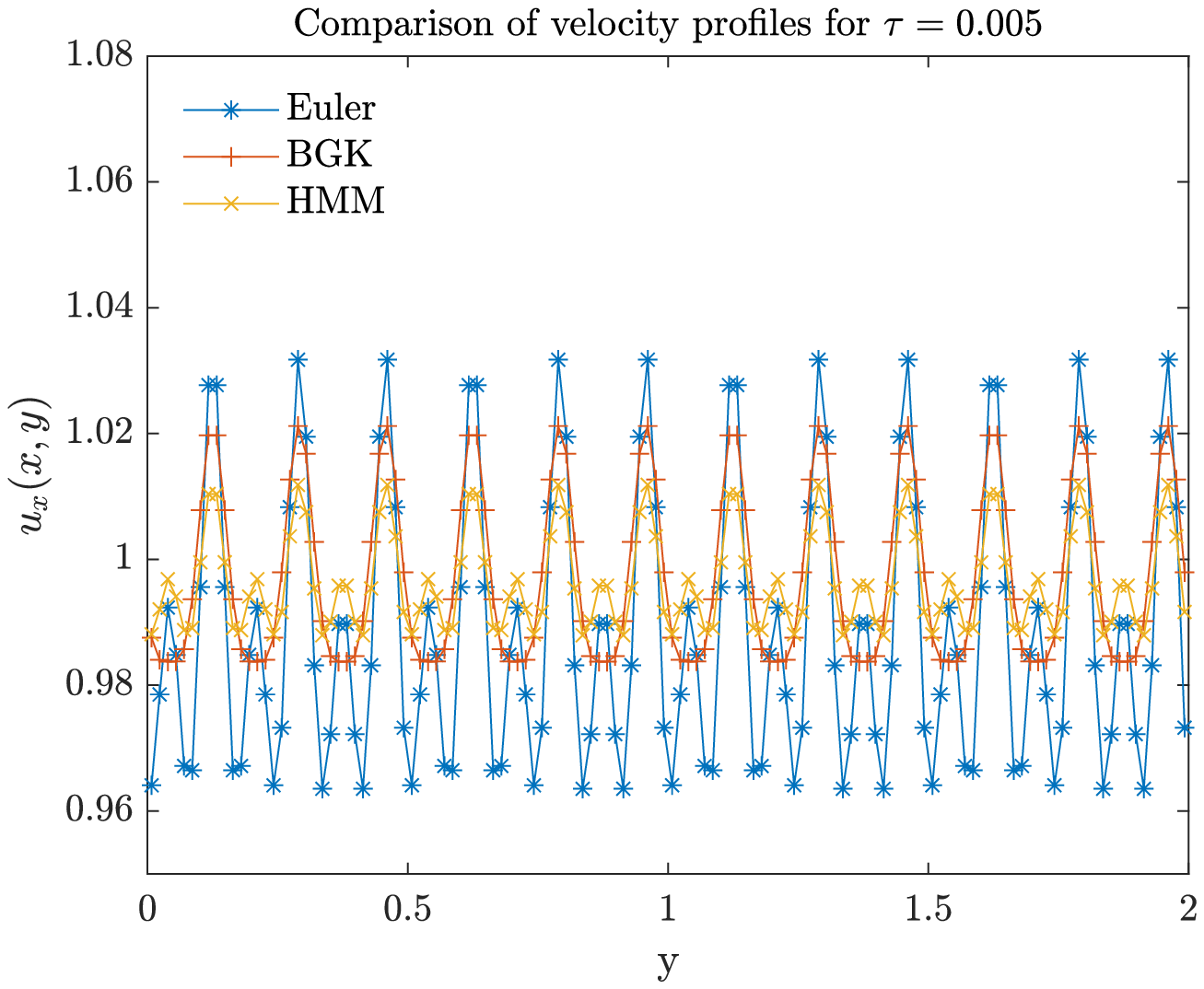}\\
	\includegraphics[width=0.45\textwidth]{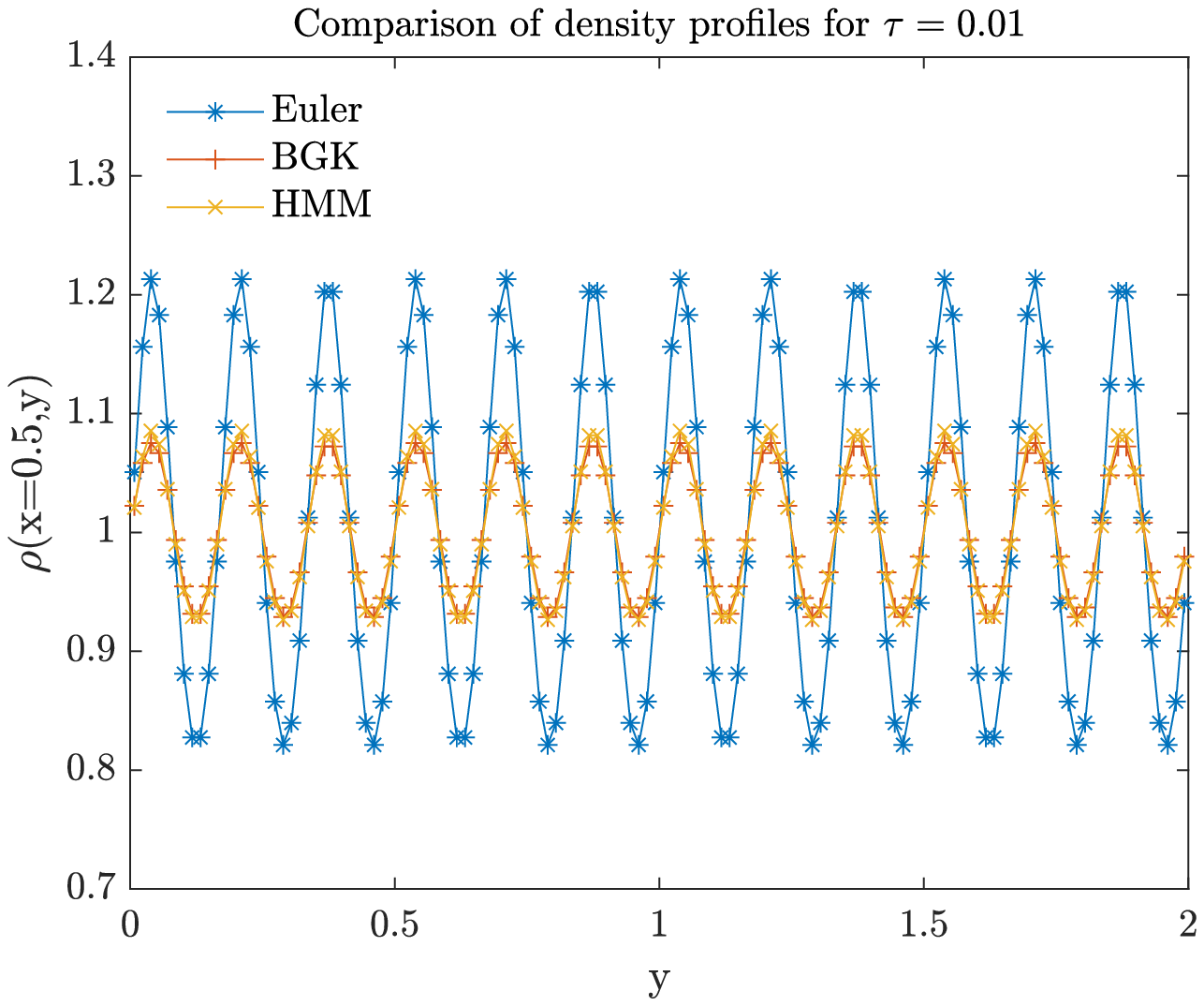}
	\includegraphics[width=0.45\textwidth]{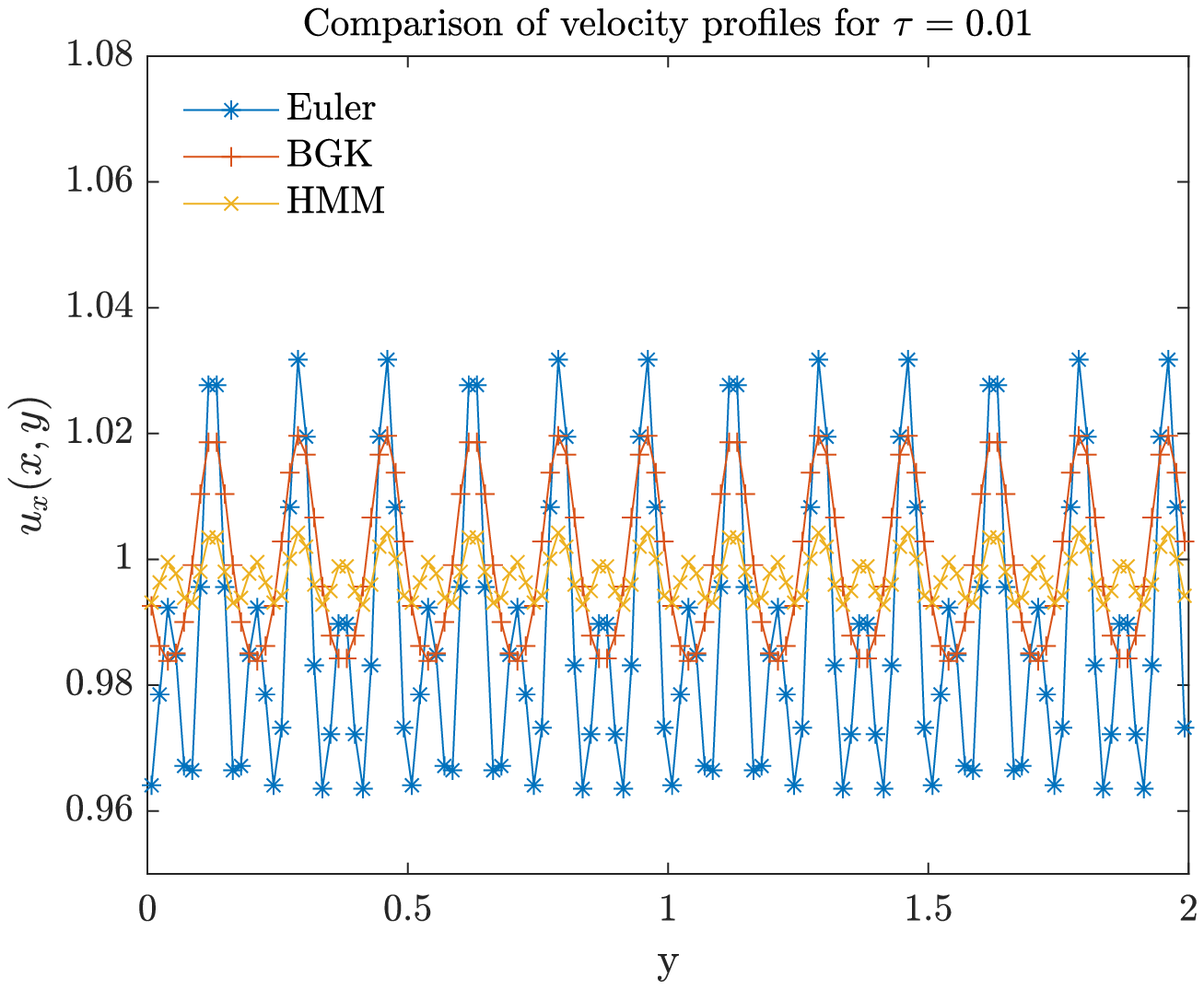}\\
	\caption{\textbf{Highly oscillating fluid}. Comparison of the density (left) and velocity $u_x(x,y,t)$ (right) profiles for $x=0.5$. The
	results for the Euler, the BGK and the HMM equations with respectively $\tau=0.001$, $\tau=0.005$ and $\tau=0.01$ are shown.	\label{fig:oscill2}}
\end{figure}

\begin{figure}[!ht]
	\centering
	\includegraphics[width=0.45\textwidth]{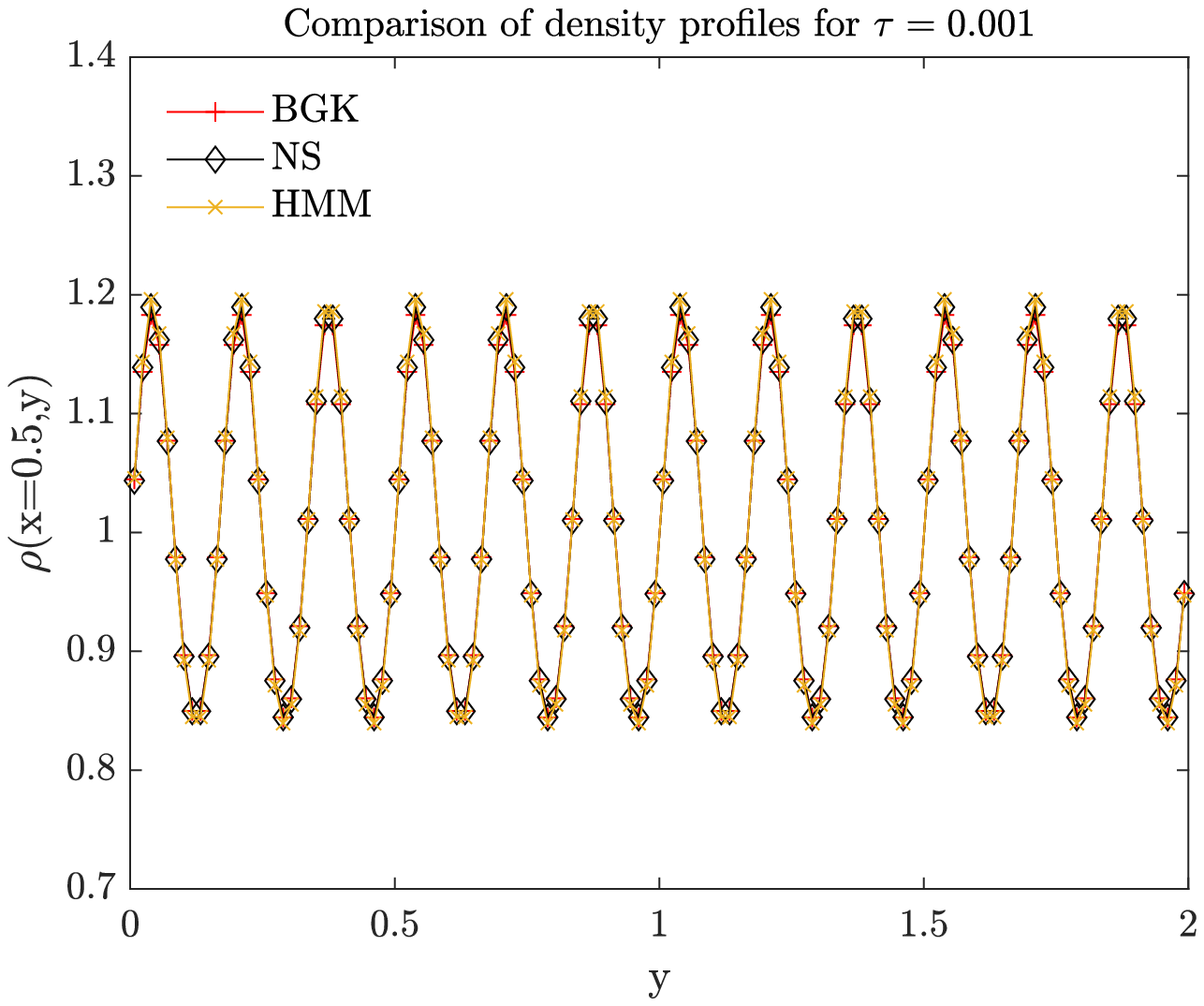}
	\includegraphics[width=0.45\textwidth]{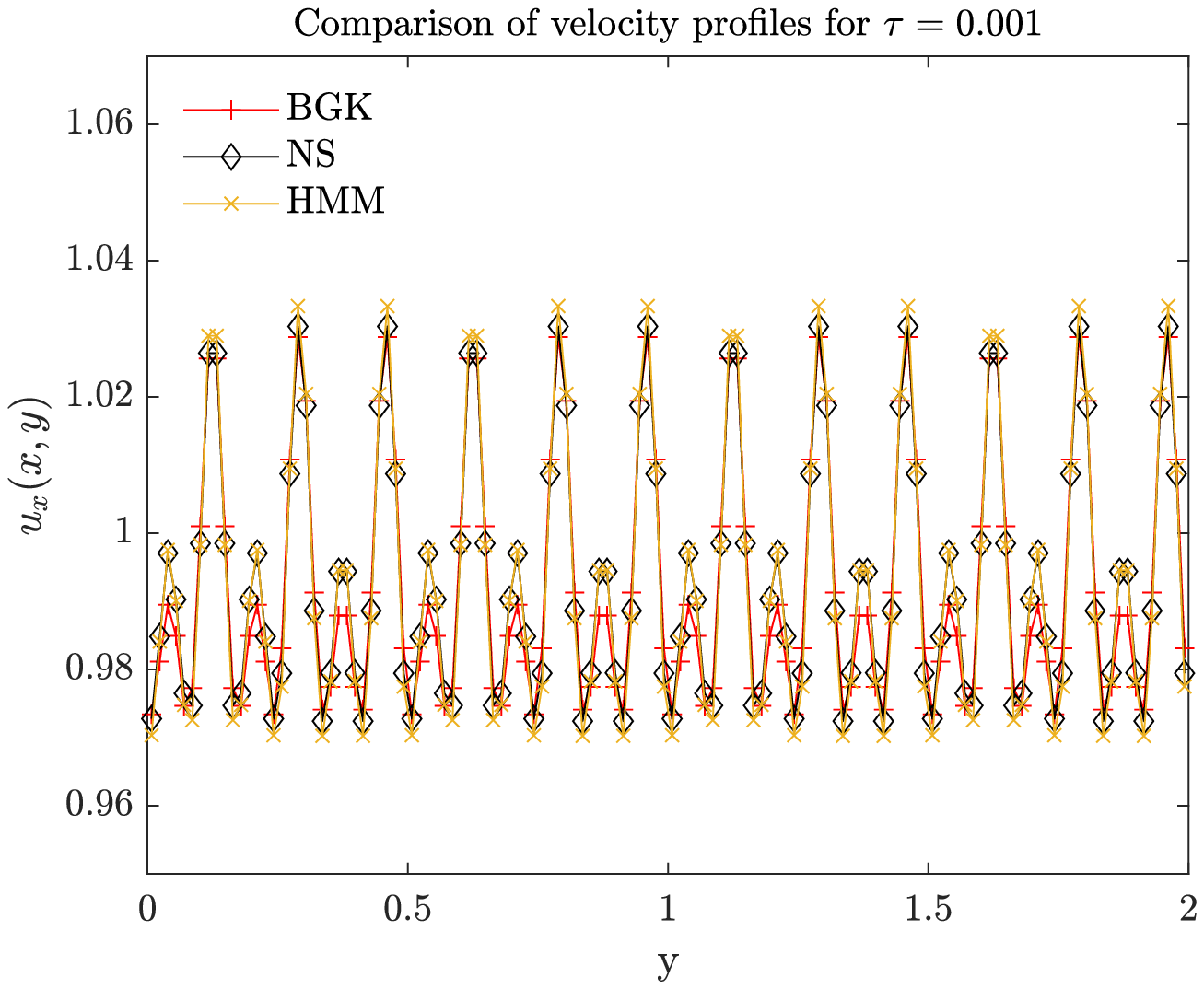}\\
	\includegraphics[width=0.45\textwidth]{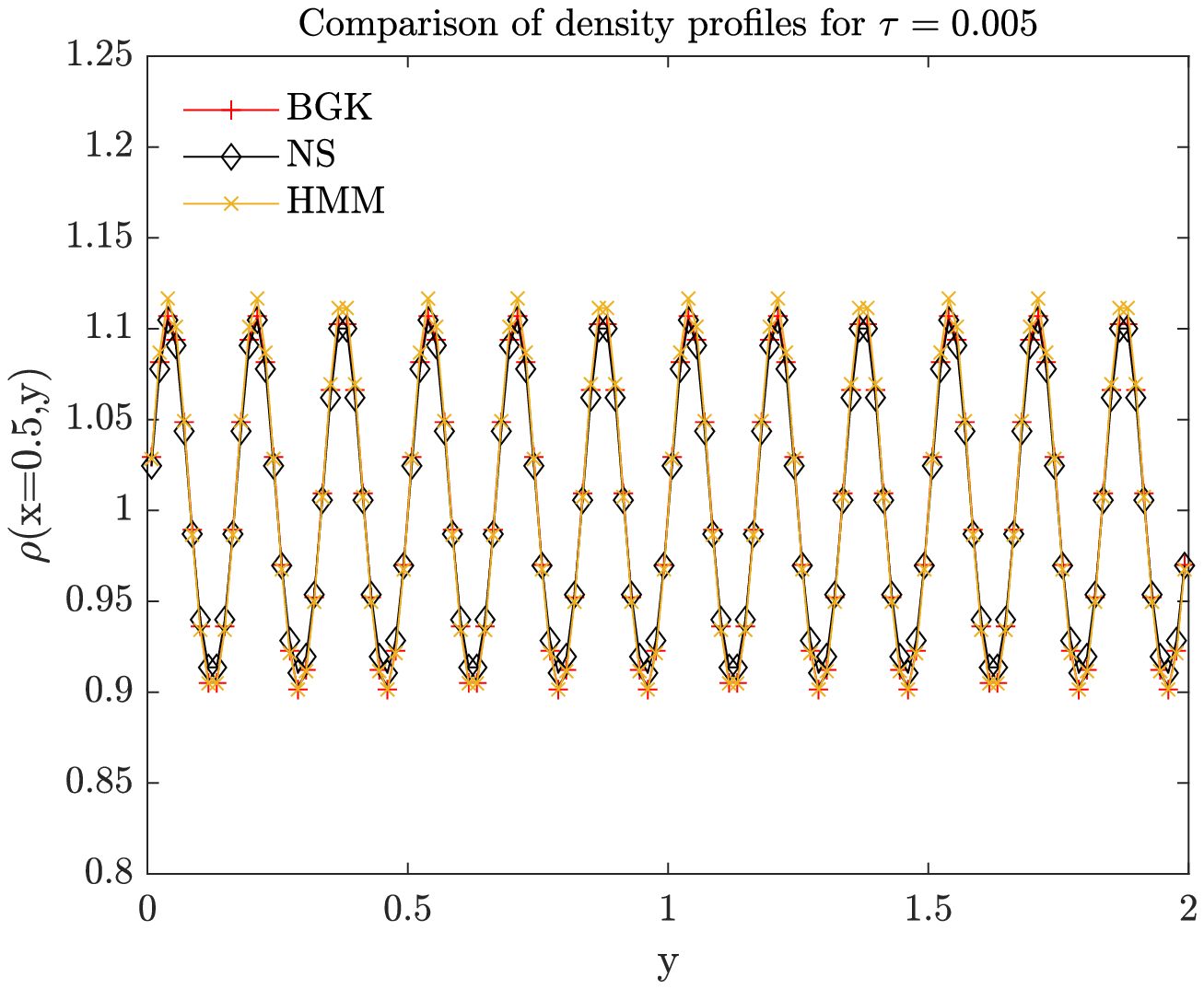}
	\includegraphics[width=0.45\textwidth]{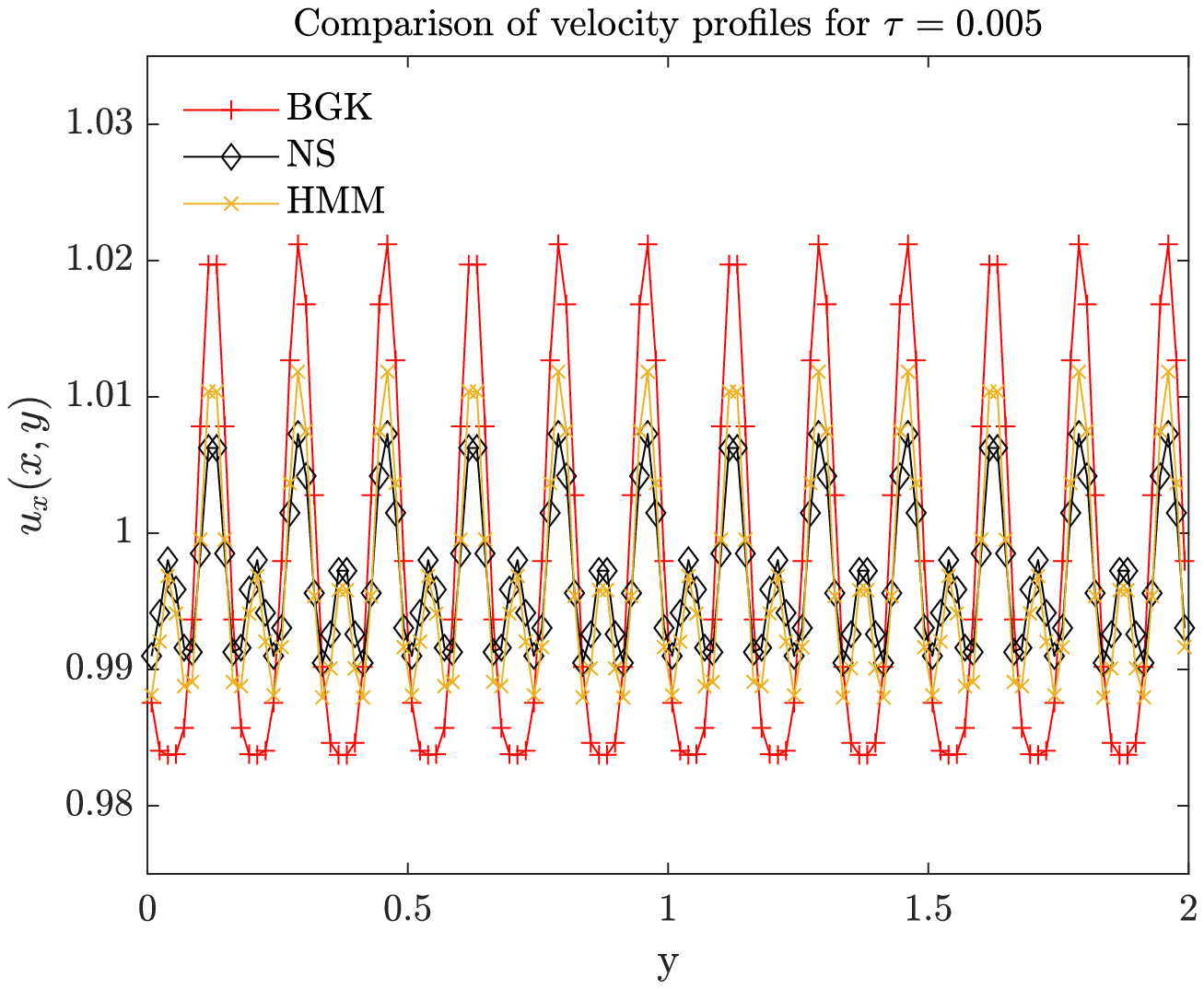}\\
	\includegraphics[width=0.45\textwidth]{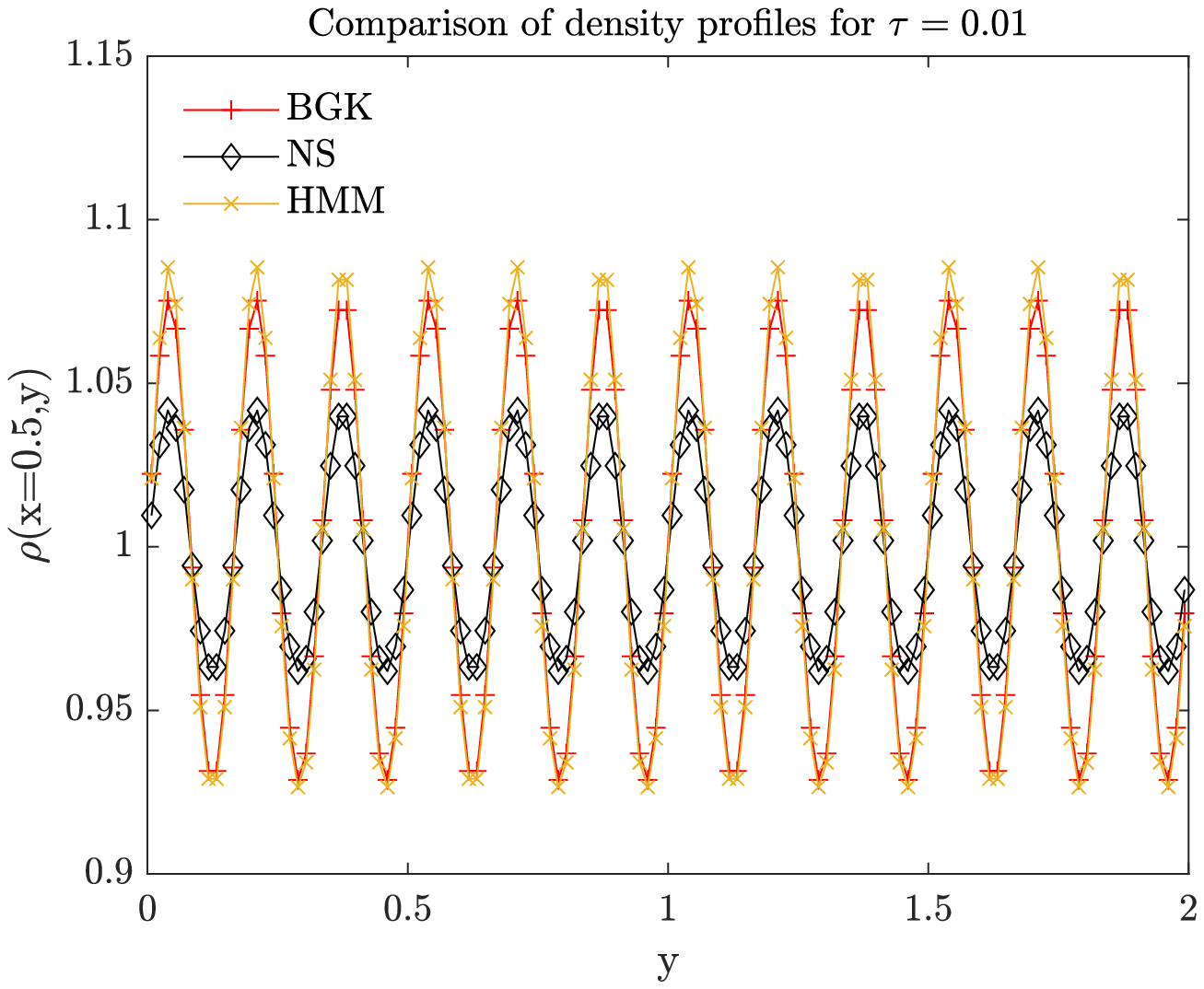}
	\includegraphics[width=0.45\textwidth]{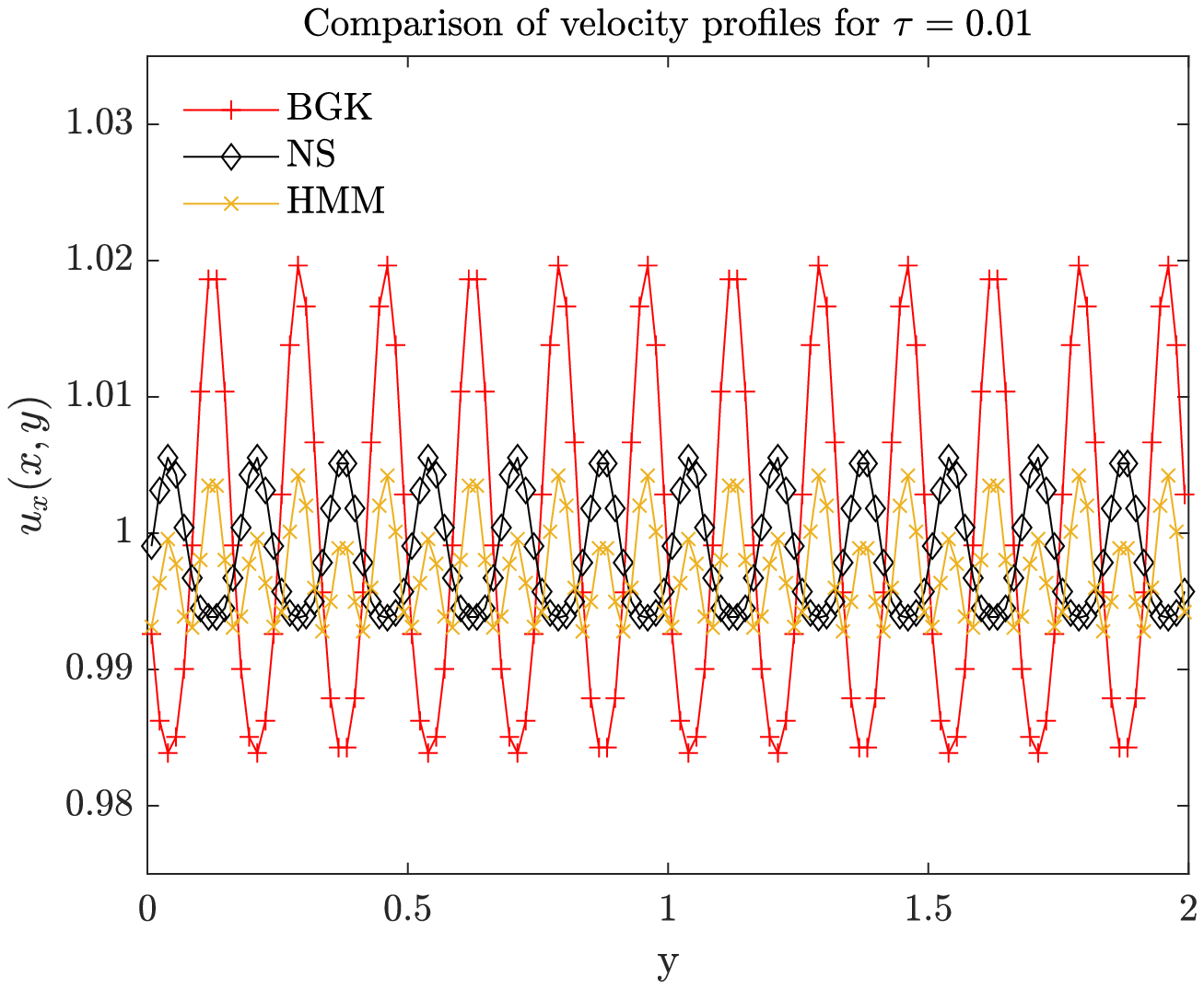}\\
	\caption{\textbf{Highly oscillating fluid}. Comparison of the density (left) and velocity $u_x(x,y,t)$ (right) profiles for $x=0.5$.
	The results for the the BGK, the Navier-Stokes and the HMM equations with respectively $\tau=0.001$, $\tau=0.005$ and $\tau=0.01$ are
	shown.\label{fig:oscill3}}
\end{figure}

We start by considering a highly oscillating initial condition for a flow in the computational domain $\Omega=[0,1]\times [0,2]$ and
assume the periodic boundary conditions to be imposed in both directions. The temperature $T$ as well as the constant $R$ are set to one.
The domain is paved with $N_x\times N_y=64\times 128$ cells in the physical space while for what concerns the BGK model we discretize the
velocity space with $20\times 20$ cells with $v_{x,max}=v_{y,max}=5$ and $v_{x,min}=v_{y,min}=-5$, where $\bm{v}=(v_x,v_y)$. The
initial distribution function is assumed to be at local equilibrium, i.e. $ f(\bm{x},\bm{v},0)=M_{\rho,\bm{u}}(\bm{x},\bm{v},0)$. The initial
density and velocities are
\begin{equation*}
	\label{initial_oscill}
	\rho(x,y,0)=1+0.2\cos(10\pi x)\sin(12\pi y),\quad \bm{u}(x,y,0)=(u_x,u_y)=(1,0).
\end{equation*}
In Figure \ref{fig:oscill1}, we plot the highly oscillating density and the velocity profiles after $100$ time iterations when the
isentropic Euler model is used. In Figure \ref{fig:oscill2}, we present the same result at a fixed value of $x=0.5$. In particular, we
compare the results obtained for the isentropic Euler, the BGK and the HMM models using different values of the scaling parameter $\tau$.
In the BGK model we choose $\tau=0.001$, $\tau=0.005$ and $\tau=0.01$, while in the HMM we set $\tau_{HMM}=\tau/3$. In Figure
\ref{fig:oscill2} we can clearly see that oscillations are damped as fast as the scaling parameter becomes larger. In these three tested
situations, the HMM seems to be able to describe a kinetic regime both for the case of the density as well as for the case of the mean
velocity $u_x(x,y,t)$. For sake of comparison, in Figure \ref{fig:oscill3}, we plot the same density and velocity profiles at fixed
$x=0.5$ where we added the results of the Navier-Stokes model for $\varepsilon=\tau$. As one can see, for small values of the damping
parameters both NS and HMM give results very close to those of the BGK model, while for $\tau=0.01$ the NS model tends to overestimate the
damping for the density and fails in describing the velocity profile. Instead, the HMM model provides a very good description of the
density and while overdamping the velocity still is able to follow the behavior of the BGK model.

\subsection{Perturbed Couette flow}\label{Couette}
\begin{figure}[!ht]
	\centering
	\includegraphics[width=0.3\textwidth]{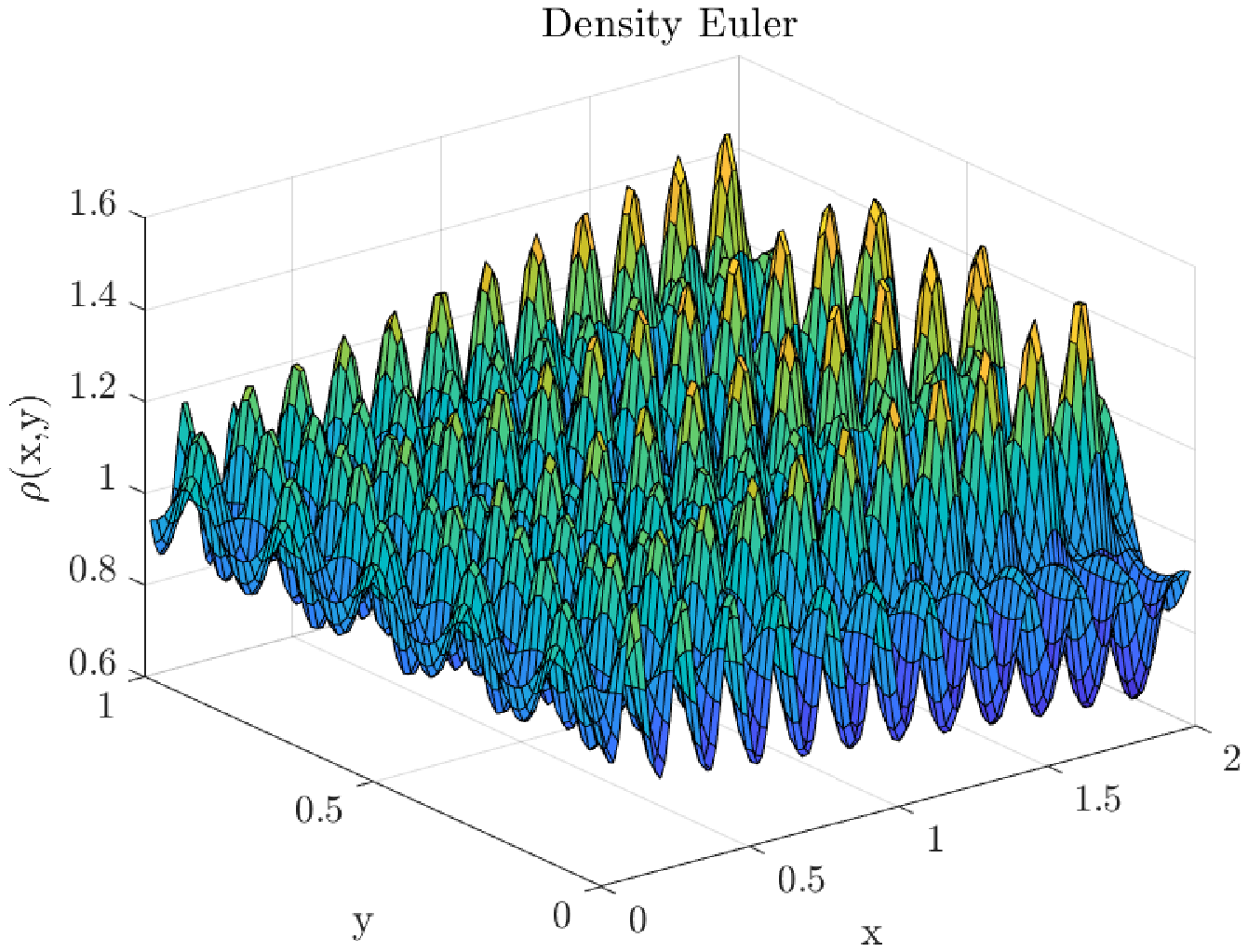}
	\includegraphics[width=0.3\textwidth]{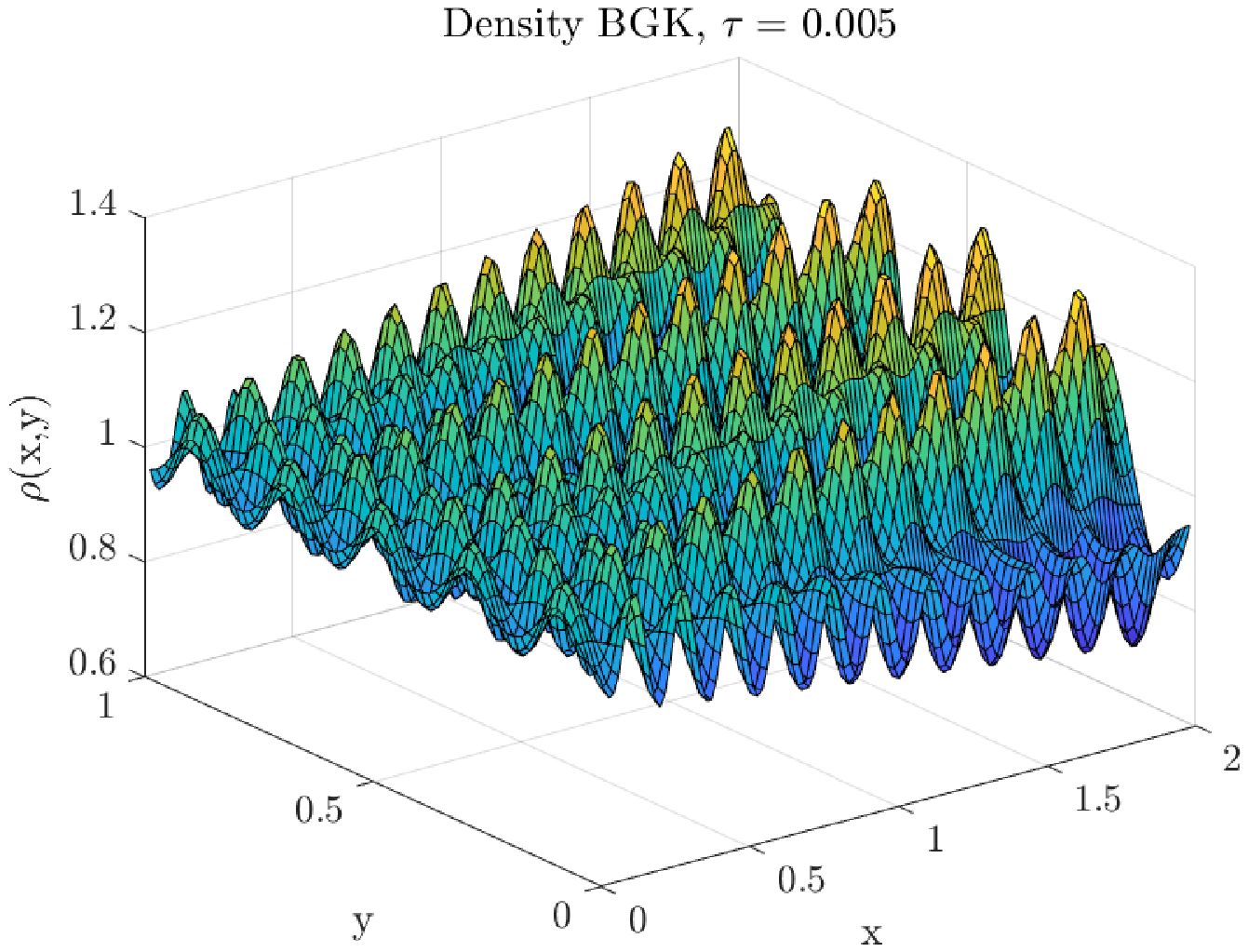}
	\includegraphics[width=0.3\textwidth]{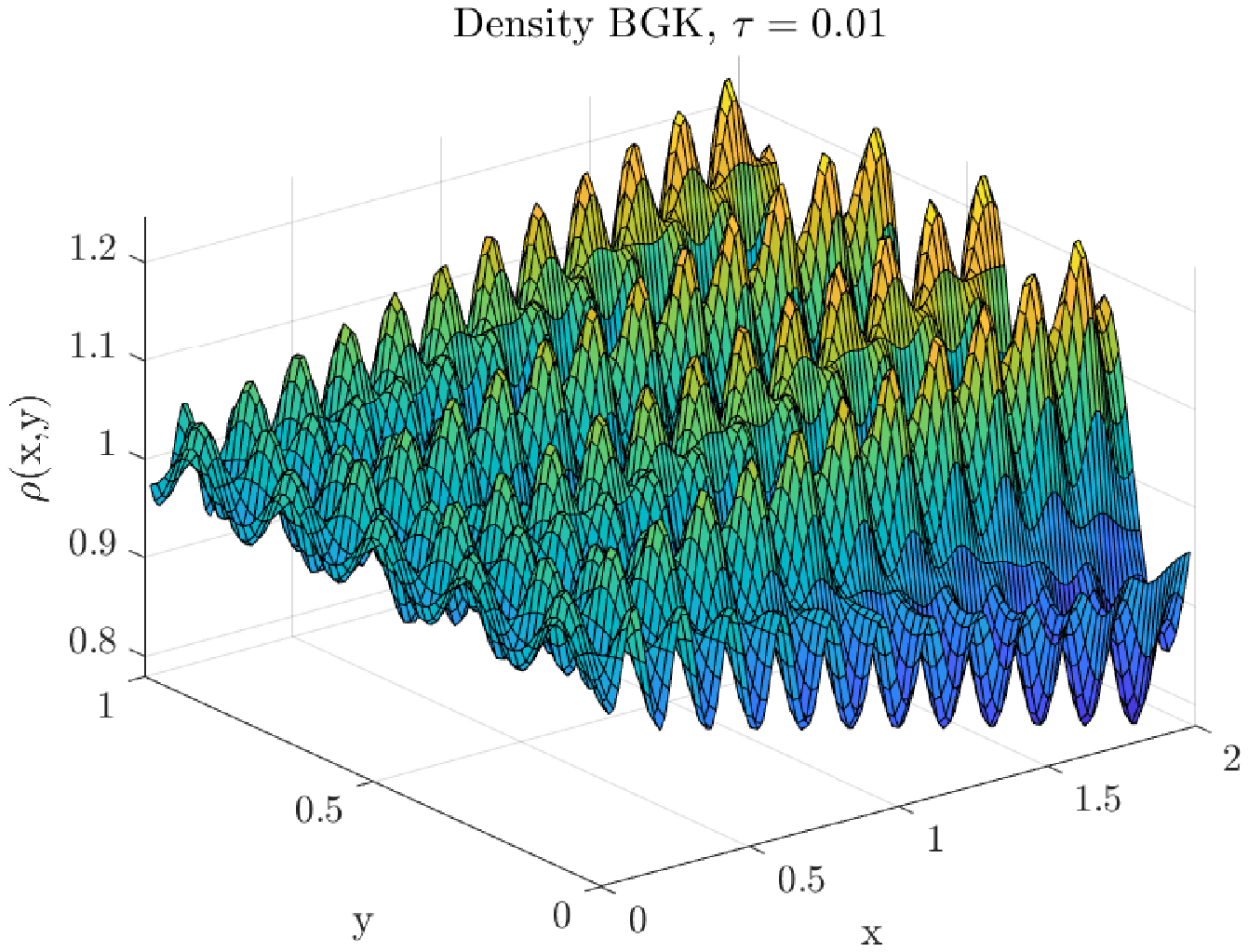}\\
	\includegraphics[width=0.3\textwidth]{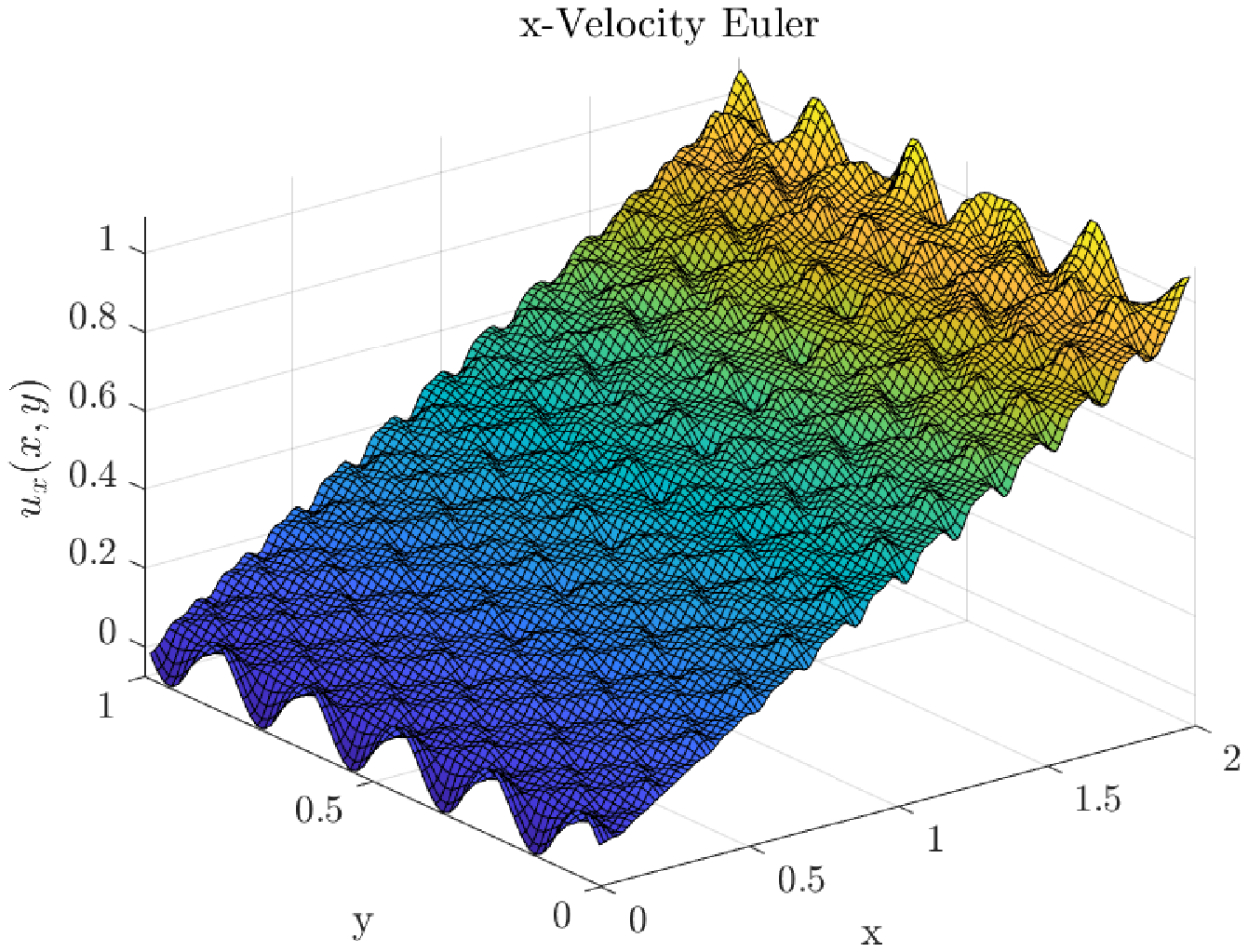}
	\includegraphics[width=0.3\textwidth]{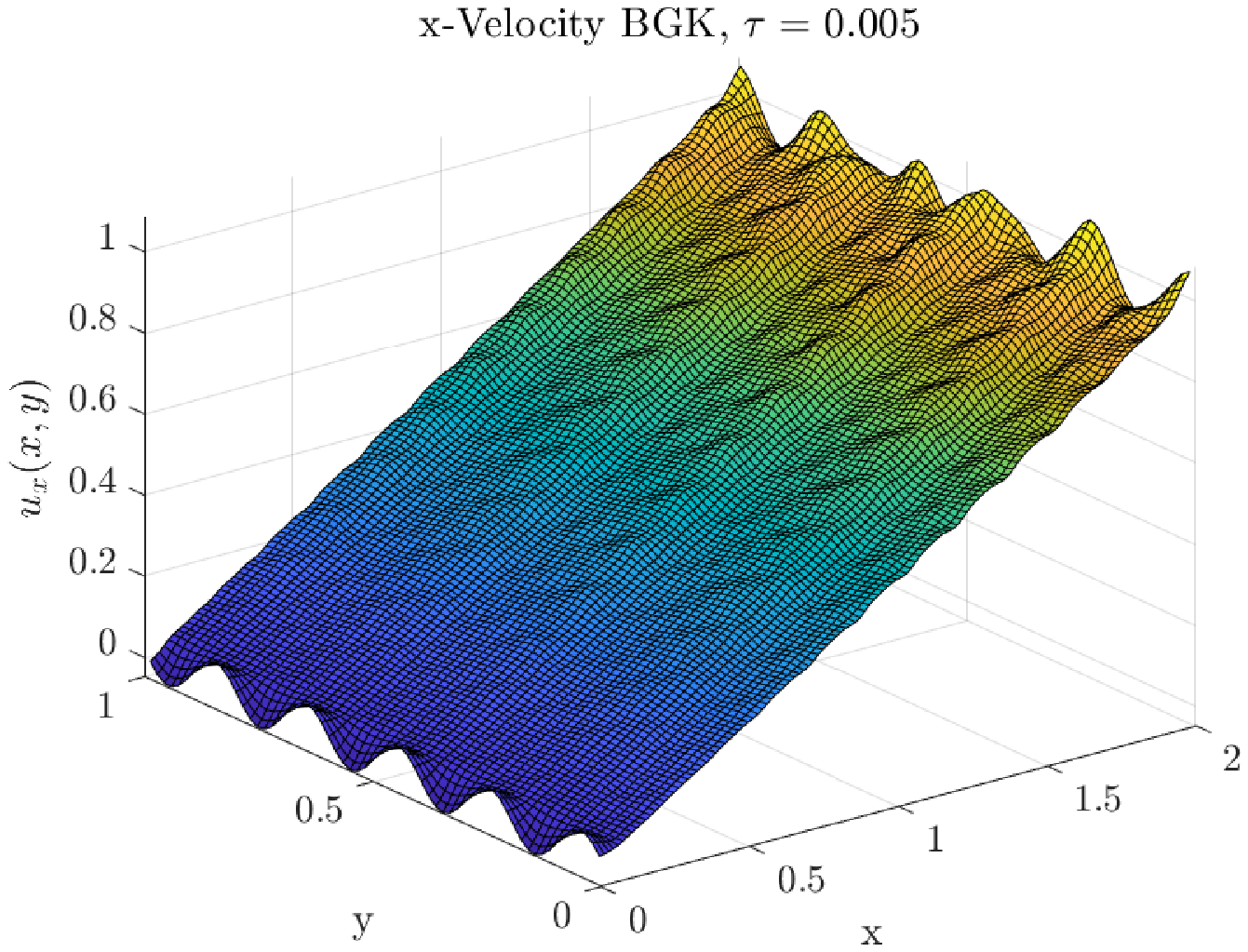}
	\includegraphics[width=0.3\textwidth]{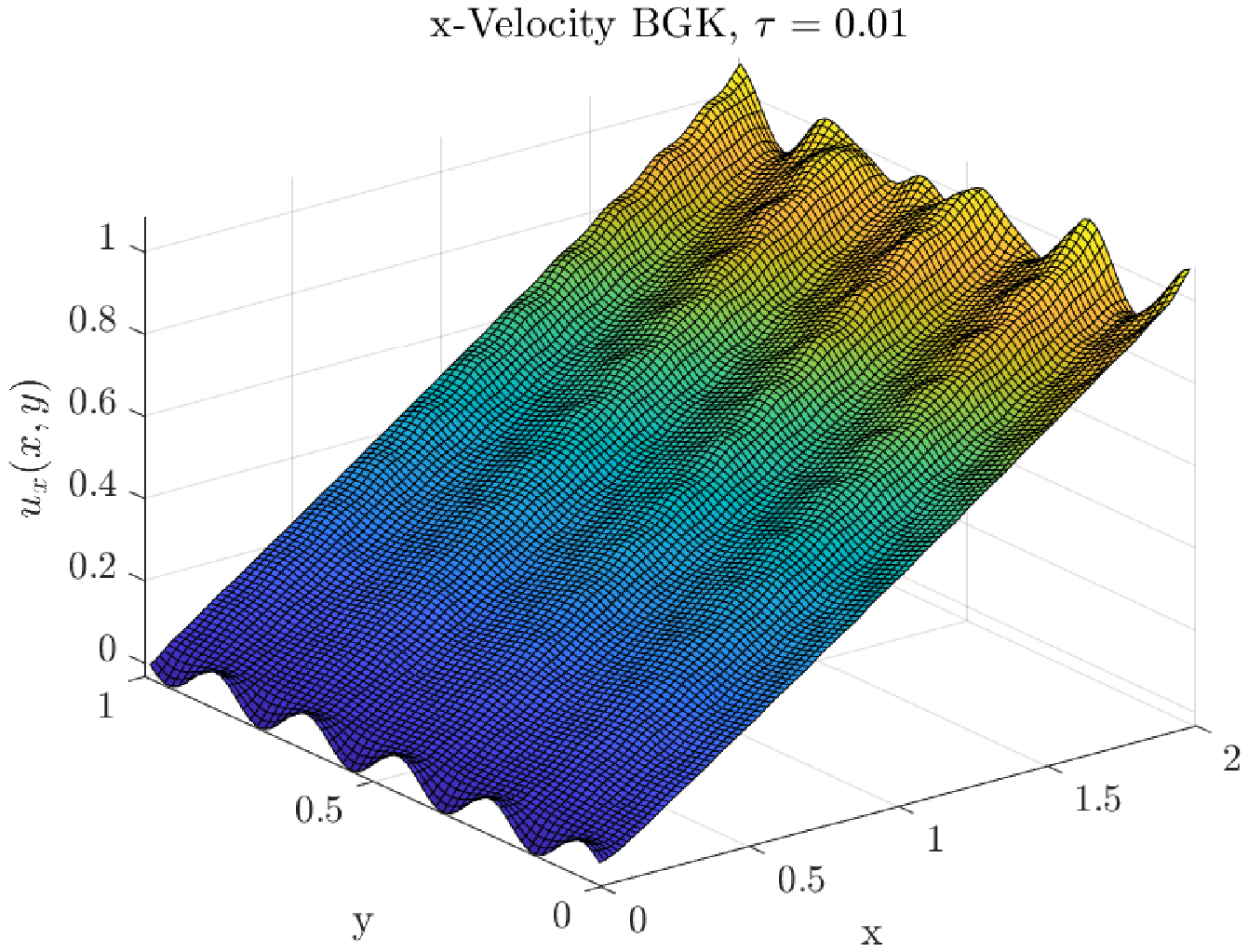}\\	
	\includegraphics[width=0.3\textwidth]{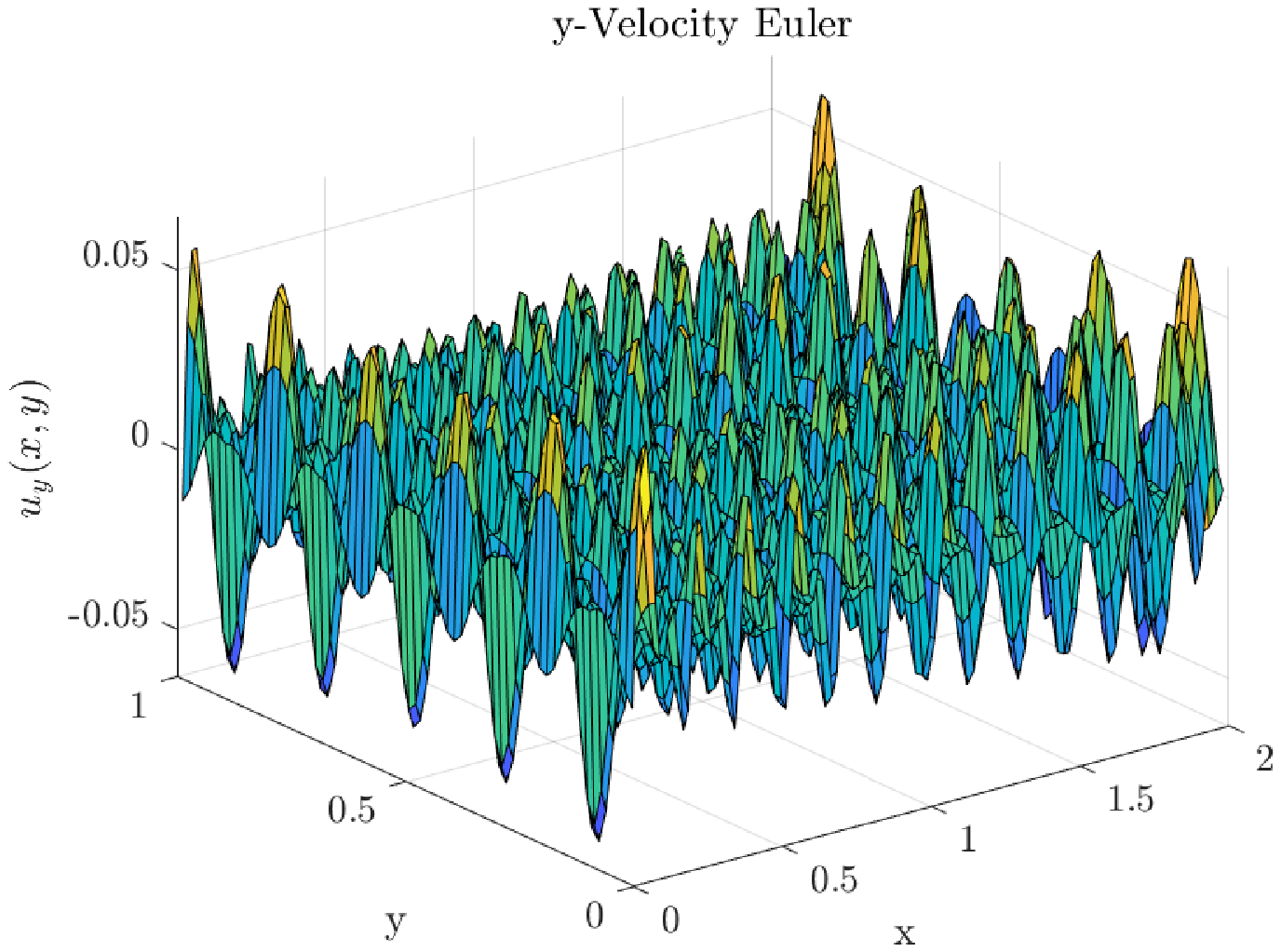}
	\includegraphics[width=0.3\textwidth]{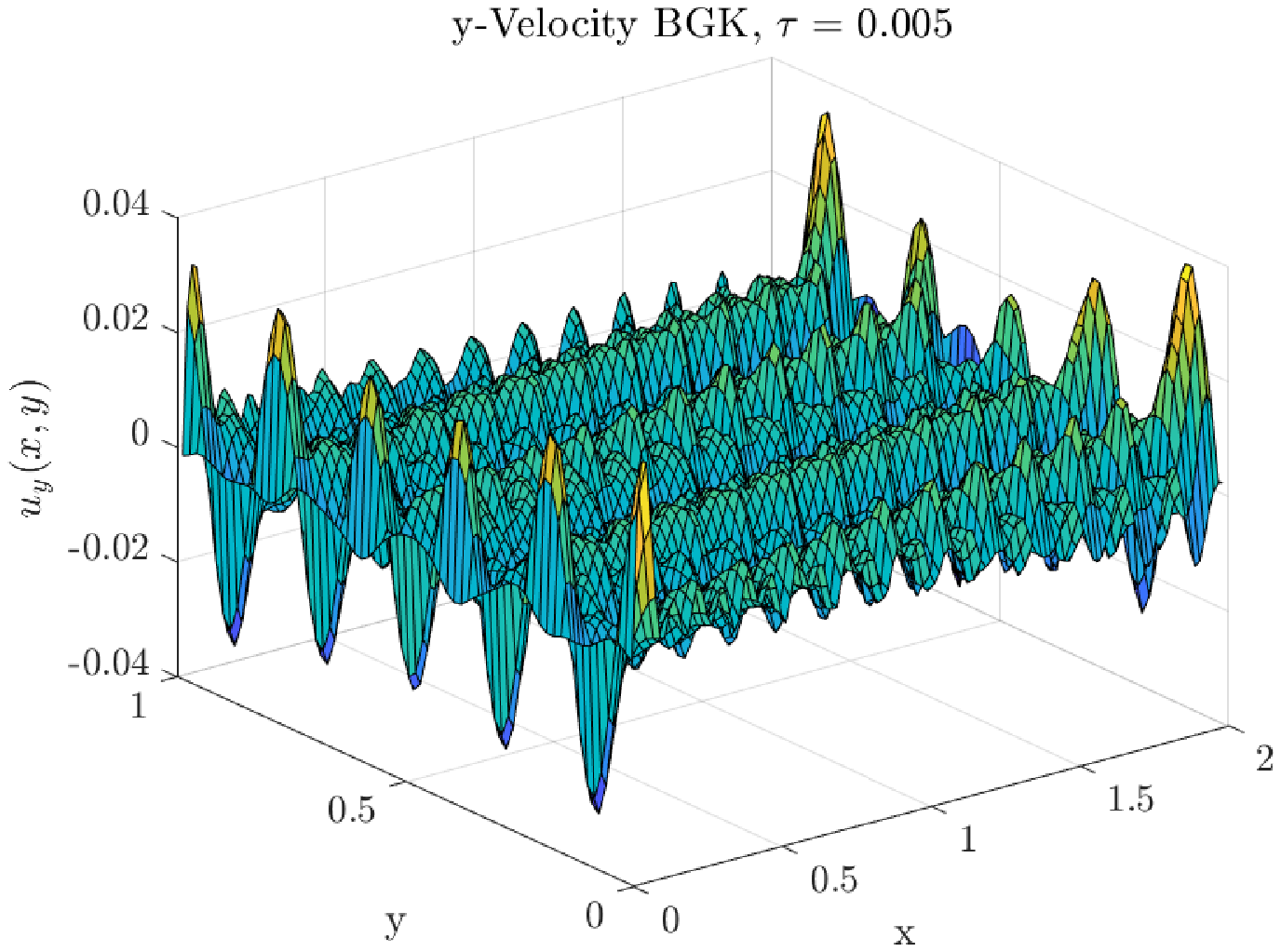}
	\includegraphics[width=0.3\textwidth]{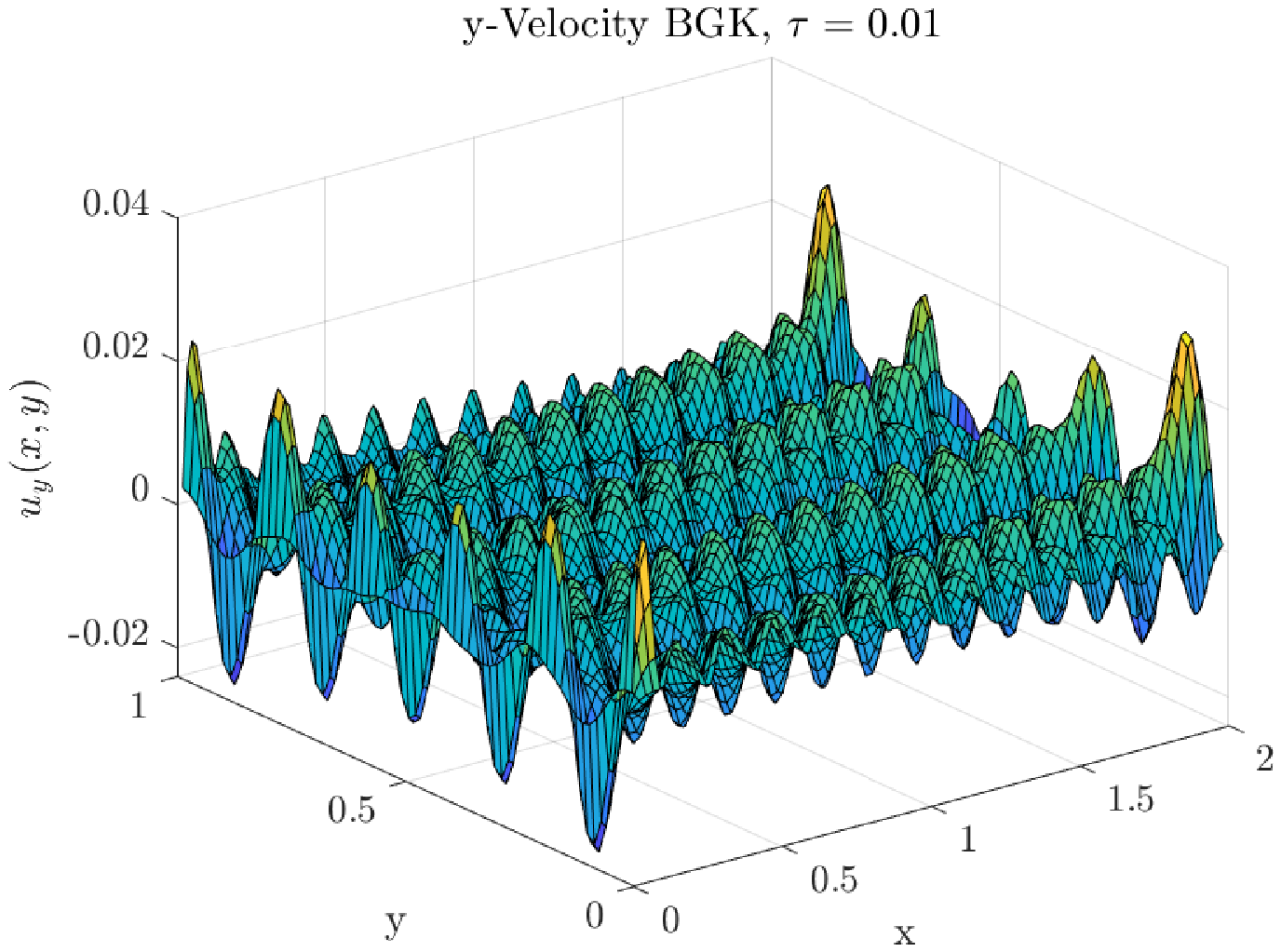}
	\caption{\textbf{Perturbed Couette flow}. Profiles of density and velocity for the Euler and the BGK equations. From top to bottom:
	Density, x-velocity and y-velocity profiles. Left $\tau=0$, middle $\tau=0.005$, right $\tau=0.01$. \label{fig:couette1}}
\end{figure}

\begin{figure}[!ht]
	\centering
	\includegraphics[width=0.3\textwidth]{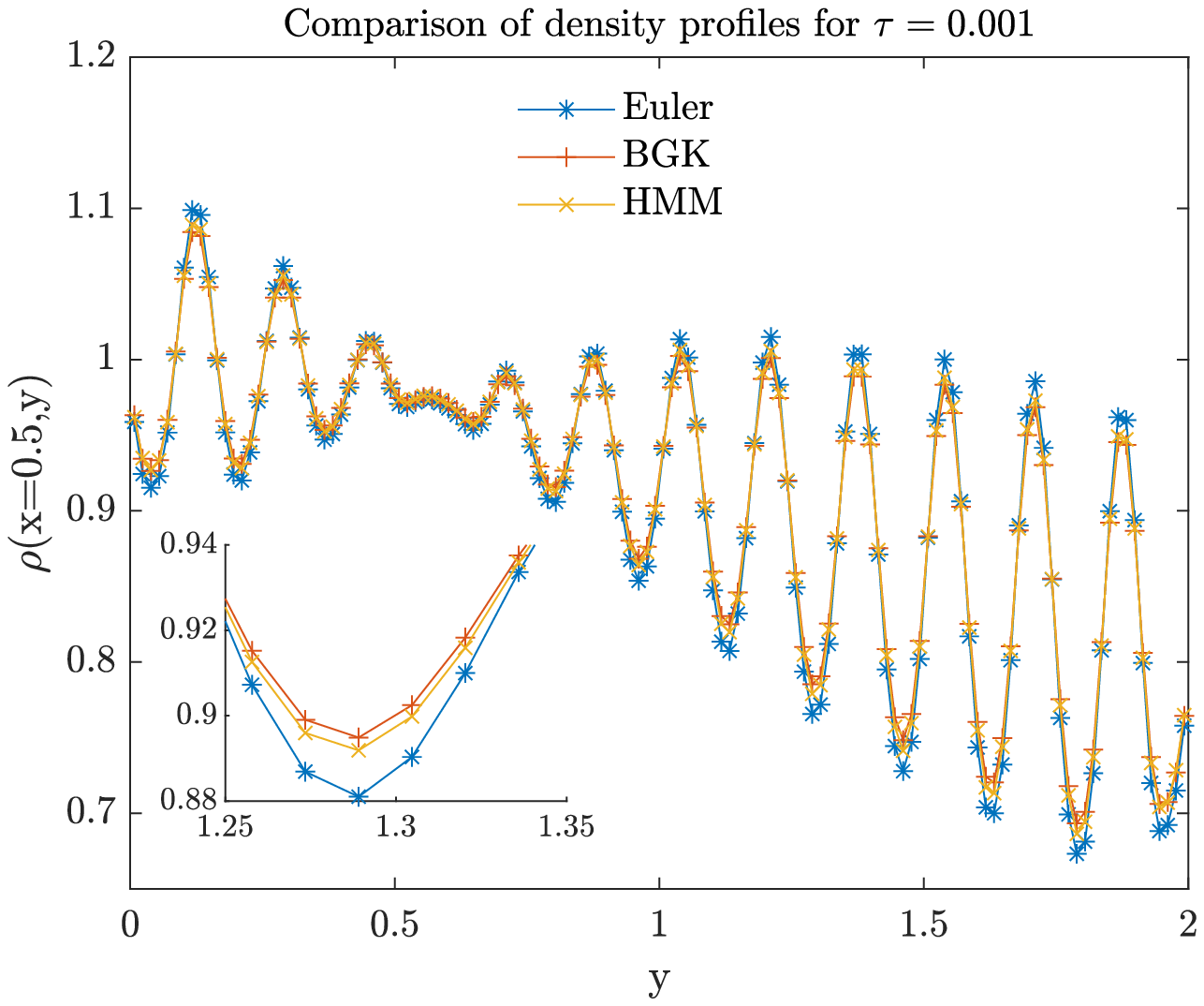}
	\includegraphics[width=0.3\textwidth]{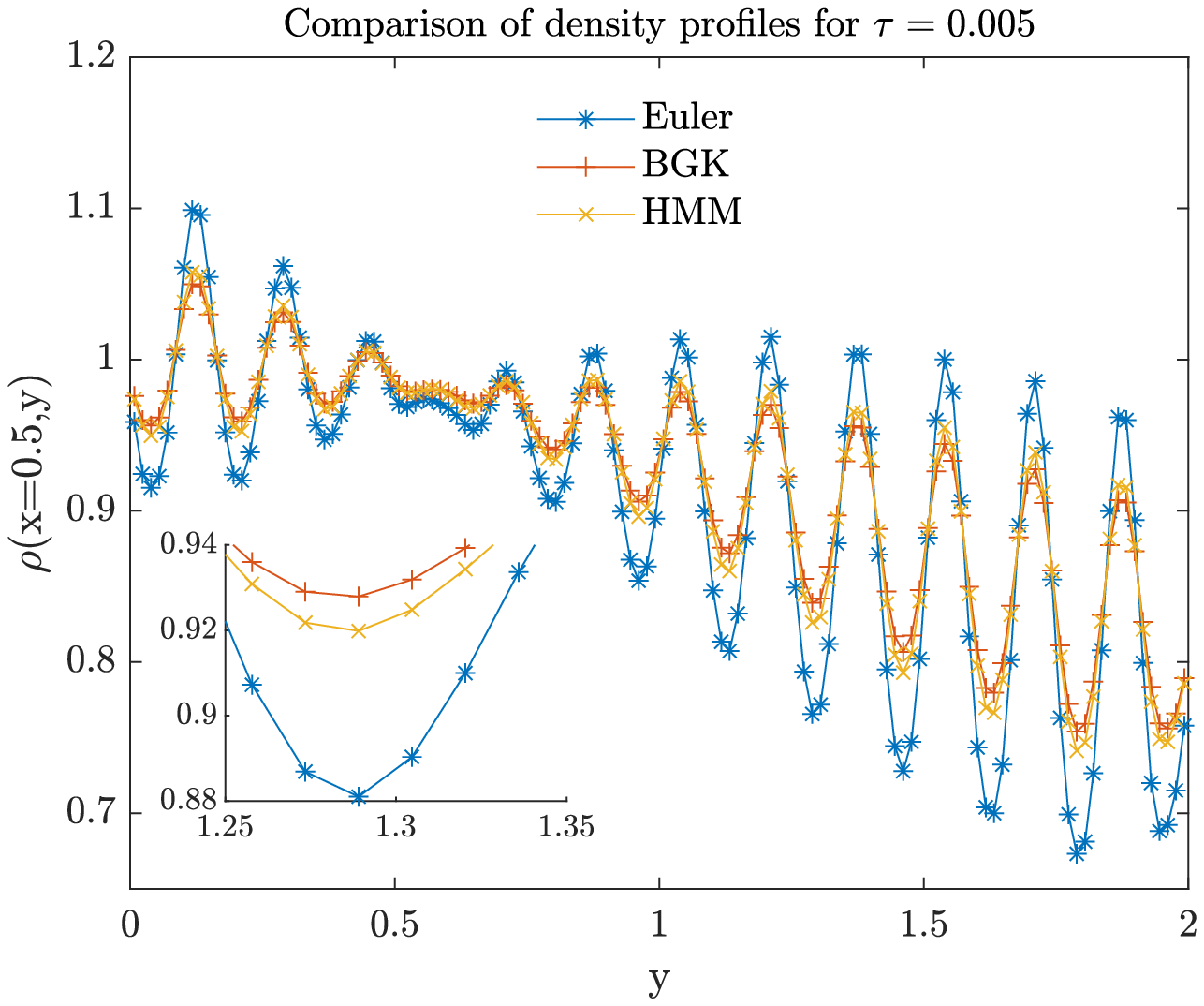}
	\includegraphics[width=0.3\textwidth]{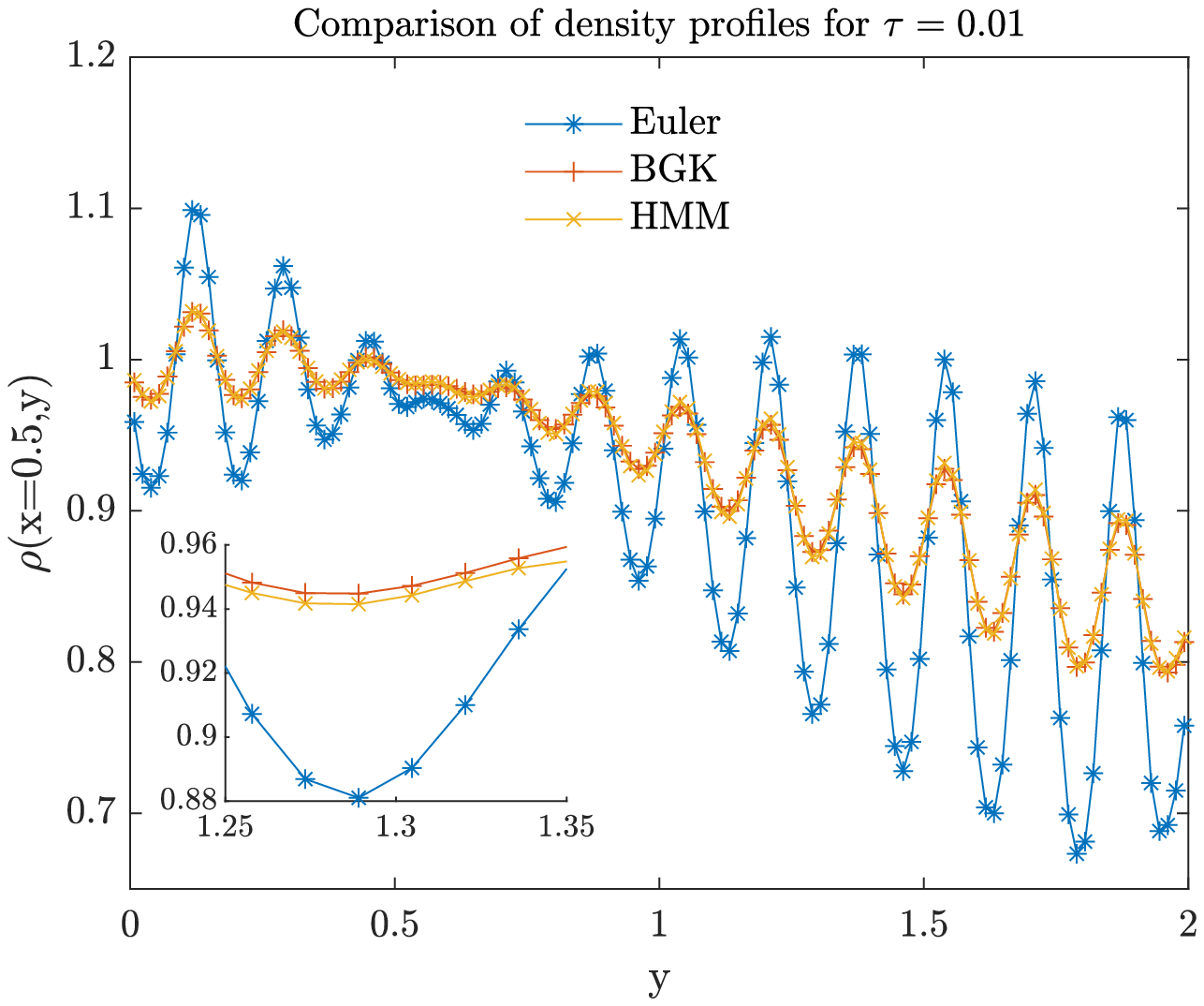}\\
	\includegraphics[width=0.3\textwidth]{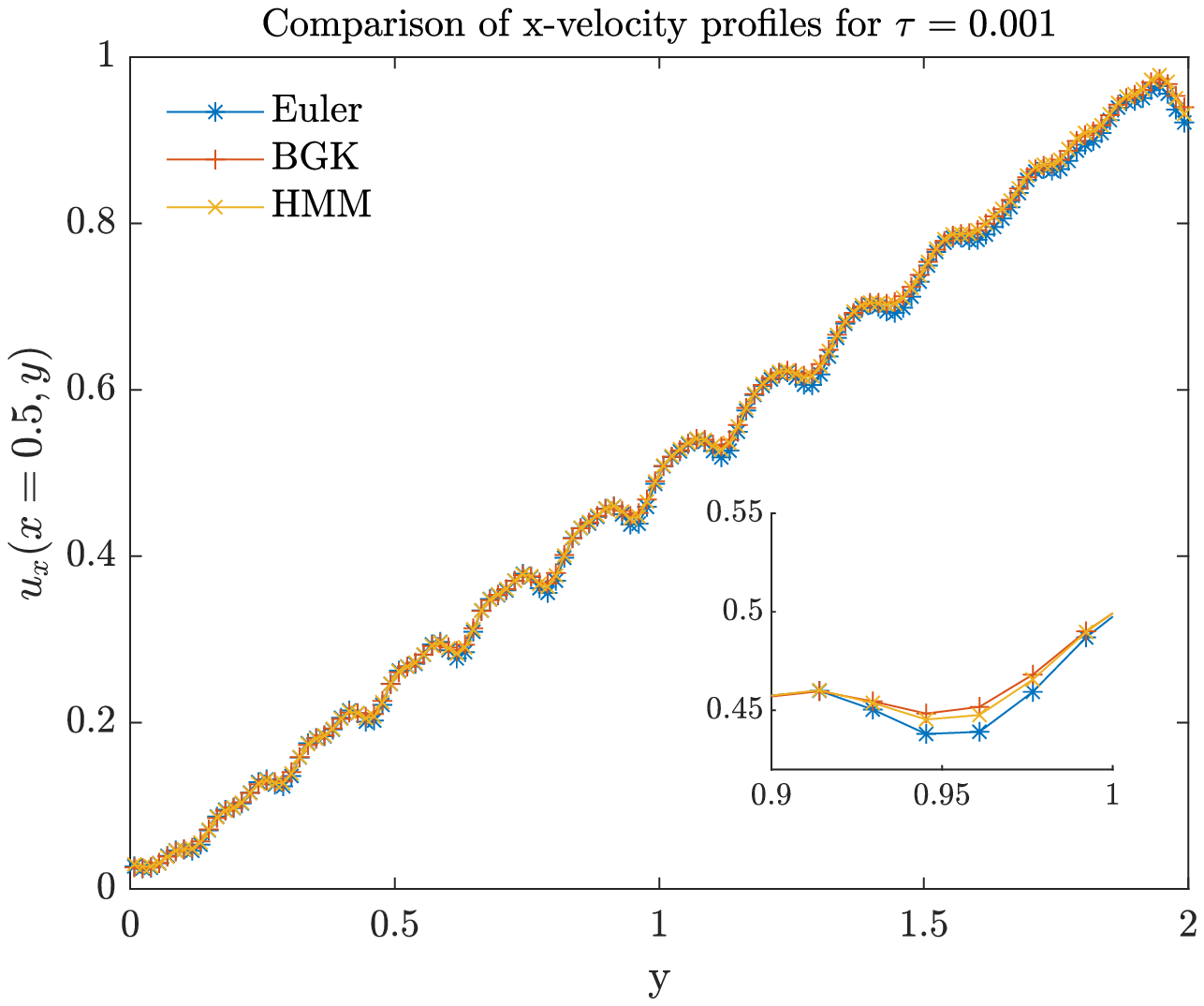}
	\includegraphics[width=0.3\textwidth]{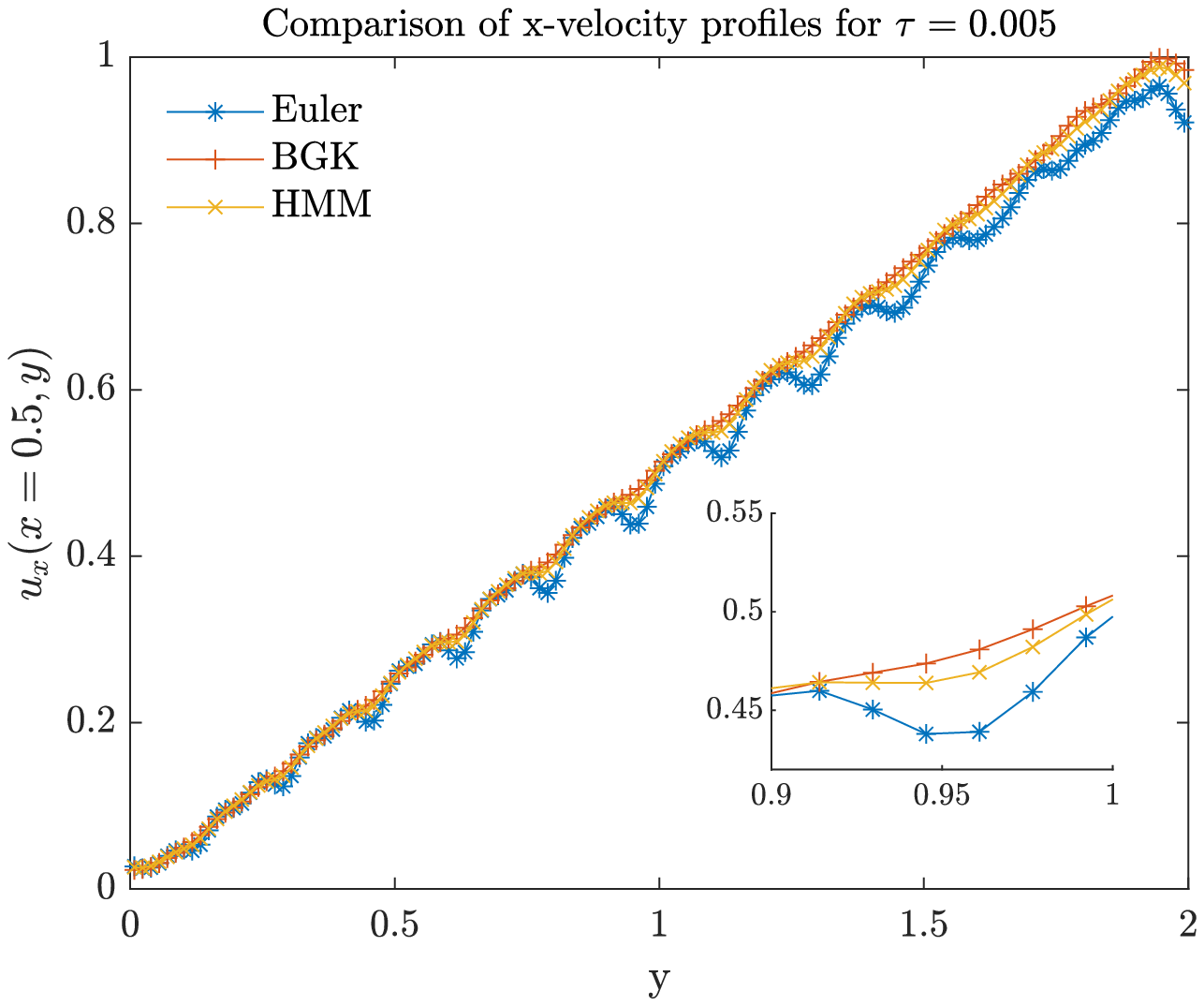}
	\includegraphics[width=0.3\textwidth]{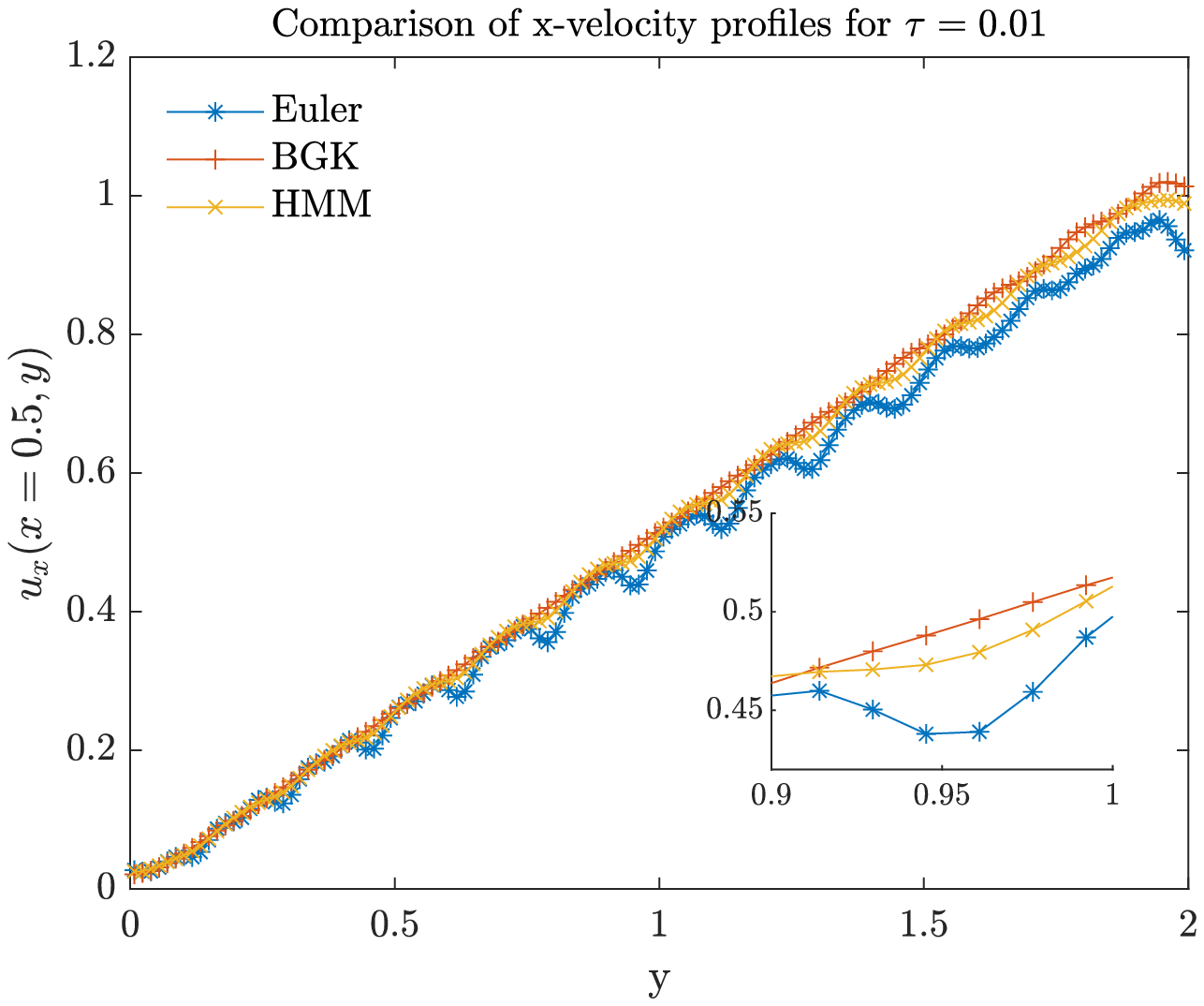}\\	
	\includegraphics[width=0.3\textwidth]{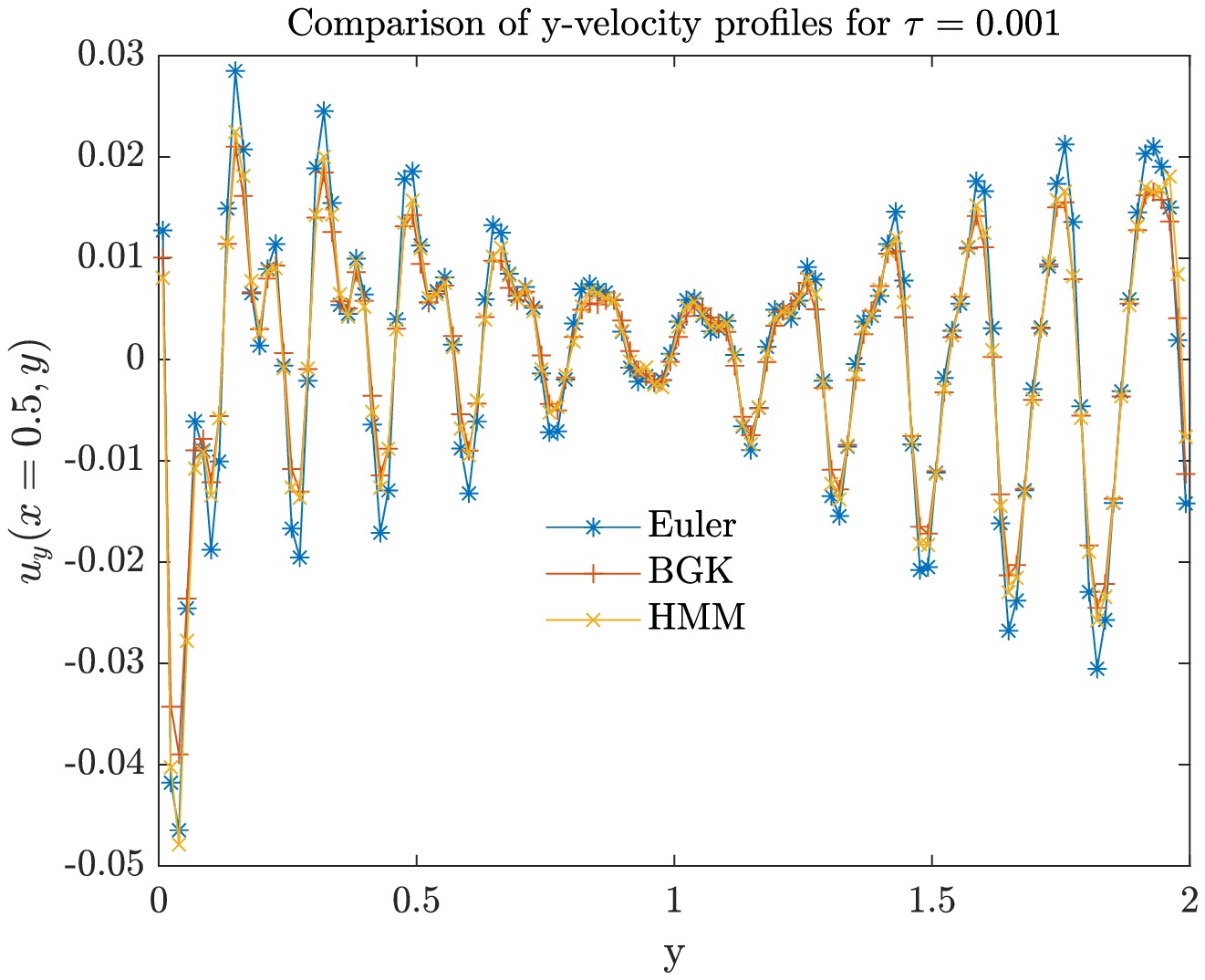}
	\includegraphics[width=0.3\textwidth]{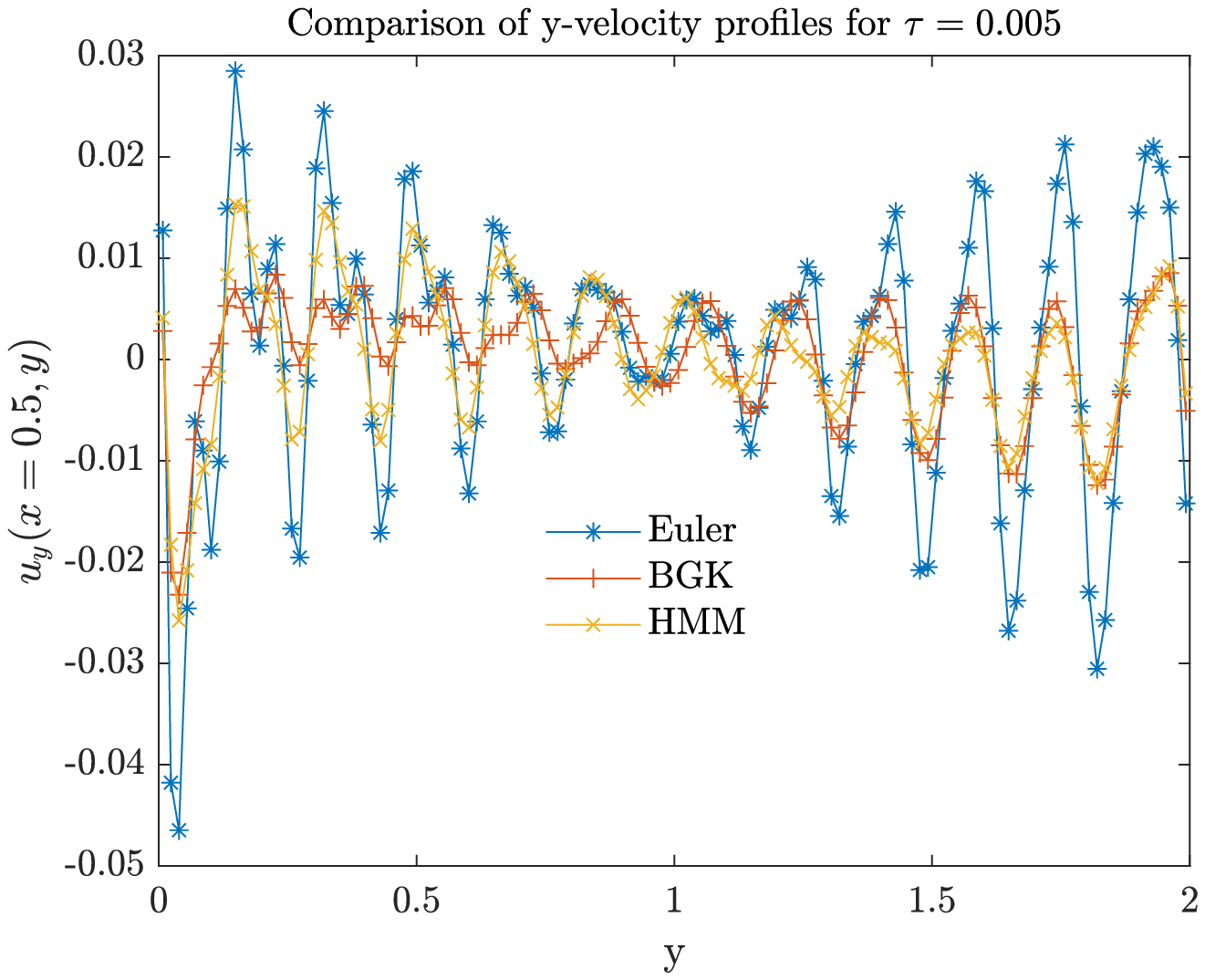}
	\includegraphics[width=0.3\textwidth]{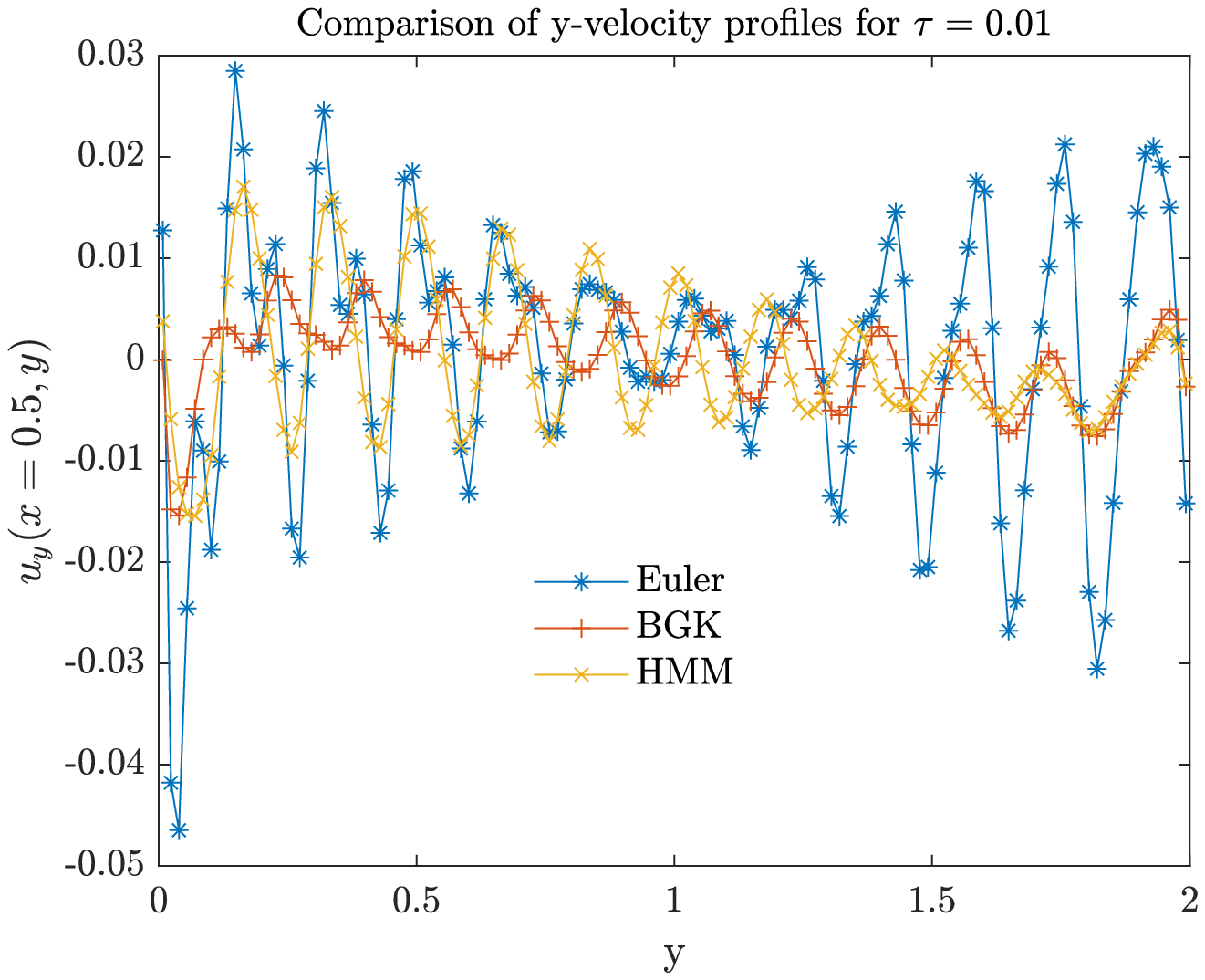}
	\caption{\textbf{Perturbed Couette flow}. Comparison of the density (left), $x$-velocity (middle) and $y$-velocity (right) profiles
	for $x=0.5$. The results for the Euler, the BGK and the HMM equations with respectively $\tau=0.001$, $\tau=0.005$ and $\tau=0.01$ are
	shown.\label{fig:couette2}}
	\end{figure}

\begin{figure}[!ht]
	\centering
	\includegraphics[width=0.3\textwidth]{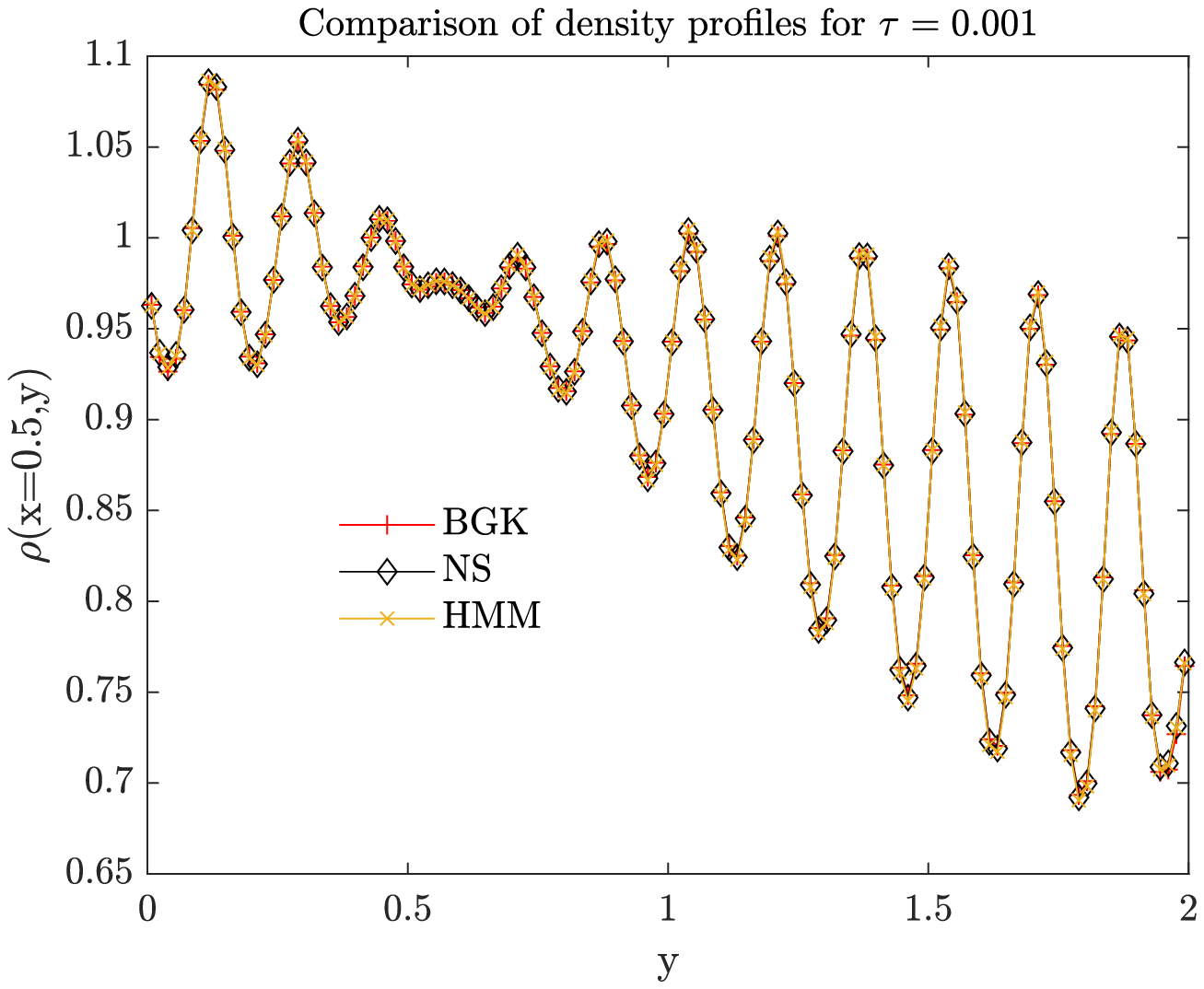}
	\includegraphics[width=0.3\textwidth]{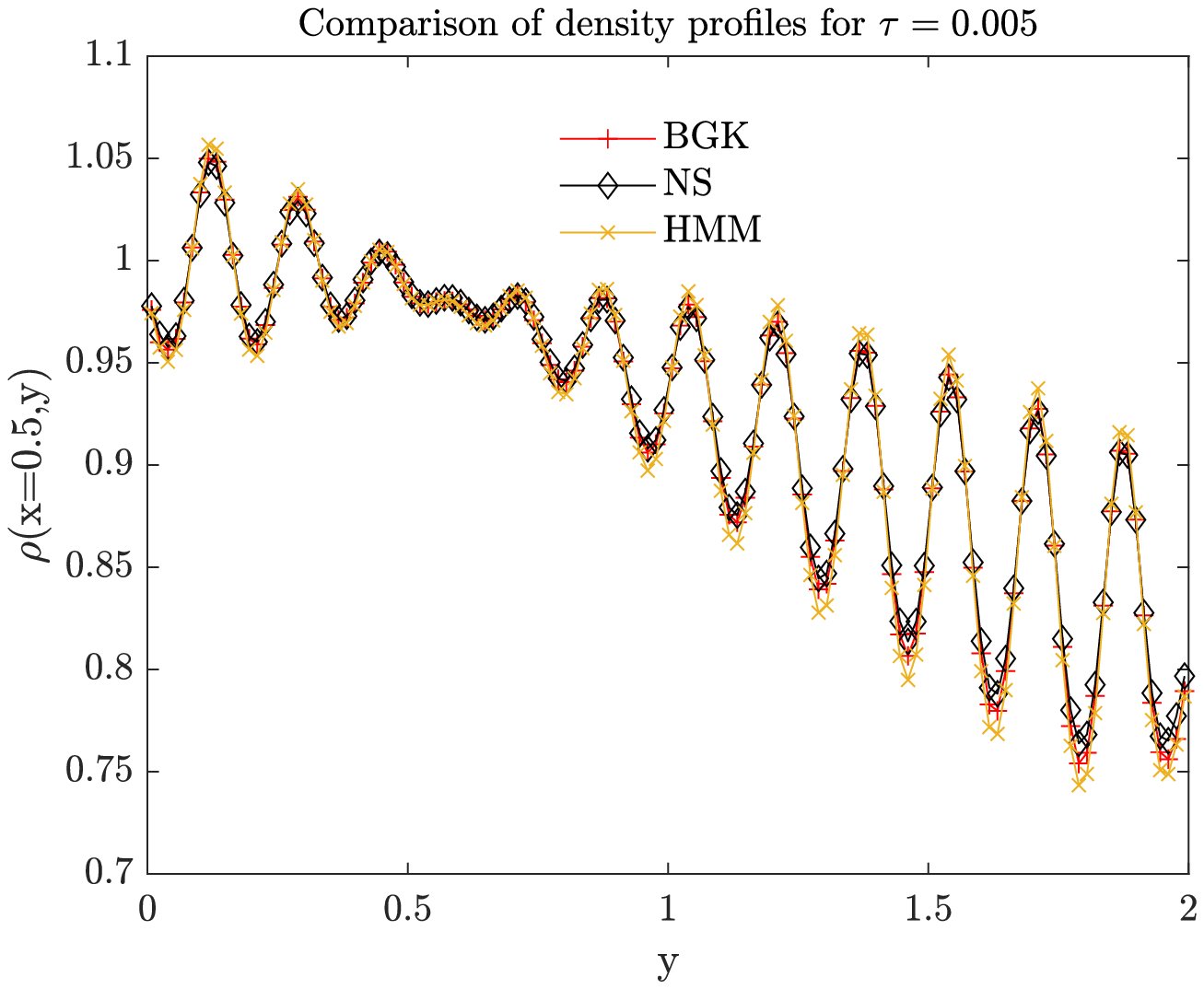}
	\includegraphics[width=0.3\textwidth]{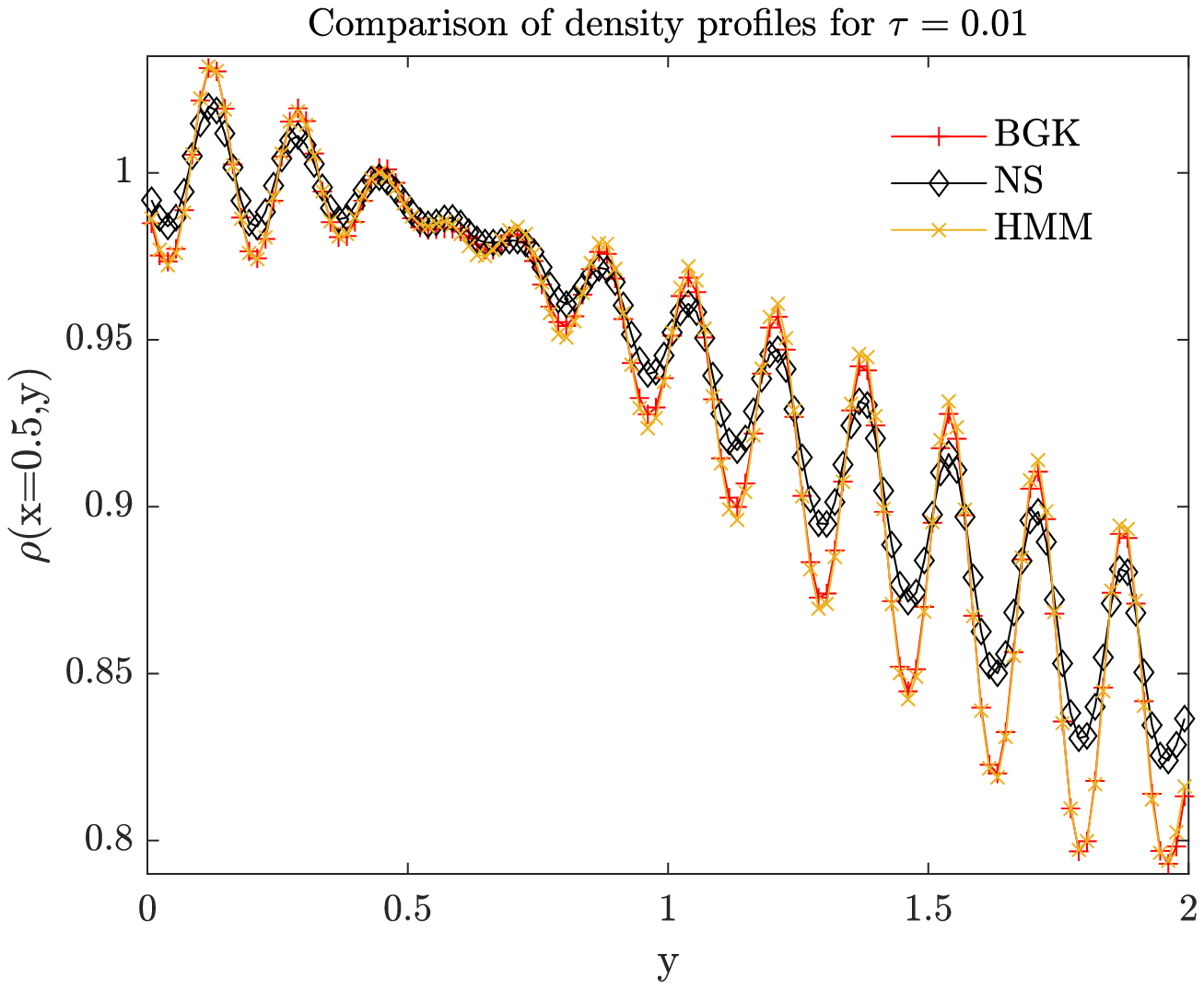}\\
	\includegraphics[width=0.3\textwidth]{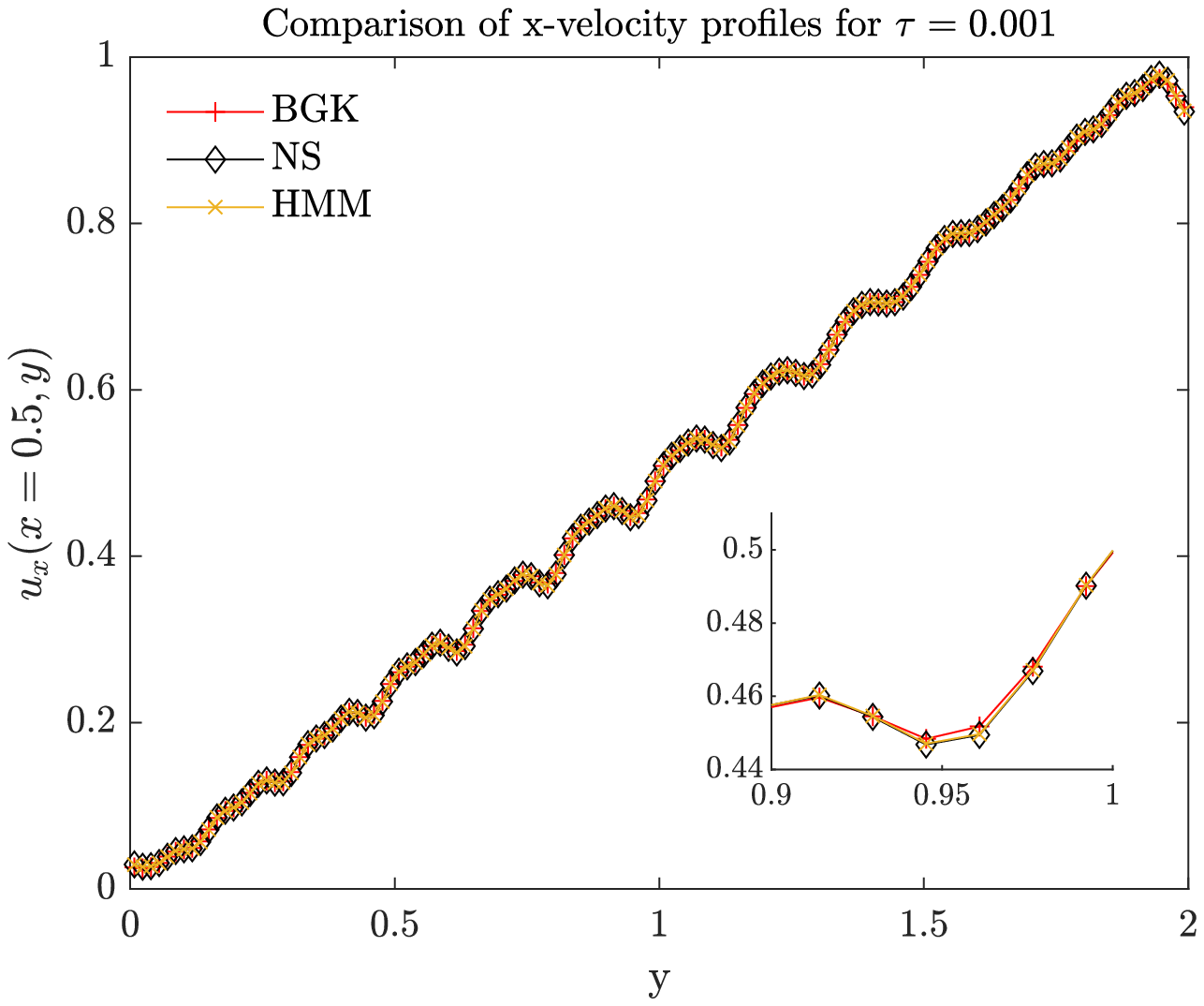}
	\includegraphics[width=0.3\textwidth]{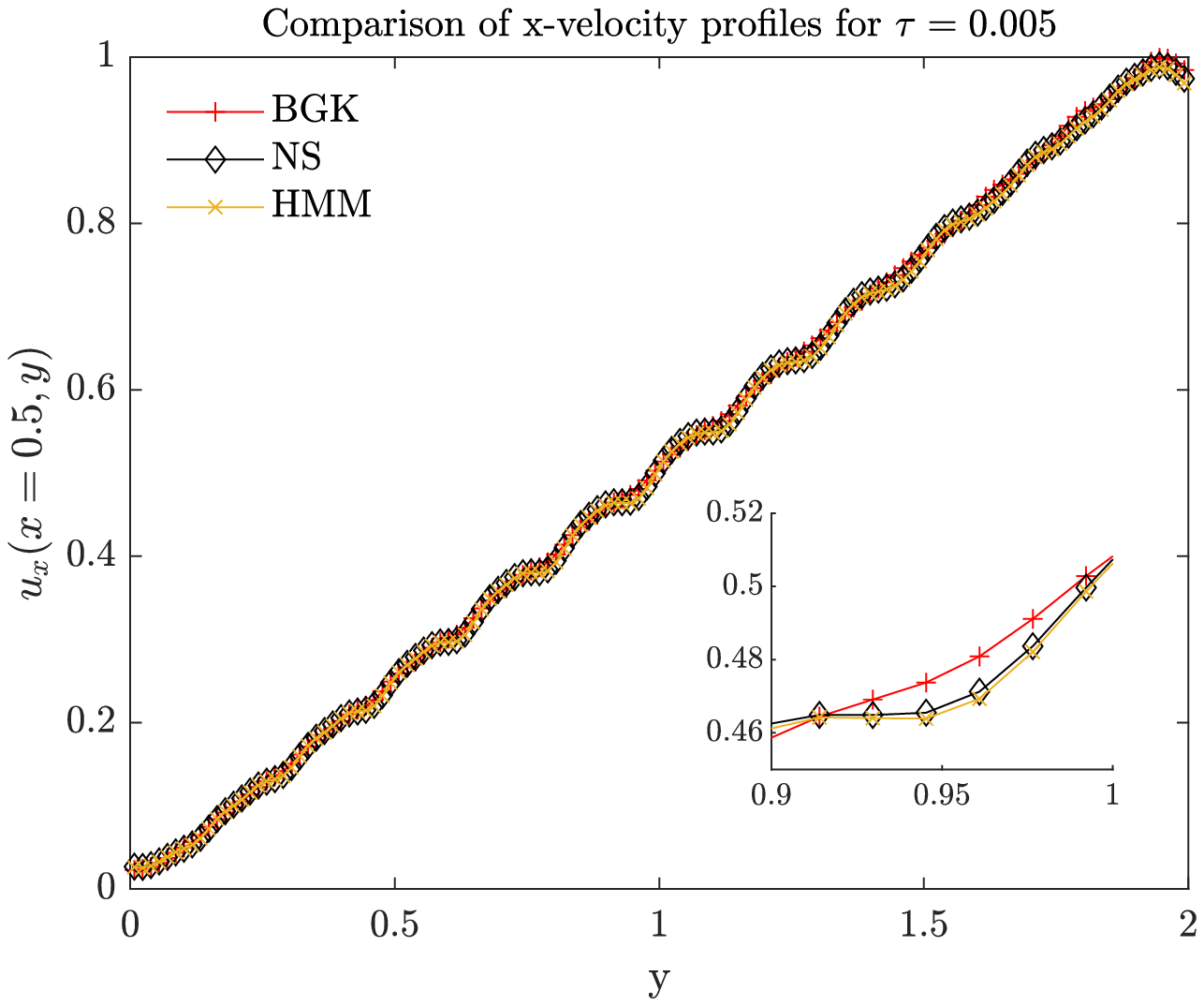}
	\includegraphics[width=0.3\textwidth]{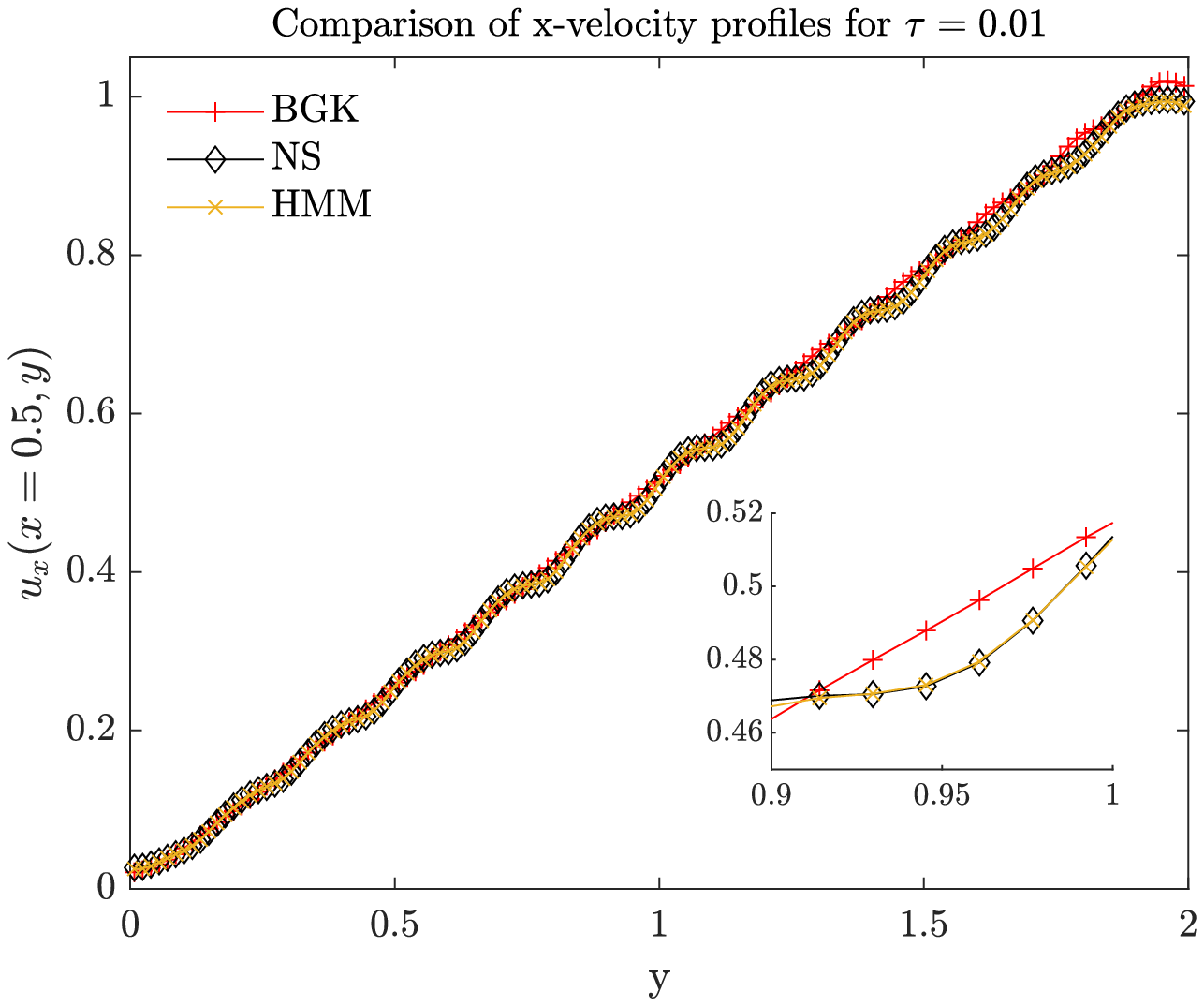}\\	
	\includegraphics[width=0.3\textwidth]{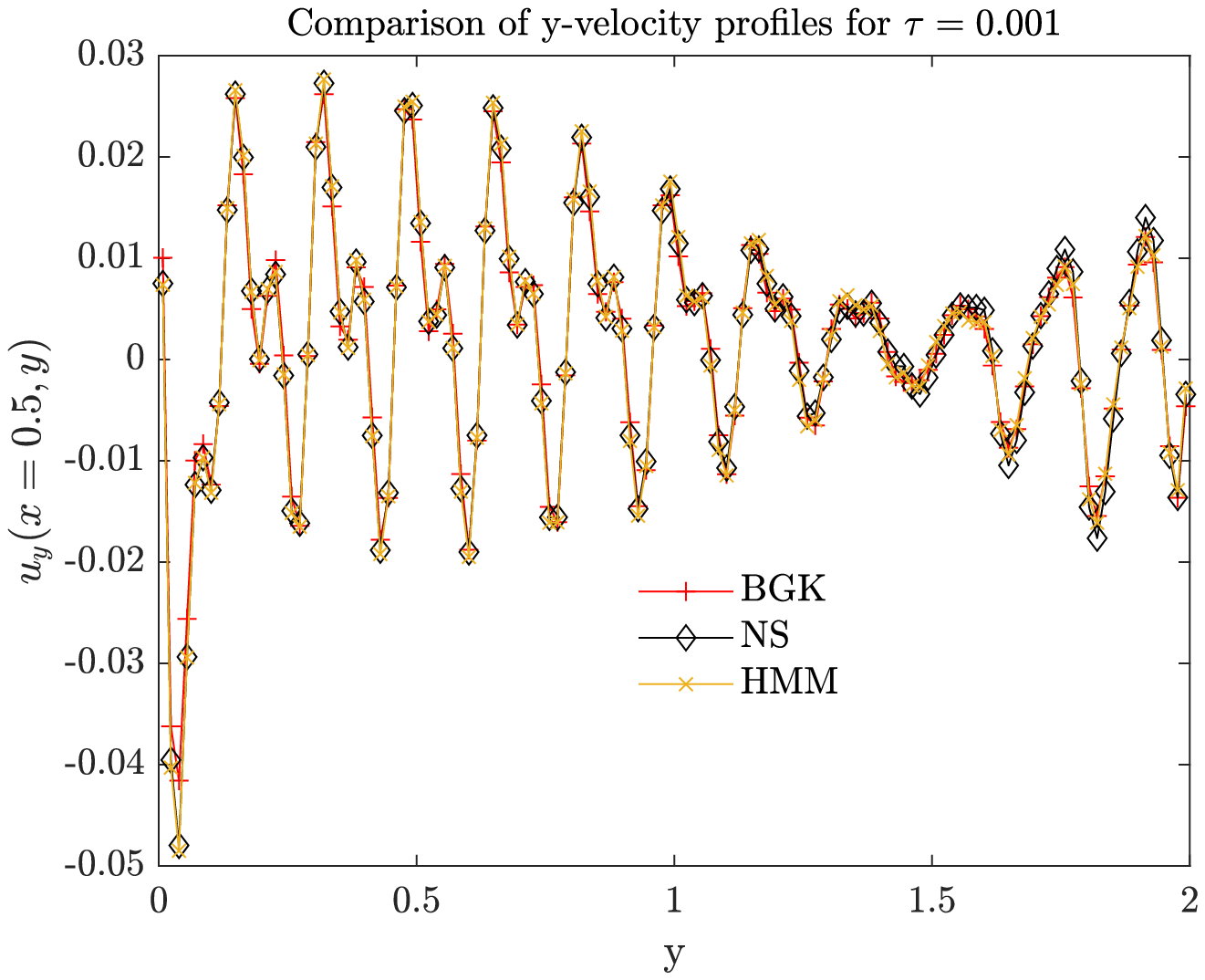}
	\includegraphics[width=0.3\textwidth]{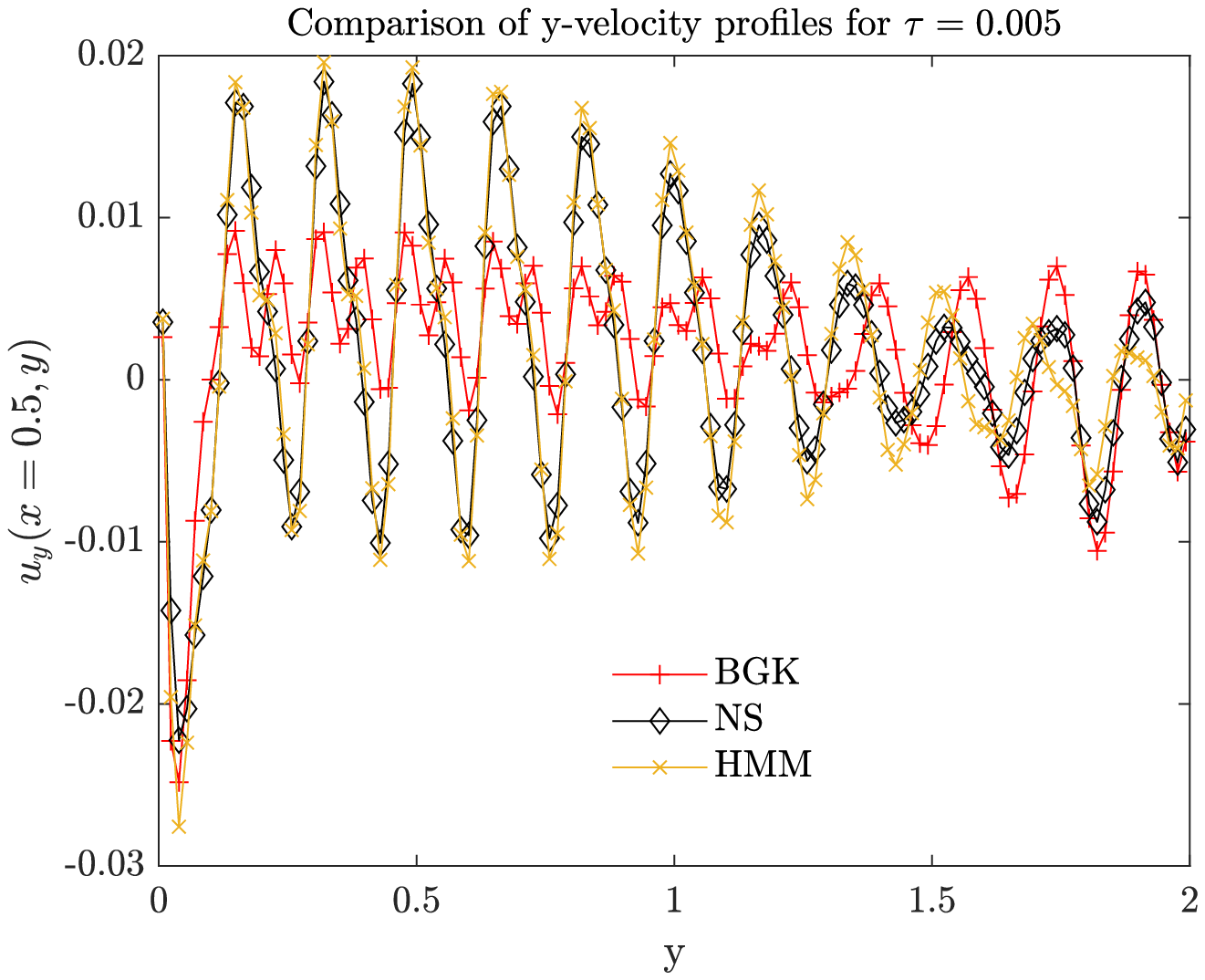}
	\includegraphics[width=0.3\textwidth]{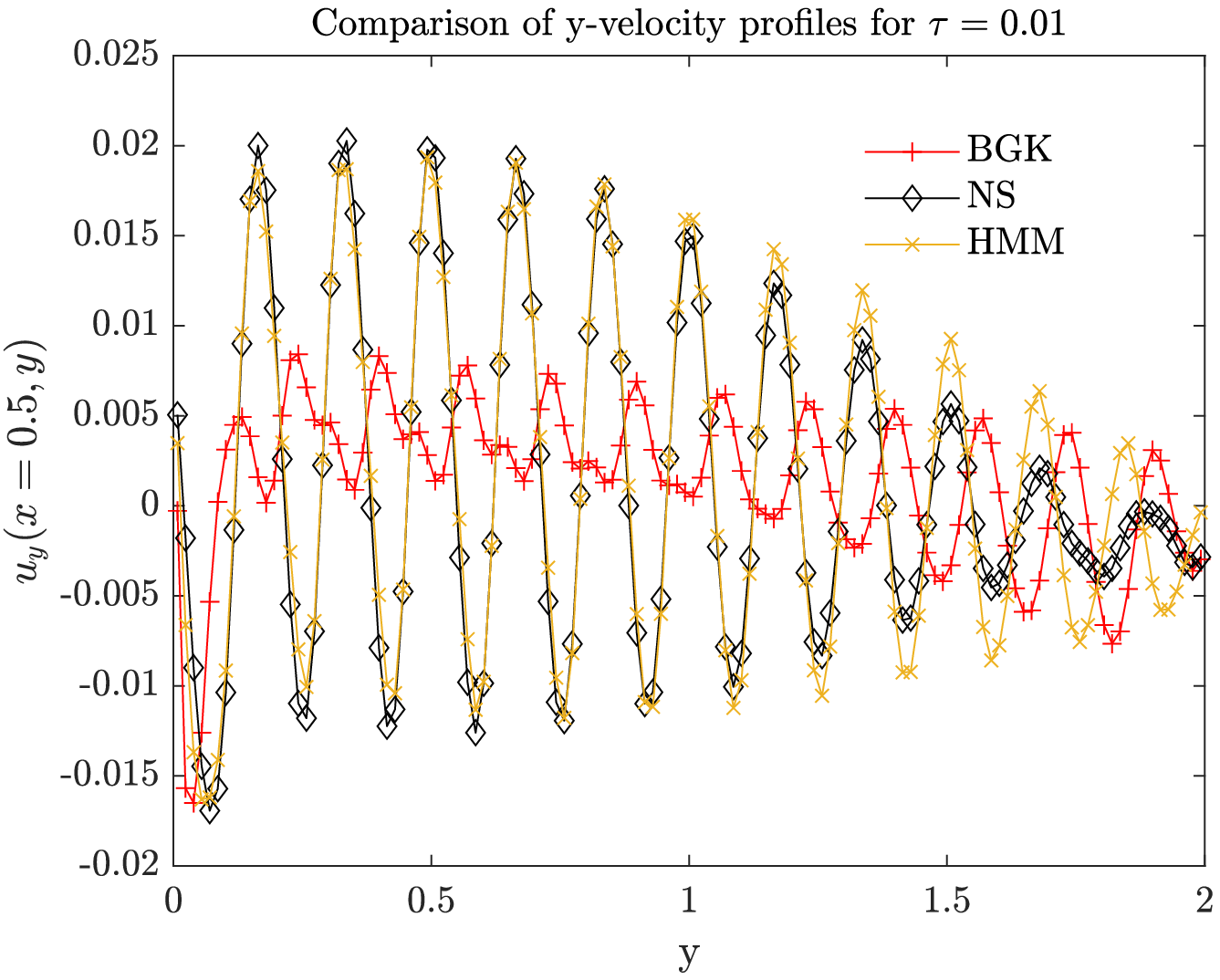}
	\caption{\textbf{Perturbed Couette flow}. Comparison of the density (left), $x$-velocity (middle) and $y$-velocity (right) profiles
	for $x=0.5$. The results for the the BGK, the Navier-Stokes and the HMM equations with respectively $\tau=0.001$, $\tau=0.005$ and $
	\tau=0.01$
	are shown.\label{fig:couette3}}
\end{figure}

In this section, we consider a two-dimensional Couette-type flow in the domain $\Omega=[0,1]\times [0,2]$ perturbed by high frequency
waves. The setup is the following. The boundary conditions in the $x$-direction are periodic while in the $y$-direction we impose the so-called no-slip boundary conditions for which the normal to the walls velocity is fixed to zero while the tangential velocity is set equal to the velocity of the walls. The spatial domain is discretized as before with $N_x\times N_y=64\times 128$ cells, while the velocity
space is approximated by using $20\times 20$ cells in the domain $[-5,5]^2$. Local equilibrium for the distribution function is supposed
at time $t=0$. The initial density and velocities are
\begin{equation*}
	\label{initial_couette}
	\rho(x,y,0)=1+0.2\cos(10\pi x)\sin(12\pi y),\ \bm{u}(x,y,0)=(u_x,u_y)=(y(1+0.125\sin(8\pi x)),0).
\end{equation*}
In other words, we simulate a steady state Couette flow with an added high frequency perturbations both in the density and the velocity.
In particular, to the top wall and in the $y$-direction it is imposed an oscillatory velocity in the $x$-direction which reads
$$
u_y(x,2,t)=1+0.25\sin(8\pi x).
$$
The temperature and the universal gas constant are set $T=1$ and $R=1$ as before. In Figure \ref{fig:couette1}, we show the density and
the velocities profiles after $100$ time iteration with a time step of $\Delta t=6.25 \ 10^{-4}$. At this time, the perturbed waves have
not been completely damped out: it is still possible to observe how the oscillations in velocity and density modify the solution as the
intensity of the relaxation parameter grows. In particular, Figure \ref{fig:couette1} shows the solutions obtained with $\tau=0$,
$\tau=0.005$ and $\tau=0.01$ for the BGK model, while for the HMM we take $\tau_{HMM}=\tau/3$ as for the first test. In Figure
\ref{fig:couette2}, we present the same results at $x=0.5$ comparing the isentropic Euler, the BGK and the HMM solutions. In this picture,
we can clearly see that oscillations are damped when moving far from the regime of validity of the isentropic Euler equations and that the
HMM captures the BGK solution very well for all regimes, considered at least for the density and the first component of the velocity
field. Finally, in Figure \ref{fig:couette3}, we show a comparison between the HMM and the NS model using again the BGK model as
reference, all presented at $x=0.5$. For the specific case of the NS model, we take $\varepsilon=\tau$ while all the others numerical
parameters are equal to the ones of the HMM model. From this figure we see that for small values of $\tau$ both models capture very well
the reference solution given by the BGK equation, while for larger $\tau$, NS tends to overestimate the damping for what concerns the
density profile. Both models is capable of capturing the first component of the velocity field $u_x$, while both fail in the description
of the second component $u_y$ when $\tau=0.01$. Let observe anyway that $u_y$ is affected by waves, which are at least one order of
magnitude smaller than the perturbation waves acting on the other macroscopic quantities and thus the error introduced in the solution
remains very small.

\subsection{A vortex test problem}\label{Vortex}
\begin{figure}[!ht]
	\centering
	\includegraphics[width=0.45\textwidth]{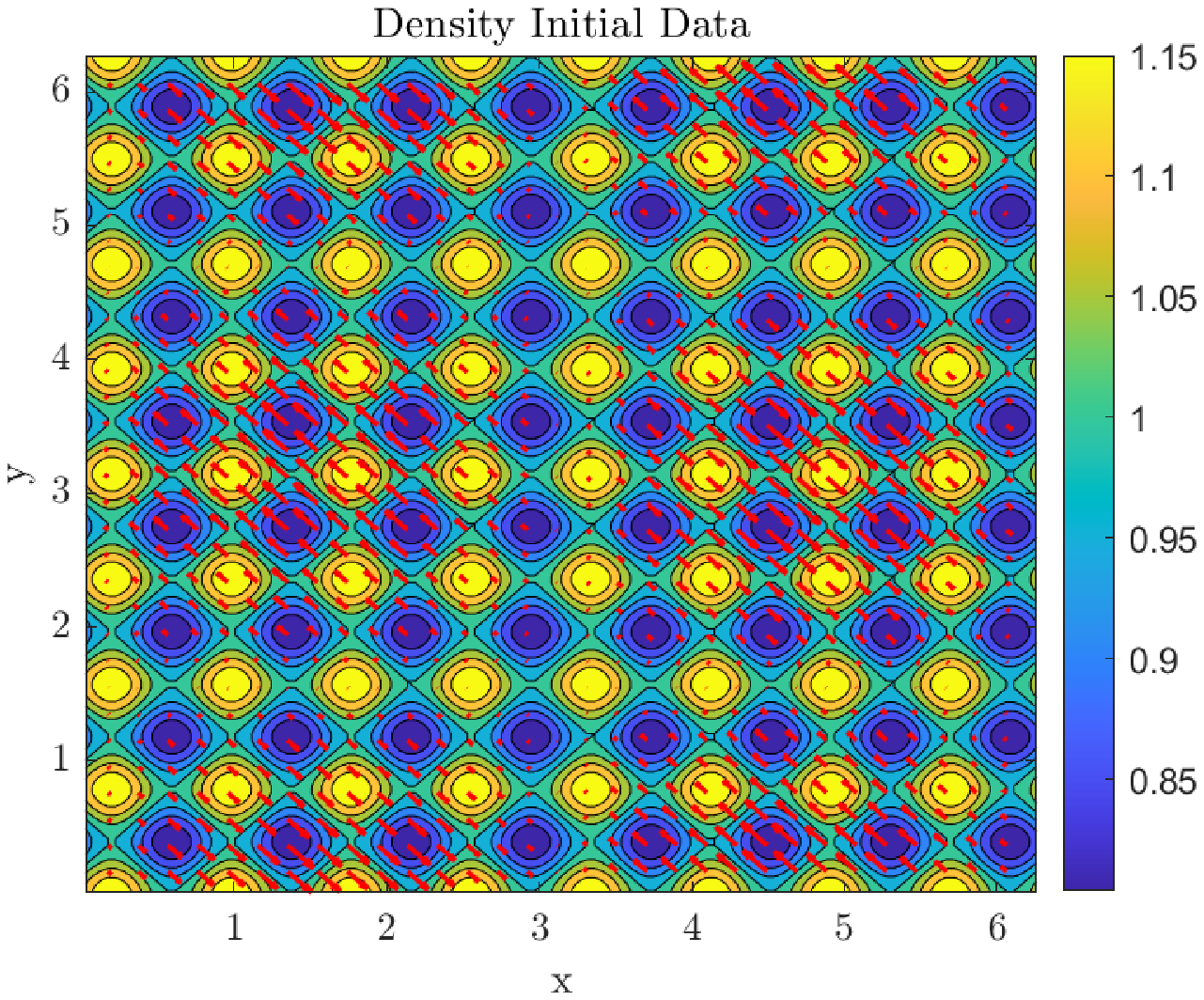}
	\includegraphics[width=0.45\textwidth]{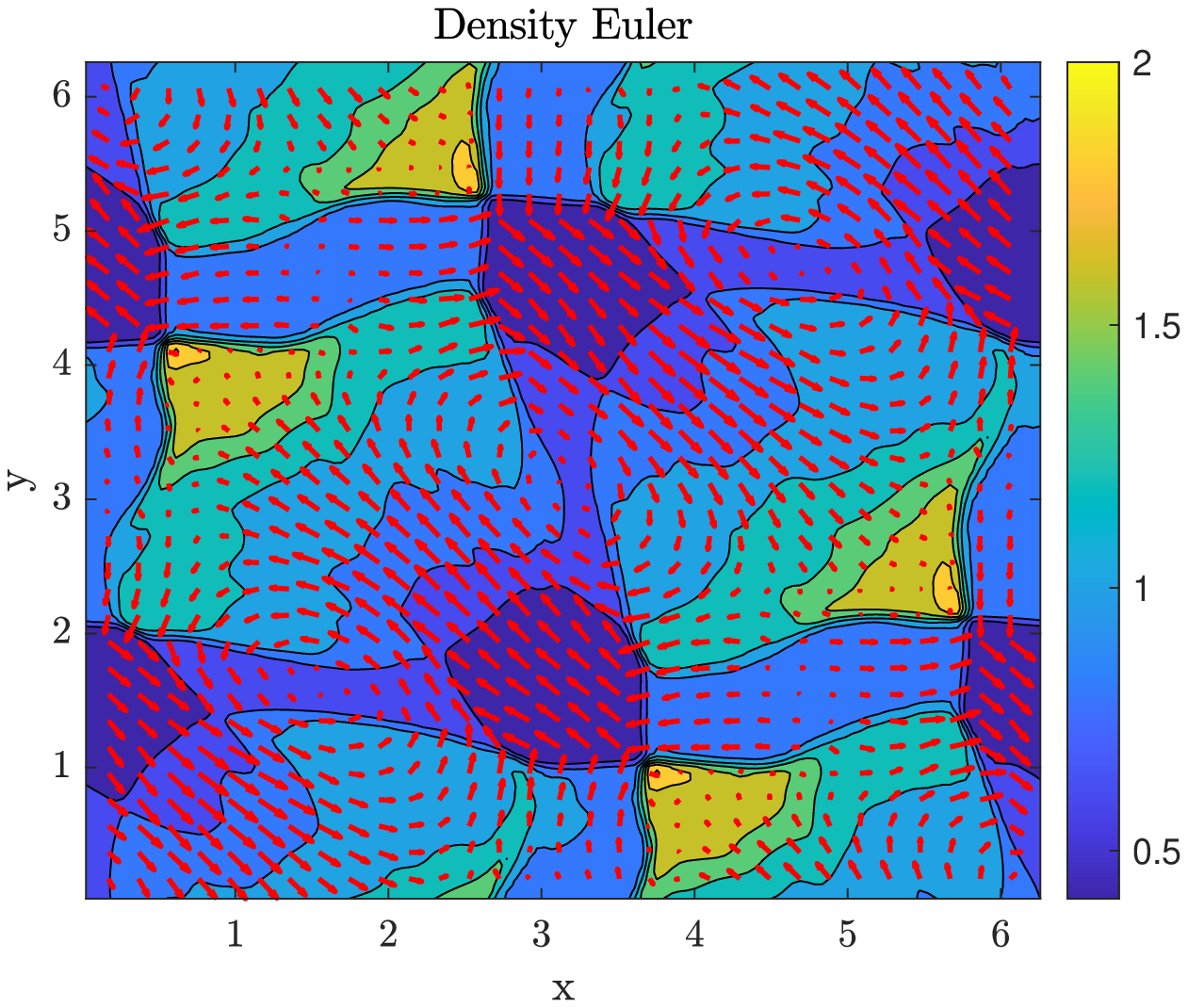}\\
	\includegraphics[width=0.45\textwidth]{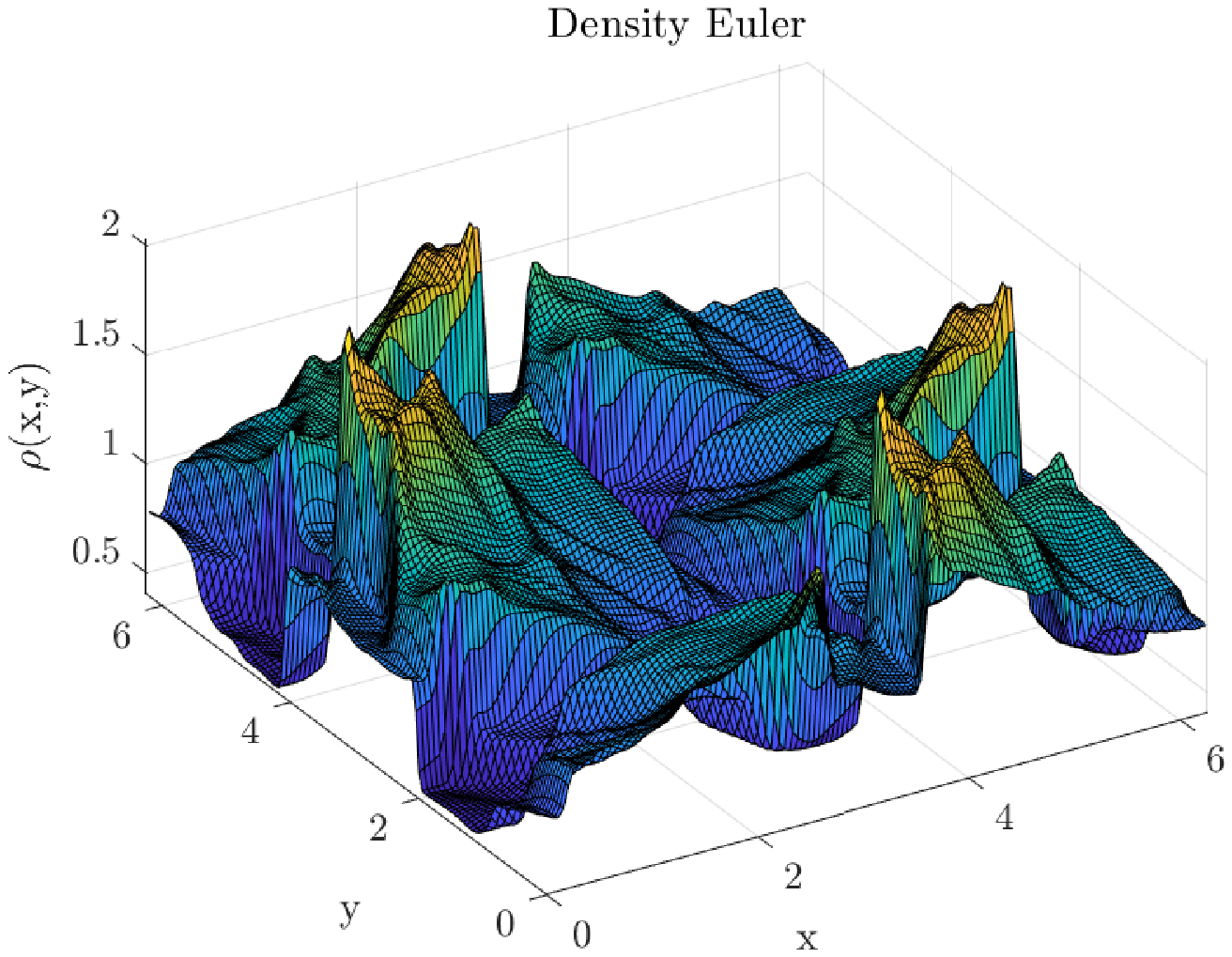}
		\includegraphics[width=0.45\textwidth]{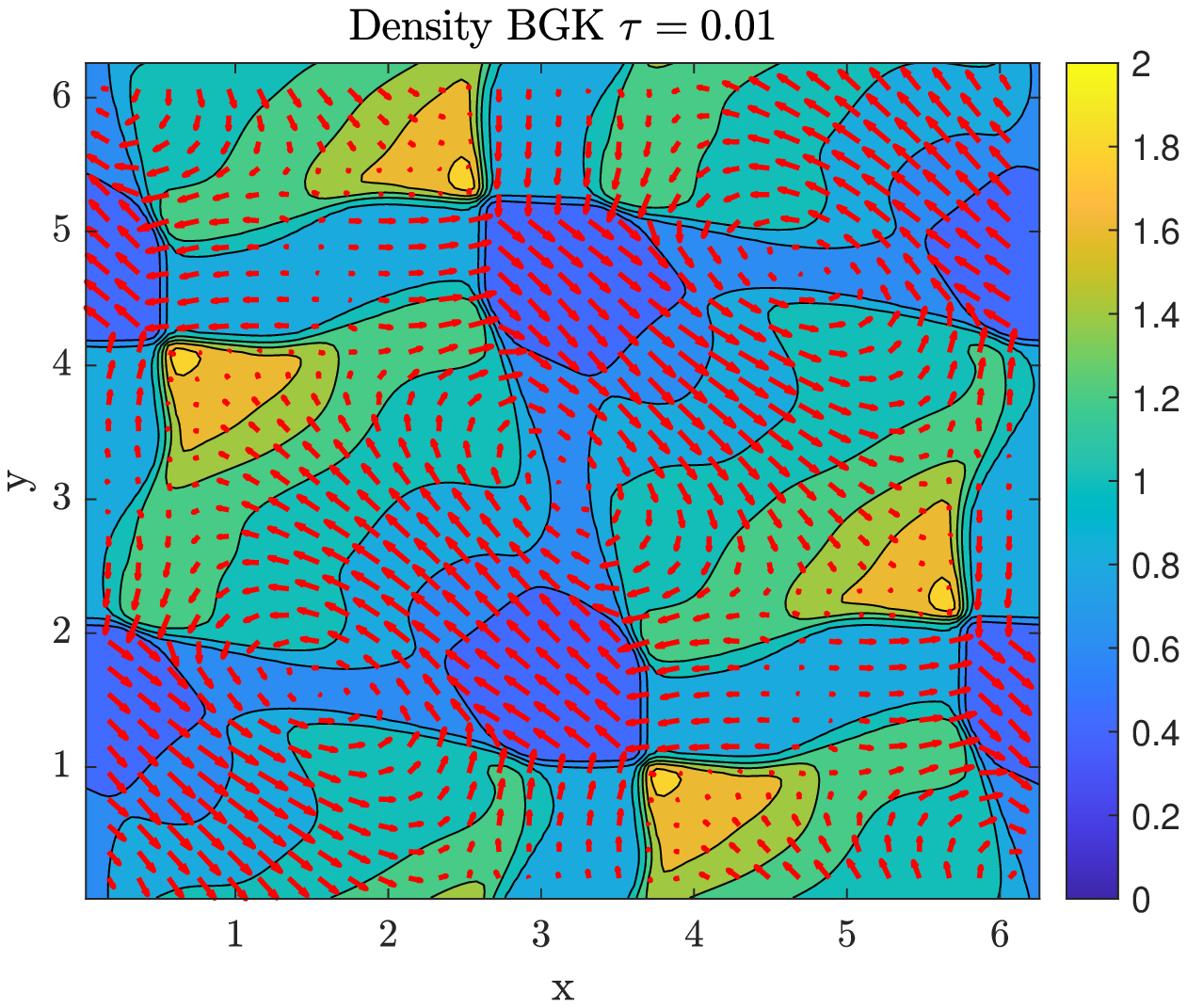}
	\caption{\textbf{Vortex flow}. Profiles of density and velocity for the Euler and the BGK equations. From top to bottom: Density initial data, density Euler solution, density Euler solution three dimensional view and density BGK solution with $\tau=0.01$. Red arrows represent the velocity vector field.  	\label{fig:vortex1}}
\end{figure}

\begin{figure}[!ht]
	\centering
	\includegraphics[width=0.3\textwidth]{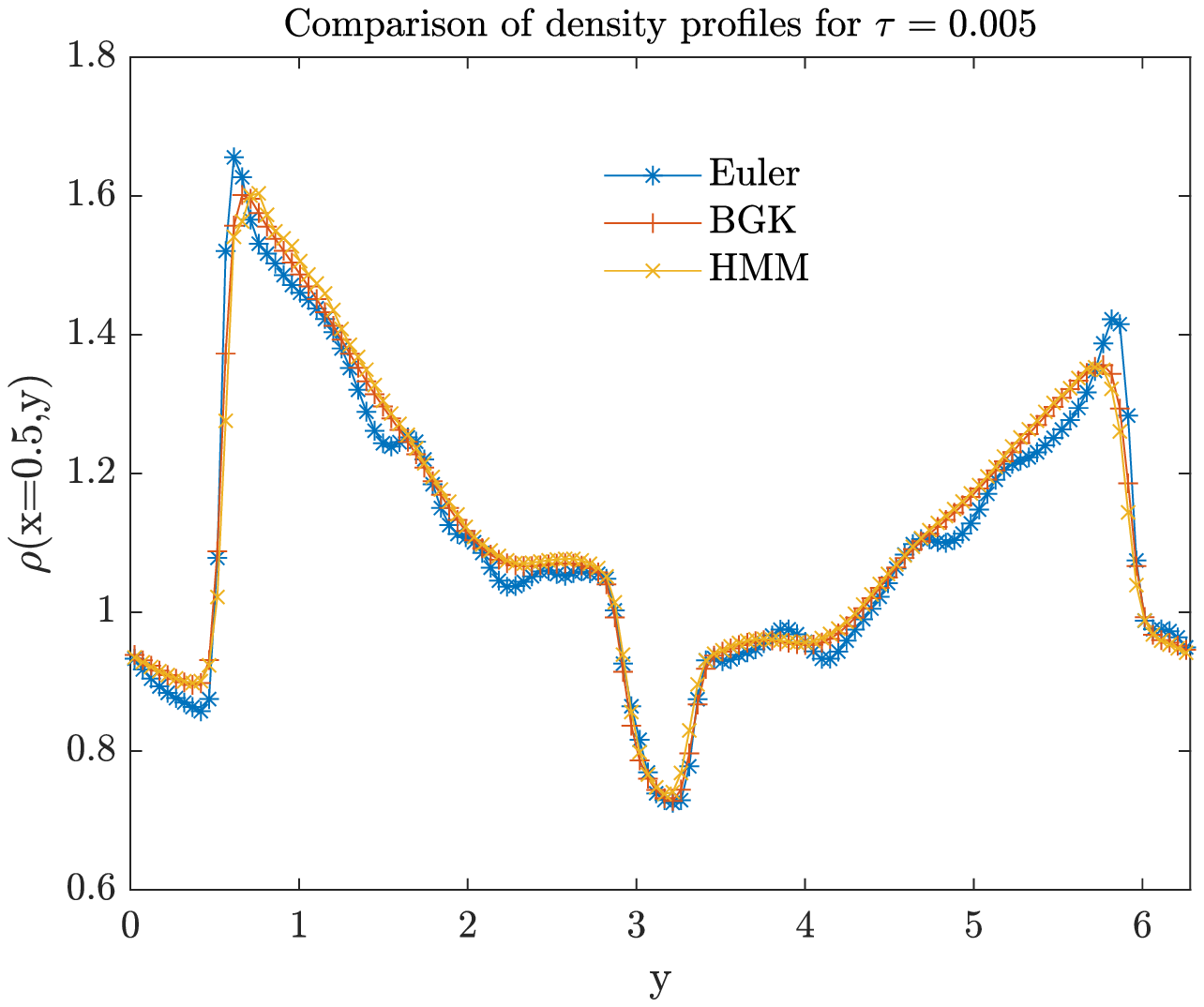}
	\includegraphics[width=0.3\textwidth]{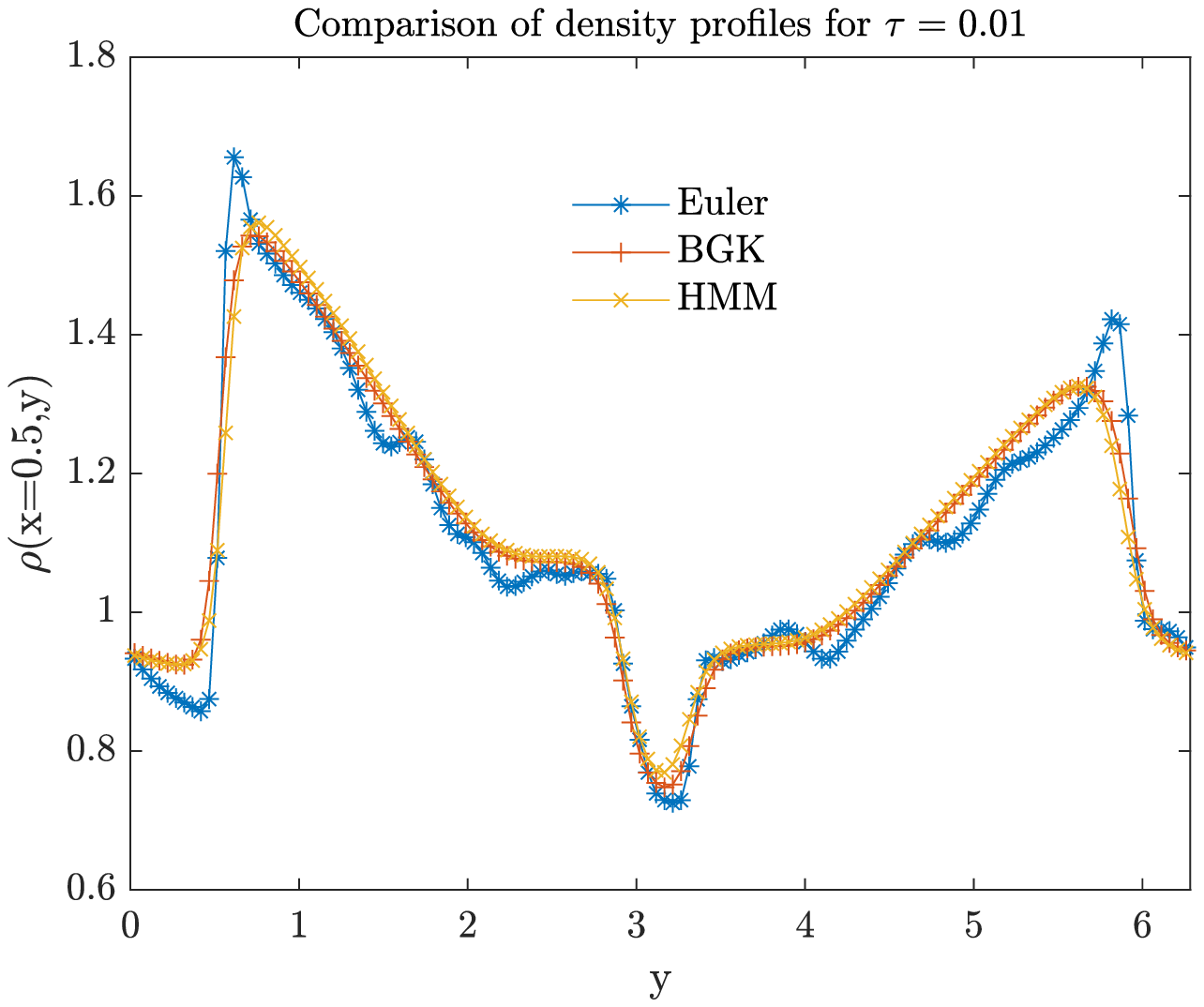}
	\includegraphics[width=0.3\textwidth]{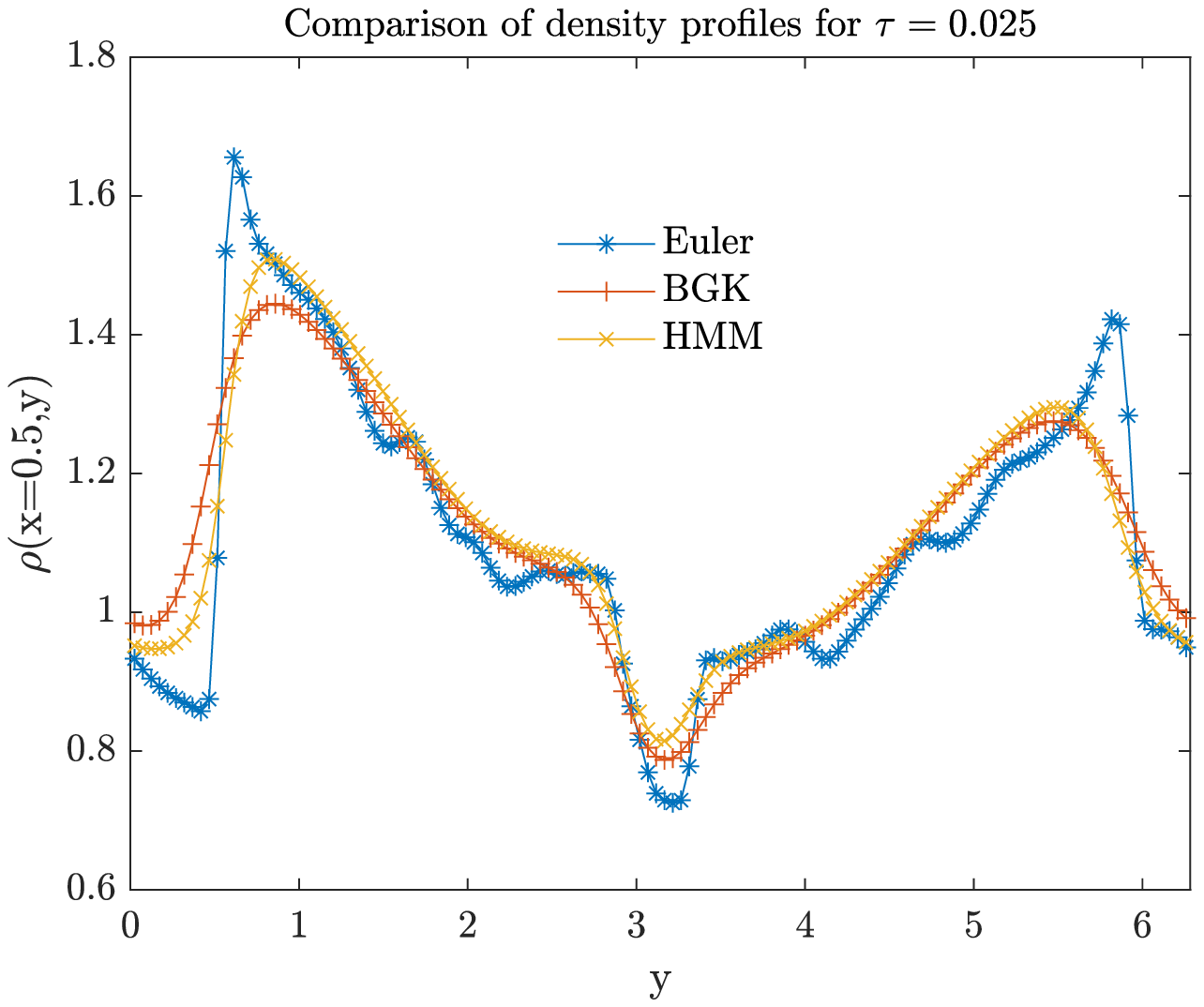}
	\\
	\includegraphics[width=0.3\textwidth]{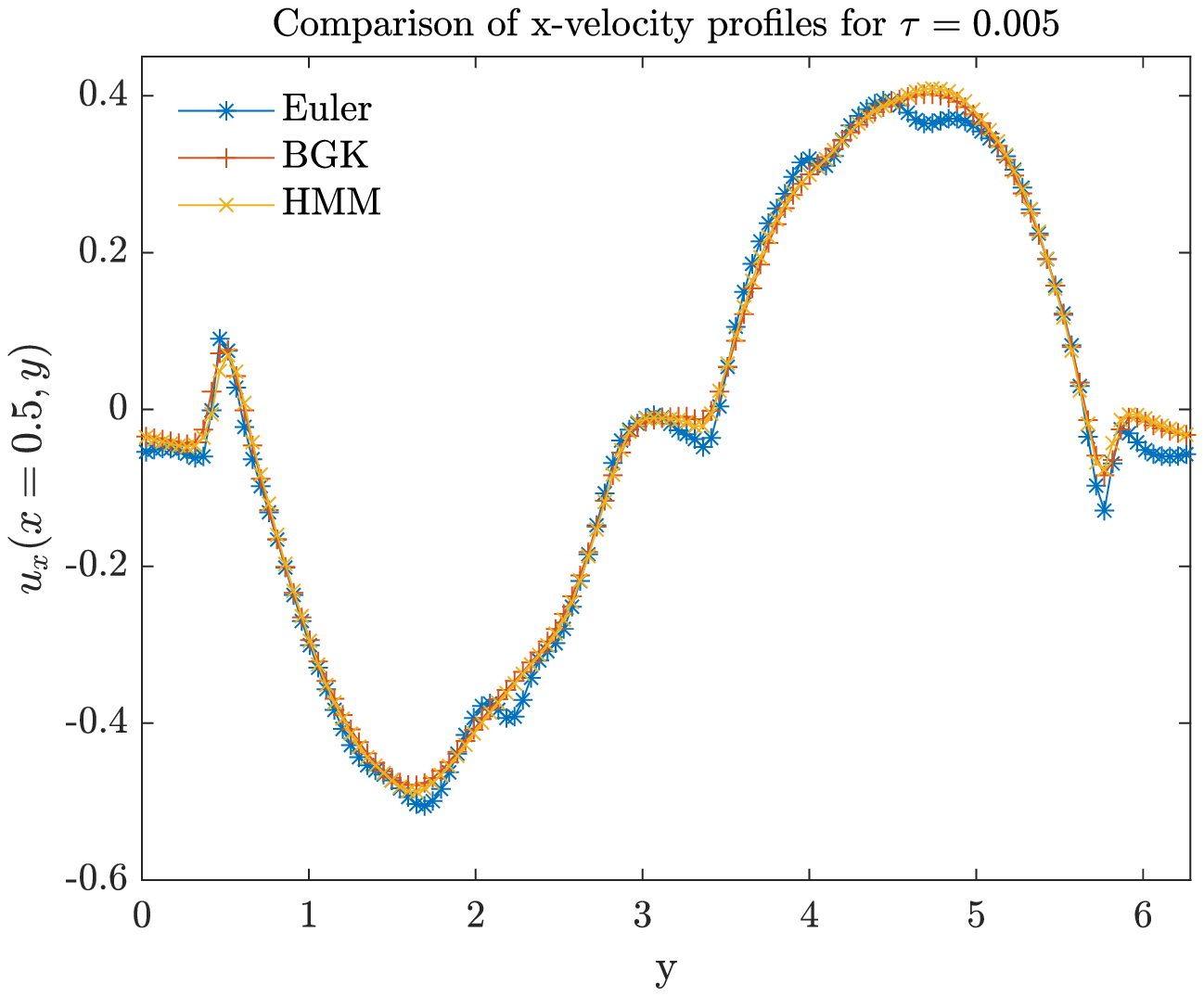}
	\includegraphics[width=0.3\textwidth]{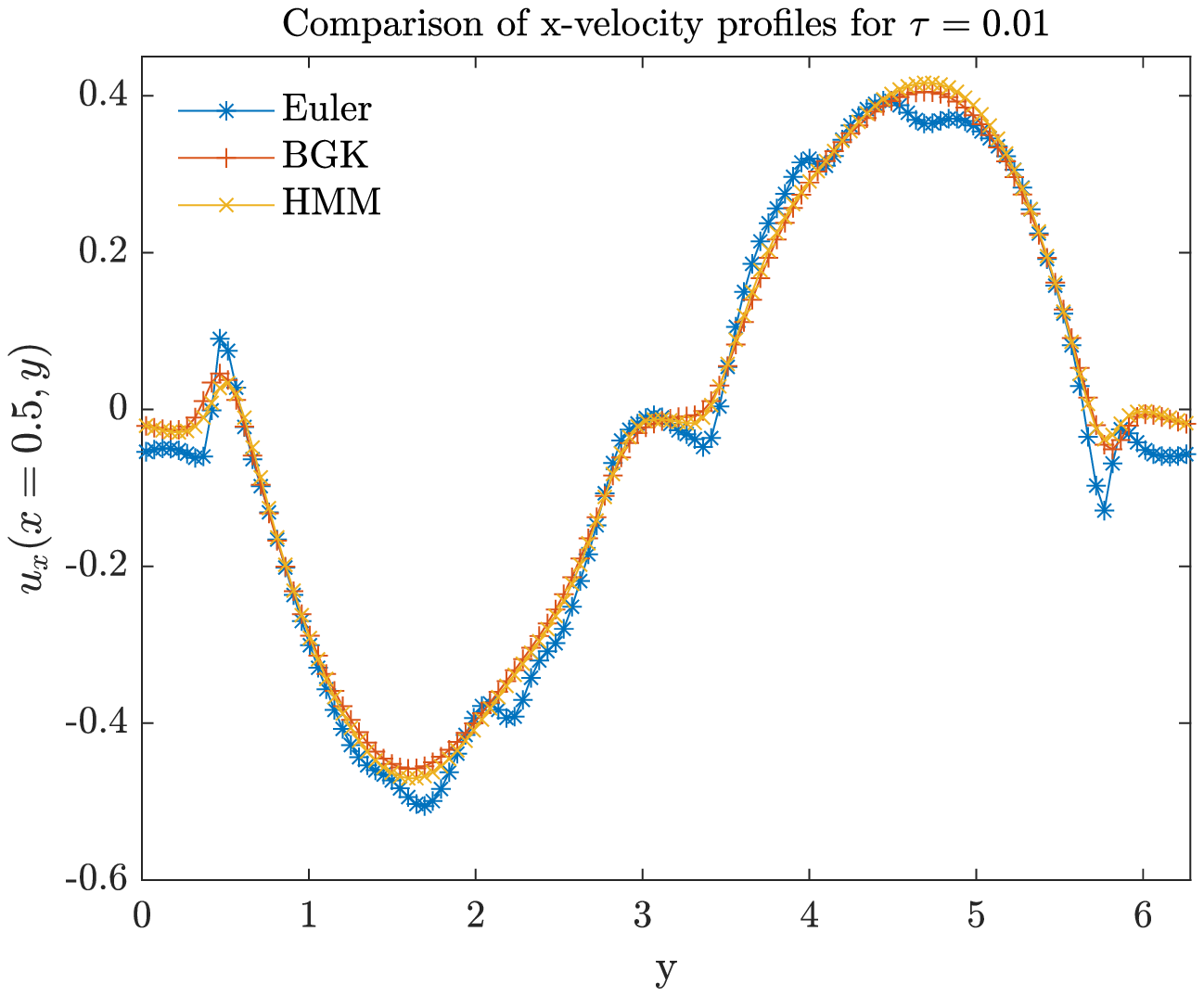}
		\includegraphics[width=0.3\textwidth]{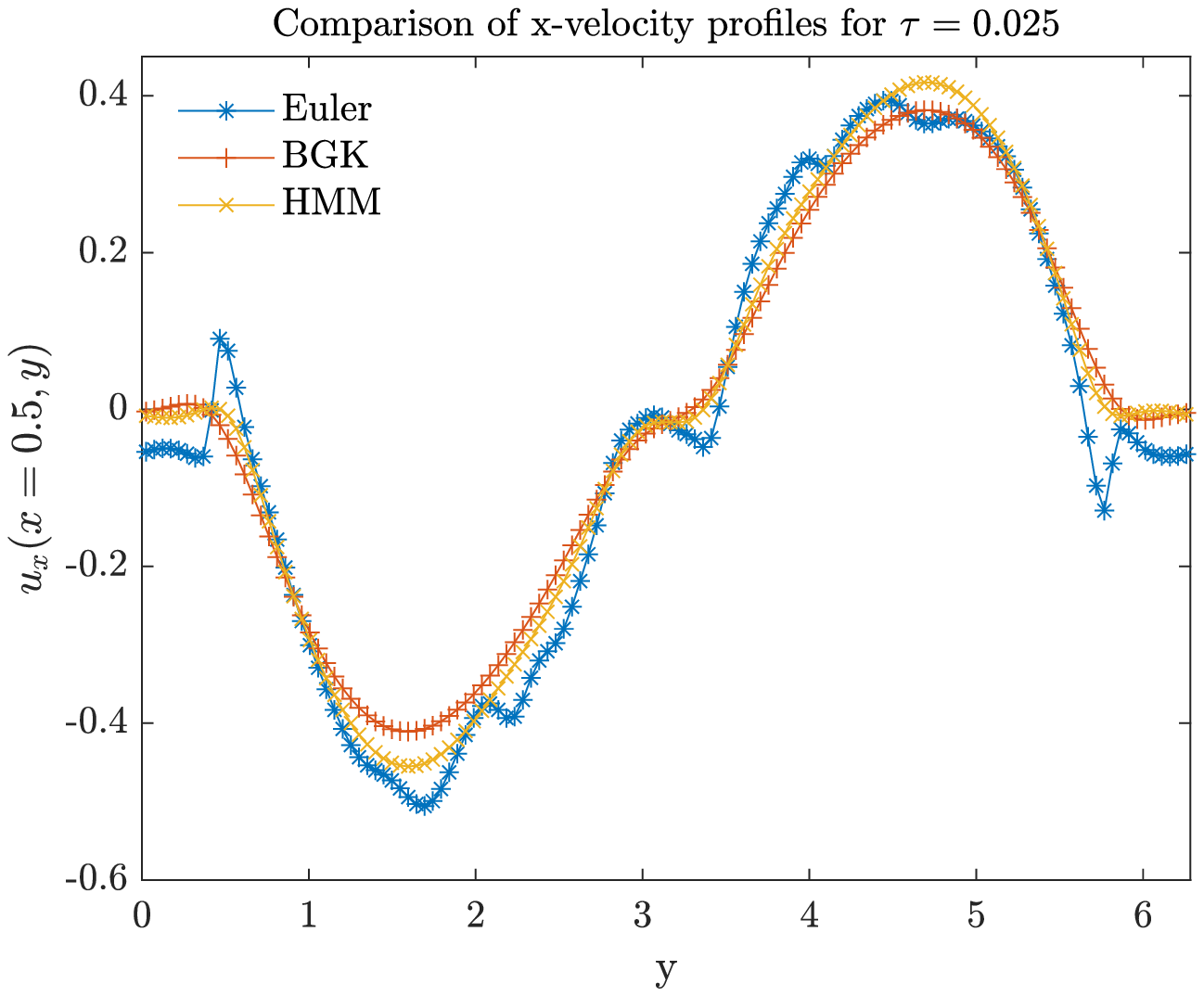}
	\\	
	\includegraphics[width=0.3\textwidth]{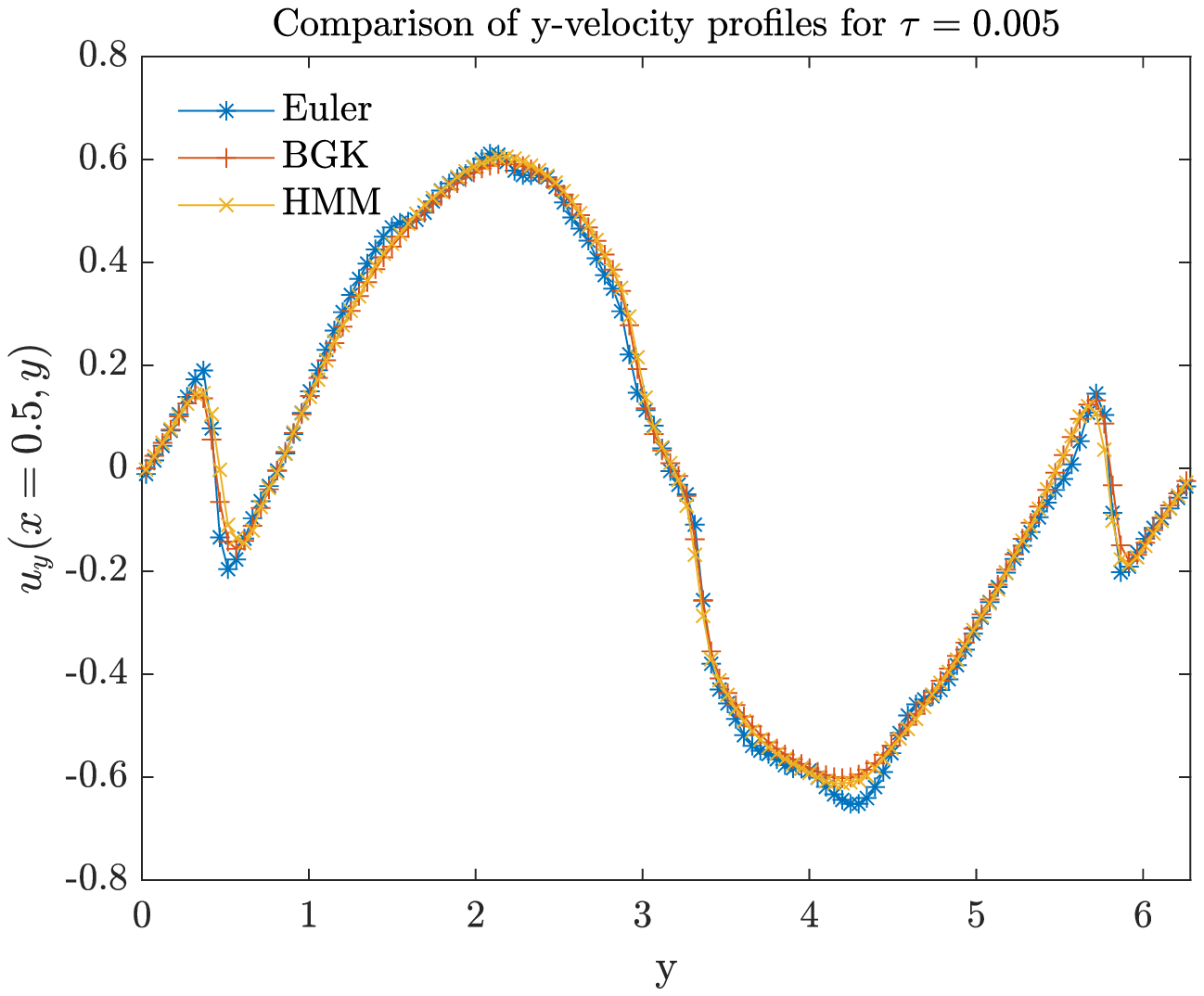}
	\includegraphics[width=0.3\textwidth]{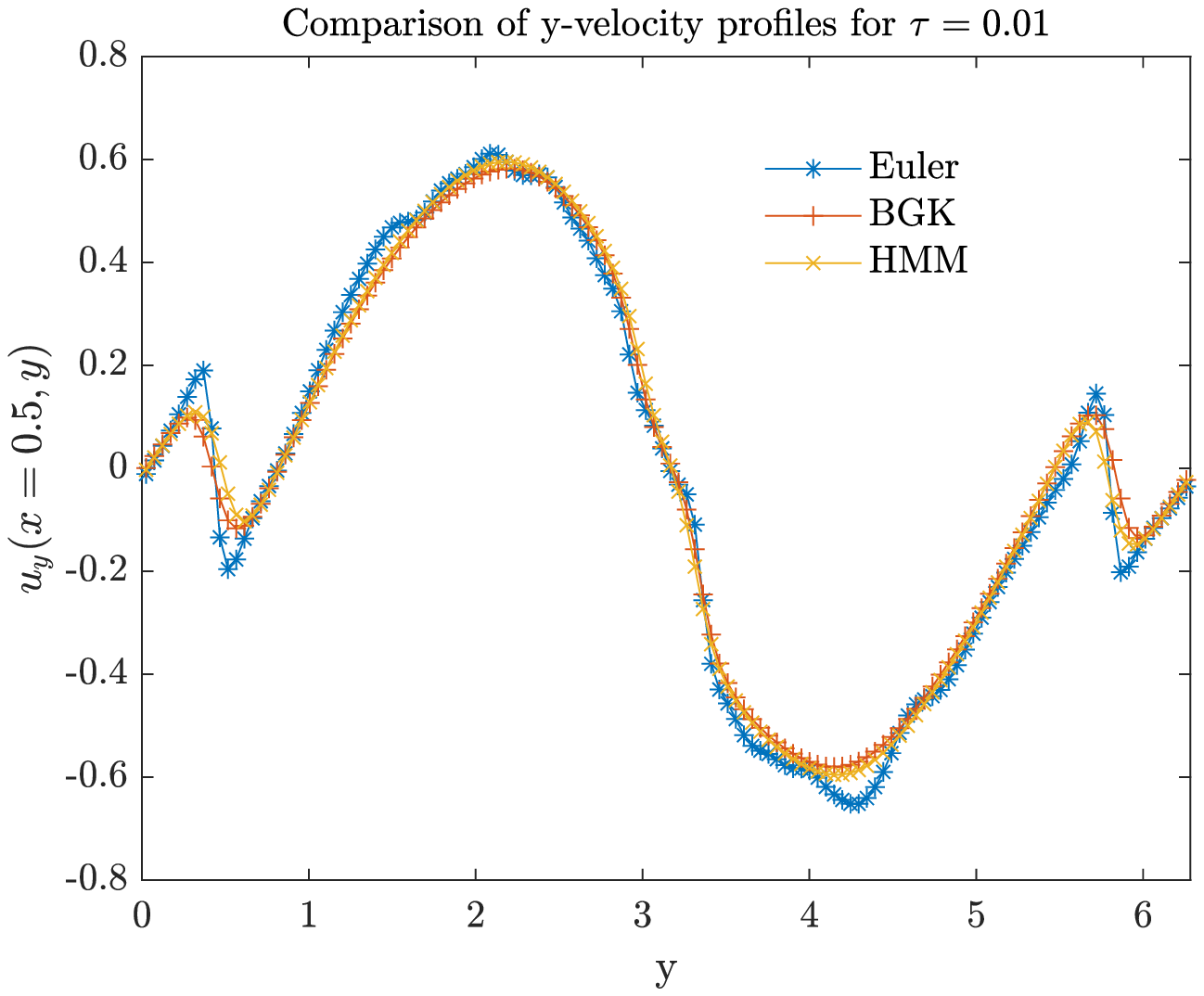}
		\includegraphics[width=0.3\textwidth]{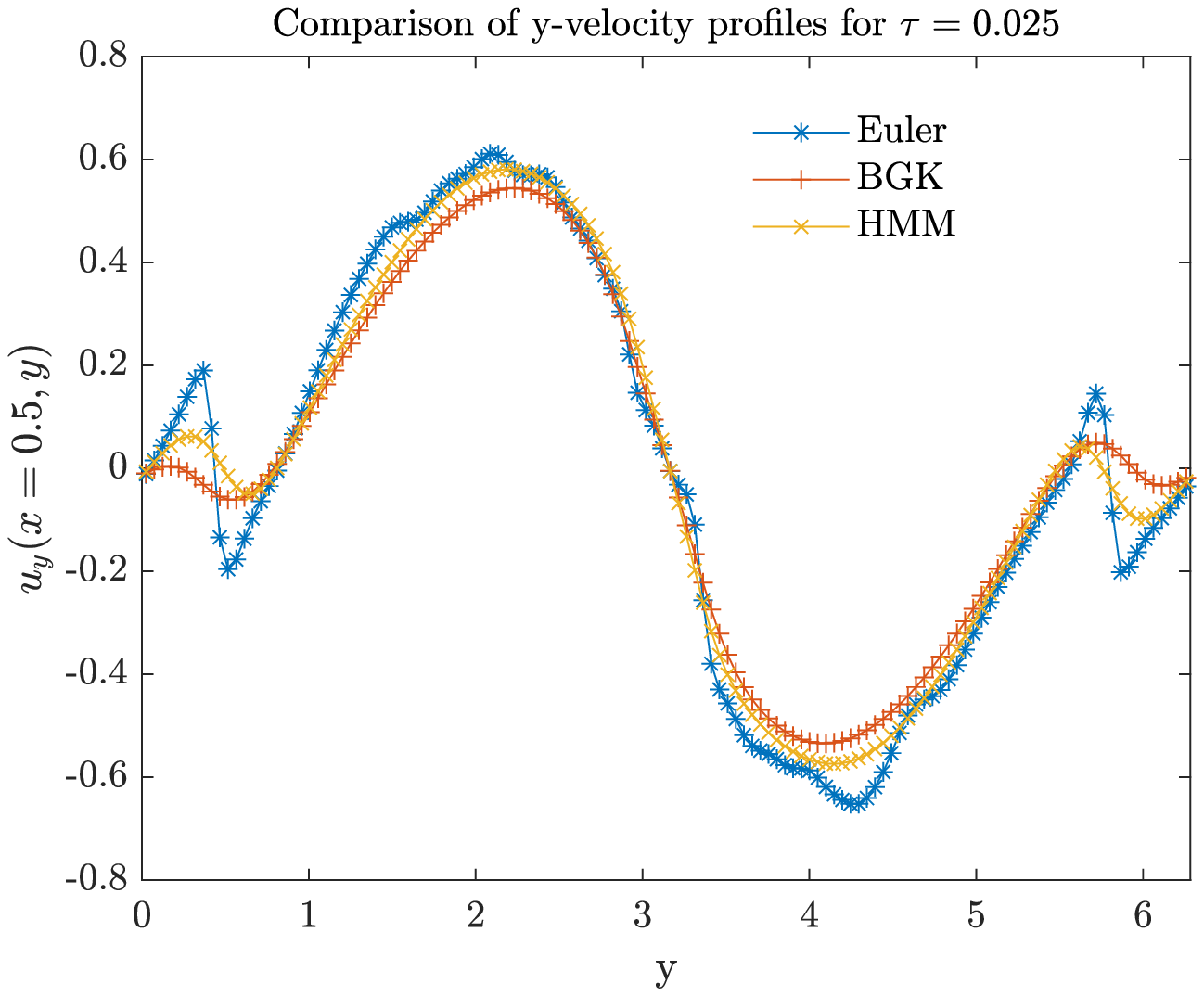}
	\caption{\textbf{Vortex flow}. Comparison of the density (left), $x$-velocity (middle) and $y$-velocity (right) profiles for $x=0.5$.
	The results for the Euler, the BGK and the HMM equations with respectively $\tau=0.005$, $\tau=0.01$ and $\tau=0.025$ are shown.	
	\label{fig:vortex2}}
\end{figure}

In this last section, we consider a Taylor-Green-type  problem. The initial data are given by
\begin{equation*}
	\label{initial_vortex}
	\rho(x,y,0)=1+0.1\cos(8\pi x)\sin(8\pi y),\ \bm{u}(x,y,0)=(u_x,u_y)=(\cos(x)\sin(y),-\cos(y)\sin(x)).
\end{equation*}
The temperature is $T=1$, the universal gas constant $R=1$, and the computational domain is $\Omega=[0,2\pi]^2$ discretized with
$N_x=N_y=128$ points in each direction. The velocity space is discretized using the same parameters as in the previous tests. We take a time
step $\Delta t=2 \ 10^{-4}$ and set a final time of $T_f=0.6$. The initial data are represented in Figure \ref{fig:vortex1} on the top
left, while the final solution for the isentropic Euler is reported in the same figure on the top right and illustrates the vortices
formation. The bottom image reports a three-dimensional view of the same solution for the density. In Figure \ref{fig:vortex2}, we
present comparisons between the isentropic Euler, the BGK and our HMM at the final time for the density and the velocities and for three
different values of the relaxation parameter $\tau=0.005$, $\tau=0.01$ and $\tau=0.025$ with $\tau_{HMM}=\tau/3$. Once again, the
obtained results indicate that the HMM is able to capture microscopic structures with a sufficiently high accuracy for the regimes under
considaration.

\section{Conclusions}
\label{sec5}
In this work, we have introduced a new hybrid multiscale model coupling the continuum and kinetic descriptions. We have applied the
Chapman-Enskog expansion on a suitably scaled stationary BGK equation and then upscaled the kinetic contribution over a microscopic box.
These microscopic regions are located at each points where the fluid equation are solved. The probability density function used to study
the kinetic evolution is related to the unknown density and velocity of the system through truncated Taylor expansions. Based on a scaling
assumption, related to the problem under consideration, a new multiscale continuum-kinetic model has been derived. A linear stability
analysis shows that new hybrid multiscale model is conditionally linearly stable in the strongly supersonic regime and unconditionally
stable in the subsonic and mildly supersonic one. Several numerical examples demonstrate that the hybrid multiscale model is more accurate
than standard fluid models and represents complex flows at different regimes more precisely.

In future work we would like to derive different hybrid multiscale models following the same general strategy outlined in this work:
couple a microscopic description with a macroscopic one through an upscaling of the microscopic information obtained by homogenization of
the micro quantities over micro sized boxes.

\begin{acknowledgments}
The work of A. Chertock was supported in part by NSF grants DMS-1818684. P. Degond holds a visiting professor association with the Department of Mathematics, Imperial College London, UK. The work of G. Dimarco was supported by the Italian Ministry of Instruction, University and Research (MIUR) under the PRIN Project 2017 (No. 2017KKJP4X). M.~Luk\'a\v{c}ov\'a-Medvid'ov\'a and A.~Ruhi were supported by the German Science Foundation (DFG) under the Collaborative Research Center TRR~146 Multiscale Simulation Methods for Soft Matter Systems (Project~C5). M.~Luk\'a\v{c}ov\'a-Medvid'ov\'a is grateful to the Gutenberg Research College and Mainz Institute of Multiscale Modelling for supporting her research.
\end{acknowledgments}

\bibliographystyle{siam}
\bibliography{cf_ref}

\end{document}